\documentclass[10pt]{article}
\usepackage{amsfonts}
\usepackage{amsmath}
\usepackage{amssymb}
\usepackage{graphicx}
\usepackage{epsfig}

\usepackage{latexsym,amssymb,amsfonts,amsmath, amscd}
\usepackage{epsfig,graphicx,graphics,latexsym,amssymb,amsfonts,amsmath, amscd}
\usepackage{changebar}
\usepackage{vmargin} %icii c la faute
\usepackage{graphicx}
\usepackage{makeidx}
\usepackage{euscript}
\date{\ }
\newtheorem{definition}{{\bf Definition}}[section]
\newtheorem{theorem}[definition]{{\bf Theorem}}

\newtheorem{corollary}[definition]{{\bf Corollary}}
\newtheorem{proposition}[definition]{\noindent {\bf Proposition}}
\newtheorem{lemma}[definition]{\noindent {\bf Lemma}}

\newtheorem{remark}[definition]{\noindent {\bf Remark}}

\def\endproof{\hfill {\kern 6pt\penalty 500
   \raise -0pt\hbox{\vrule \vbox to5pt {\hrule width 5pt
  \vfill\hrule}\vrule}}}
\newcommand{\bi}{\begin{itemize}}
\newcommand{\ei}{\end{itemize}}
\newcommand{\be}{\begin{enumerate}}
\newcommand{\ee}{\end{enumerate}}

\newcommand{\E}{{\mathcal E}}

\begin{document}
%j'ai commencé la rédaction le 15/avril/2011
%accepté pour publication le 19/juillet/2013

\title{ Graphs whose indecomposability graph is $2$-covered}
\author{Rim Ben Hamadou $^{(a)}$, Imed Boudabbous $^{(b)}$\\\\
(a) Facult\'e des Sciences de Sfax, Universit\'e de Sfax, Tunisie.\\E-mail:
rim.benhamadou@gmail.com \\ (b) Institut Pr\'eparatoire aux \'Etudes
d'Ing\'enieurs de Sfax, Universit\'e de Sfax, Tunisie. \\E-mail: imed.boudabbous@gmail.com
} \maketitle \vskip 0.2 cm

\noindent{\bf Abstract.} Given a graph $G=(V,E)$, a subset $X$ of
$V$ is an interval of $G$ provided that for any $a, b\in X$ and $
x\in V \setminus X$, $\{a,x\}\in E$ if and only if $\{b,x\}\in E$.
For example, $\emptyset$, $\{x\}(x\in V)$ and $V$ are intervals of
$G$, called trivial intervals. A graph whose  intervals  are trivial
is indecomposable; otherwise, it is decomposable. According to Ille, the indecomposability graph of an undirected indecomposable graph $G$ is the graph $\mathbb I(G)$ whose vertices are those of $G$ and edges are the unordered pairs of distinct
vertices $\{x,y\}$ such that the induced subgraph $G[V \setminus \{x,y\}]$ is indecomposable. We
characterize the indecomposable graphs $G$ whose $\mathbb I(G)$
admits a vertex cover of size $2$.

\noindent {\bf Keywords:} Graphs, indecomposable, interval,
indecomposability graph, partially critical.

\vspace{3mm} \setcounter{page}{1}

\section{Introduction and presentation of the results}
Over the years, the concept of indecomposability has become fundamental in the study of finite structures. Pioneered by T. Gallai in the theory of graphs with his seminal paper (\cite{Gallai}) and independently by R.Fra\"{\i}ss\'e (\cite{Fraiss}) in the
 theory of relations, this concept was developed in several papers e.g (\cite{EhR90, Fraisse, MH, Ille97, Kelly, SchTro93, Sumner}), and is now presented in a book by Ehrenfeucht, Harju and Rozenberg (\cite{EHR99}). Properties of the indecomposable substructures of a given indecomposable structures were developed by Schmerl and Trotter (1993) in their fundamental paper. Several papers along these lines have then appeared (\cite{BB, BBE, BouchBE, ChS, EHR99, Ille97, Pouzag09, Sayar, SchTro93}). Ille (1993) introduced the notion of indecomposability graph associated with a binary relation. This graph is an important tool of many research studies
 on the indecomposability. For example, Ille used that graph to show that indecomposable binary relations can be recognized, while there are not
  necessarily reconstructible.\\ %\cite{BBE, BoudIm11, BIlle09, BCI}.\\
This paper is about the indecomposability graph of undirected indecomposable graphs. The {\it indecomposability graph} of an indecomposable graph $G$ is the graph denoted by $\mathbb I(G)$, whose vertices are those of $G$ and the edges are the pairs
$\{x,y\}$ of distinct vertices such that $G[V \setminus \{x,y\}]$ is indecomposable. Given a graph $G$, with vertex set $V$ and edge set $E$, a {\it vertex cover} of $G$ is a subset $X$ of $V$ such that for each edge  $e
\in E(G)$, $e \cap X\neq \emptyset$.
%We also say that $G$ is {\it $X$- covered}. Consider a positive integer $k$.
We say that $G$ is {\it $k$-covered} if it is $X$-covered for a subset
$X$ of $V$ with $\mid X\mid =k$. \\
In this paper we give a description of indecomposability graphs $G$ such that $\mathbb I(G)$ is $2$-covered (a question raised in \cite{BCI}). A description of indecomposable tournaments $T$ such that $\mathbb I(T)$ is $2$-covered was given by the second author in (\cite{BoudIm11}).\\
A first reason to look at $2$-covered graphs is because for an indecomposable graph $G$ with $v(G) \geq 11$, $\mathbb{I}(G)$ is not $1$-covered. That is,\\
\begin{theorem} {\rm (\cite{Ille97})} \label{propnon1couvert}
Let $G = (V,E)$ be an indecomposable graph, with $v(G) \geq 11$. For
every $x \in V$, there are $y \neq z \in V \setminus \{x\}$ such
that $G[V \setminus \{y,z\}]$ is indecomposable.
\end{theorem}

The starting point of this result was due to Schmerl and Trotter proving that an indecomposable graph with $n \geq 6$ vertices contains an indecomposable induced subgraph on $n-2$ vertices.

 \begin{theorem}{\rm (\cite{SchTro93})}\label{thSchTro} Let $G=(V,E)$ be an indecomposable graph where $v(G)\geq 6$. Then,
there exist $x \neq y \in V$ such that $G[V \setminus \{x,y\}]$ is
indecomposable.
\end{theorem}

Theorem \ref{thSchTro} was improved as follows.
\begin{theorem} \label{thillemoins2}{\rm (\cite{Ille97})}
Let $G=(V,E)$ be an indecomposable graph, $X$ be a subset of $V$
such that $\mid X\mid \geq 4$ and $G[X]$ is indecomposable. If $\mid
V \setminus X\mid \geq 6$, then there exist $x \neq y \in V\setminus X$
such that $G[V \setminus \{x,y\}]$ is indecomposable.
\end{theorem}

%An important chapter of the theory of graphs concerns the indecomposability of graphs. For the studies along the lines, we can mention
%This paper is about non-oriented graphs. A {\it graph} (or a non-oriented graph) $G$ is defined by a finite vertex set $V(G)$ and an edge set $E(G)$, where an edge is an unordered pair of distinct vertices. Such a graph is denoted by $(V(G),E(G))$ or
%simply $(V, E)$. The {\it cardinality} of $G$ is that of $V(G)$ denoted by
%$v(G)$. With each subset $X$ of $V$, the subgraph of $G$ induced by $X$ is denoted by $G[X]$. We are interested in the study of indecomposable big-sized subgraphs in an indecomposable graph. The starting point was the result of J.H. Schmerl and W.T. Trotter, in $1993$ presented below. In fact, they prove that an indecomposable graph with $n \geq 6$ vertices contains an indecomposable subgraph of $n-2$ vertices.
Since, for each vertex $x$ of an indecomposable graph $G = (V,E)$, with $v(G) \geq 4$, there exists $X \subseteq V$ such that $x \in X$,
$\mid X \mid =4$ or $5$ and $G[X]$ is indecomposable \cite{CIlle98}, the above theorem follows.\\

A second reason for looking at $2$-covered indecomposability graphs of indecomposable graphs is that those graphs are arise in the study of $(-2)$-recognition, and more generally to progress toward the knowledge of the structure of the indecomposability graph of indecomposable graph.

%Hence the following problem follows.

%\begin{problem} (\cite{BCI}) Characterize the indecomposable
%graphs $G$ such that $\mathbb I(G)$ is $2$-covered.
%\end{problem}

%In this paper, we solve this problem which is a natural question in the study of the
%$(-2)$-recognition (\cite{BCI}). Notice that the characterization of indecomposable 'tournaments' $T$ such that $\mathbb I(T)$ is $2$-covered was obtained by I. Boudabbous \cite{BoudIm11}.\\

The major tool in our description of $2$-covered indecomposability graph is the notion of {\it minimal graph} defined as follows. Given two
distinct vertices $x$ and $y$ of an indecomposable graph $G$ of
cardinality $\geq 4$, we say that $G$ is {\it minimal} for
$\{x,y\}$, or {\it $\{x,y\}$-minimal}, whenever for each proper
subset  $X$ of $V(G)$, if $x, y \in X$ and $\mid \! X\! \mid \geq
3$, then $G[X]$ is decomposable. The minimal graphs for two vertices
were characterized by A. Cournier and P. Ille \cite{CIlle98}. In order to recall
this characterization, we introduce the following graphs $P_{n}$ and
$Q_{n}$, where $\mathbb{N}_{n}=\{1,\ldots,n\}$ for $n\geq 1$.\\

$\bullet$ For $n\geq 1$, the graph $P_{n}$ (see Figure $1$) is defined as follows. $V(P_{n})=\mathbb{N}_{n}$ and for $i \neq j \in \mathbb{N}_{n}, \,
\{i,j\} \in E(P_{n})$ if $\mid i-j \mid = 1$.\\

\begin{figure}[h!]
\centering
\includegraphics[width=8cm]{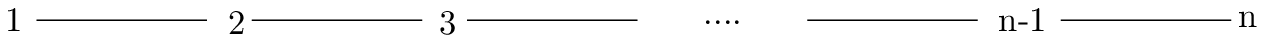}%\vspace{3.75 cm}
\caption{$P_{n}$}
\end{figure}

$\bullet$ For $n\geq 4$, the graph $Q_{n}$ (see Figure $2$) is defined as follows. $V(Q_{n})=\mathbb{N}_{n}$ and \\
$E(Q_{n})=E(P_{n-2}) \cup \{\{n-1,i\}: i \in \mathbb{N}_{n} \setminus \{n-1,n-2\}\}$.\\

\begin{figure}[h!]
\centering
\includegraphics[width=8cm]{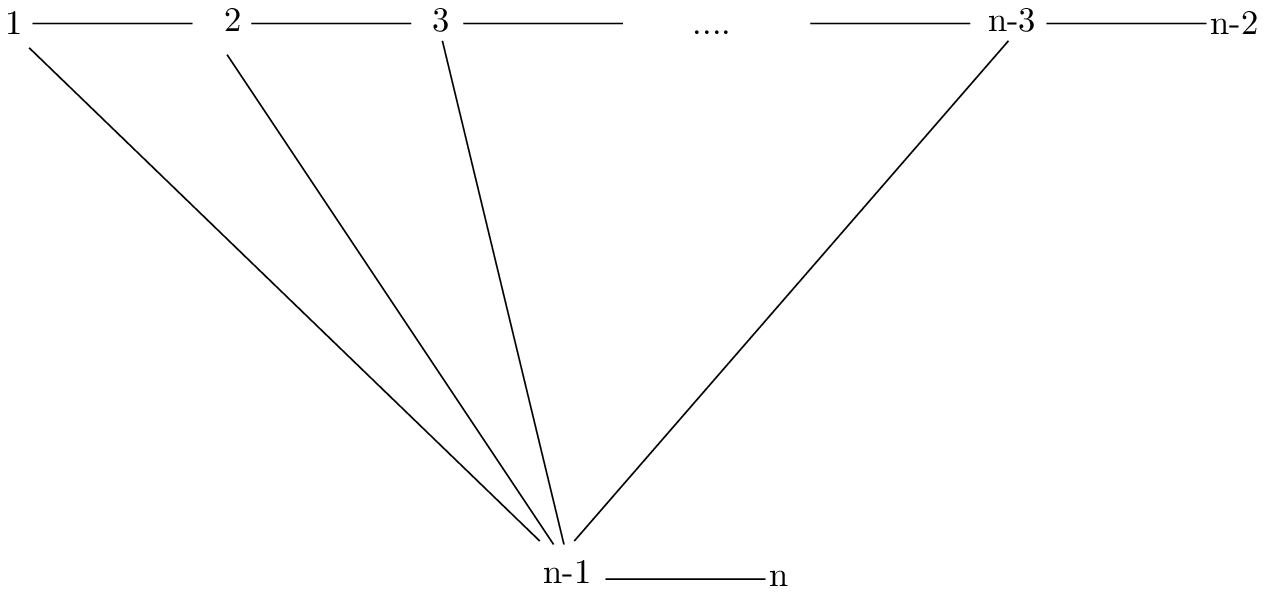}%\vspace{3.75 cm}
\caption{$Q_{n}$}
\end{figure}

For  $n\geq 4$, the graphs $P_{n}$, $Q_{n}$ and their complements are indecomposable
and $\{1,n\}$-minimal. Conversely :

\begin{theorem}{\rm (\cite{CIlle98})}
\label{THM} Given a graph $G$, with $v(G)\geq 4$, consider two
vertices $a \neq b \in V(G)$. The graph $G$ is $\{a,b\}$-minimal if and
only if there is an isomorphism $f$ from $G$ or $\overline{G}$ onto
$P_{v(G)}$ or $Q_{v(G)}$ such that $f(\{a,b\})= \{1,v(G)\}$.
\end{theorem}

We can easily check that:
\begin{remark}
$\bullet$ For $n\geq 6$, $\mathbb I (P_{n})= (\mathbb{N}_{n},\,  \{
\{1,2\}, \{1,n\}, \{n-1,n\} \})$.\\

$\bullet$ $ \mathbb I (Q_{n}) = \left \{
                               \begin{array}{ll}
                                 (\mathbb{N}_{n},\, \{\{1,2\}, \{2,n\}, \{n-1,n\} \}) & \hbox{$when \,\,\, n = 6,7$} \\
                                 (\mathbb{N}_{n},\,  \{ \{1,2\}, \{2,n\}, \{n-1,n\}, \{1,n\} \}) & \hbox{$when \,\,\,n\geq 8$.}
                               \end{array}
                             \right.
$
%For $n = 6,7$, $\mathbb I (Q_{n})= (\mathbb{N}_{n},\, \{
%\{1,2\}, \{2,n\}, \{n-1,n\} \})$.\\
%\hspace *{0.42cm} For $n\geq 8$, $\mathbb I (Q_{n})=
%(\mathbb{N}_{n},\,  \{ \{1,2\}, \{2,n\}, \{n-1,n\}, \{1,n\} \})$.\\
\end{remark}

Our results are presented below. Undefined terminology concerning graphs will be explained in Section 2.\\
Our first result is elementary:
\begin{proposition} \label{theofond}
Let $G$ be an indecomposable graph with $v(G)\geq 6$. Given  $a\neq b\in V(G)$, if
$\mathbb I (G)$ is $\{a,b\}$-covered, then $G$ contains an
$\{a,b\}$-minimal induced subgraph of cardinality $v(G)$, $v(G)-1$, $v(G)-3$
or $v(G)-5$.
\end{proposition}

The above Proposition leads to the description of the graphs $G$ whose
 $\mathbb I (G)$ is $\{a,b\}$-covered, from the $\{a,b\}$-minimal induced subgraphs embedding into $G$. We introduce the following classes of graphs.

\begin{itemize}
\item $\mathcal{P}$ is the set of $P_{n}$ for some  $n\geq 9$.

\item $\mathcal{Q}$ is the set of $Q_{n}$ for some  $n\geq 9$.

\item $\mathcal{P}_{-1}$ is the set of indecomposable graphs
 $G$ defined on  $\mathbb{N}_{n}$ for some  $n \geq 9$, such that  $\mathbb I (G)$ is  $\{1, n-1\}$-covered
 and   $G-n = P_{n-1}$.

 \item $\mathcal{Q}_{-1}$ is the set of indecomposable graphs
 $G$ defined on  $\mathbb{N}_{n}$ for some $n \geq 11$, such that  $\mathbb I (G)$ is  $\{1, n-1\}$-covered
 and   $G-n = Q_{n-1}$.

\item $\mathcal{P}_{-3}$ is the set of indecomposable graphs
 $G$ defined on $\mathbb{N}_{n}$ for some
$n \geq 12$, such that $\mathbb I (G)$ is  $\{1, n-3\}$-covered and
$G-\{n, n-1,n-2 \}= P_{n-3}$.

\item  $\mathcal{Q}_{-3}$ is the set of indecomposable graphs
 $G$ defined on $\mathbb{N}_{n}$ for some
$n \geq 10$, such that $\mathbb I (G)$ is  $\{1, n-3\}$-covered and
$G-\{n, n-1,n-2 \}= Q_{n-3}$.

\item $\mathcal{P}_{-5}$ is the set of indecomposable graphs
 $G$ defined on $\mathbb{N}_{n}$ for some
$n \geq 14$, such that $\mathbb I (G)$ is  $\{1, n-5\}$-covered and
$G-\{n, n-1,n-2, n-3, n-4 \}= P_{n-5}$.

\item  $\mathcal{Q}_{-5}$ is the set of indecomposable graphs
 $G$ defined on $\mathbb{N}_{n}$ for some
$n \geq 12$, such that $\mathbb I (G)$ is  $\{1, n-5\}$-covered and
$G-\{n, n-1,n-2, n-3, n-4 \}= Q_{n-5}$.

\end{itemize}

\begin{remark} \label{rqTHEOREME}
It is clear that $\mathcal{P}$ is a subset of $\mathcal{P}_{-1}$. Since $Q_{n} - 1 \simeq Q_{n-1}$ for $n \geq 5$, and the graph $Q_{n}$ is $\{1, n\}$-covered
 and $\{2, n\}$-covered, then each element of $\mathcal{Q}$ is isomorphic to an element of $\mathcal{Q}_{-1}$. Thus, we can say that,
  up to isomorphism, $\mathcal{Q} \subseteq \mathcal{Q}_{-1}$.
\end{remark}

Our description is done by the following result :
\begin{theorem}\label{THEOREME}
Given an indecomposable graph $G$ with $v(G) \geq 14$, $\mathbb I (G)$ is  $2$-covered if and only if $G$ or $\overline{G}$
is isomorphic to an element $G'$ of $\mathcal{P}_{-1} \cup
\mathcal{P}_{-3} \cup \mathcal{P}_{-5} \cup \mathcal{Q}_{-1} \cup
\mathcal{Q}_{-3} \cup \mathcal{Q}_{-5}$ with $v(G') \geq 14$.
\end{theorem}

\noindent{\em Proof\/.} First, assume that $G$ or $\overline{G}$ is isomorphic to an element $G'$ of $\mathcal{P}_{-1} \cup \mathcal{P}_{-3} \cup \mathcal{P}_{-5} \cup \mathcal{Q}_{-1} \cup
\mathcal{Q}_{-3} \cup \mathcal{Q}_{-5}$. As, by definition, $\mathbb{I}(G')$ is $2$-covered, then $\mathbb{I}(G)$ or $\mathbb{I}(\overline{G})$
is $2$-covered. In addition, $\mathbb{I}(G) = \mathbb{I}(\overline{G})$ because $G$ and $\overline{G}$ share the same intervals. It follows that
$\mathbb I(G)$ is  $2$-covered.\\

Conversely, consider a graph $G=(V,E)$ with $v(G) \geq 14$ and assume that $\mathbb I(G)$ is  $\{a,b\}$-covered where $a \neq b \in V$. Using
Proposition \ref{theofond}, there is a subset $X$ of $V$ containing $\{a,b\}$ such that $G[X]$ is $\{a,b\}$-minimal and $\mid X\mid = v(G)-k$ where
 $k \in \{0,1,3,5\}$. If $k=0$, then by Theorem \ref{THM}, we obtain that $G$ or $\overline{G}$ is isomorphic to an element of $\{P_{v(G)},Q_{v(G)}\}$.
  Thus, by Remark \ref{rqTHEOREME}, $G$ or $\overline{G}$ is isomorphic to an element of $\mathcal{P}_{-1} \cup \mathcal{Q}_{-1}$. If $k \in \{1,3,5\}$,
  we pose $V \setminus X = \{y_{i}, v(G)-k+1 \leq i \leq v(G)\}$. By Theorem \ref{THM}, there is an element $K$ of $\{G,\overline{G}\}$ and an isomorphism $f$ from $K[X]$ onto an element $H$ of $\{P_{v(G)-k},Q_{v(G)-k}\}$ with $f(\{a,b\})= \{1,v(G)-k\}$. Consider the graph $G'$ defined on $\mathbb{N}_{n}$ as follows.\\
   $G'[\mathbb{N}_{v(G)-k}]= H$, for $i \neq j \in \mathbb{N}_{n} \setminus \mathbb{N}_{v(G)-k}, \, \{i,j\} \in E(G')$ if and only if
    $\{y_{i},y_{j}\} \in E(G)$ and for $p \in  \mathbb{N}_{n} \setminus \mathbb{N}_{v(G)-k}$ and $q \in \mathbb{N}_{v(G)-k}$, $\{p,q\} \in E(G')$ if and only if $\{y_{p},f^{-1}(q)\} \in E(G)$. \\
    Clearly, by construction, the bijection $g$ from $V$ onto $\mathbb{N}_{n}$ defined by: $g _{/X} = f$ and $g(y_{i})=i$ for
    $y_{i} \in V \setminus X$ is an isomorphism from $K$ onto $G'$. As $K$ is indecomposable, $\mathbb I(K)$ is  $\{a,b\}$-covered and $g(\{a,b\})=\{1,v(G)-k\}$ then $G'$ is indecomposable and $\mathbb I(G')$ is $\{1,v(G)-k\}$-covered. Thus, $G$ or $\overline{G}$ is isomorphic to $G' \in \mathcal{P}_{-k} \cup \mathcal{Q}_{-k}$; which allows to conclude. \\
{\hspace*{\fill}$\Box$\medskip}

Section $4$ is devoted to describe each of the
classes  $\mathcal{P}_{-1}$, $\mathcal{P}_{-3}$, $\mathcal{P}_{-5}$,
$\mathcal{Q}_{-1}$, $\mathcal{Q}_{-3}$ and $\mathcal{Q}_{-5}$.\\

The following result is a direct consequence.
\begin{corollary}
Let $G=(V,E)$ be an indecomposable graph, with $v(G) \geq 14$. Then, for $a\neq b \in V$, one of the following
assertions is satisfied.
\begin{itemize}

\item There is $X \subset V$ such that $a, b\in X$, $\mid X\mid =
n-2$ and $G[X]$ is indecomposable.

\item There is an isomorphism $f$ from $G$ or $\overline{G}$ onto an element of $\mathcal{P}\cup \mathcal{Q}$ such that $f(\{a,b\})=\{1,n\}$.

\item There is an isomorphism $f$ from $G$ or $\overline{G}$ onto an element of
$\mathcal{P}_{-k}\cup \mathcal{Q}_{-k}$ such that
$f(\{a,b\})=\{1,n-k\}$, where $k \in \{1,3,5\}$.
\end{itemize}
\end{corollary}

This paper is organized as follows. Section $2$ contains the material needed about indecomposable graphs. Section $3$ contains the characterization of critical and partially critical graphs.  The description of the
classes  $\mathcal{P}_{-1}$, $\mathcal{P}_{-3}$, $\mathcal{P}_{-5}$,
$\mathcal{Q}_{-1}$, $\mathcal{Q}_{-3}$ and $\mathcal{Q}_{-5}$ is obtained in Section $4$.

\section{Prerequisites}
\subsection{Basic terminology and notation for graphs}
We denote by $\mathcal{P}_{2}(V)$, the set of pairs of distinct elements of a set $V$. Given a graph $G=(V, E)$, the {\it complement} of $G$ is the
graph $\overline{G} = (V,\mathcal{P}_{2}(V)\setminus E)$. The
\emph{neighbourhood of $x$ in $G$}, denoted by $N_{G}(x)$ or simply
$N(x)$, is the set $N_{G}(x) = \{y \in V \setminus \{x\}: \{x,y\}
\in E\}$. The graph $G$ is \emph{complete} (resp. \emph{empty}) if,
$E = \mathcal{P}_{2}(V)$ (resp. $E = \emptyset$).\\
Let $G=(V,E)$ be a graph. For every two vertices $x$, $y$ of $V$, the
notation $x\rule[0.4mm]{0.6cm}{0.1mm}\,\, y$ signifies that $\{x,y\}
\in E$, and $x \ldots y$ (or $x \quad \,\,  y$) signifies that
$\{x,y\} \notin E$. For each two disjoint subsets $I$ and $J$ of
$V$, we denote by $I \rule[0.4mm]{0.6cm}{0.1mm}\,\, J$ whenever for
each $(x, y) \in I \times J$, $x \rule[0.4mm]{0.6cm}{0.1mm}\,\,
y$. Similarly, for each $x \in V$ and for each $Y \subseteq V
\setminus \{x\}$, $x \rule[0.4mm]{0.6cm}{0.1mm}\,\, Y$ (resp. $x
.\ldots Y$) signifies that $x \rule[0.4mm]{0.6cm}{0.1mm}\,\, y$
(resp. $x .\ldots y$) for each $y \in Y$. Furthermore, $x \sim Y$
means $x \rule[0.4mm]{0.6cm}{0.1mm}\,\, Y$ or $x .\ldots Y$.
The negation is denoted by $x \not\sim Y$.
%Let $G=(V,E)$ be a graph. For each $x \neq y$ and $x' \neq y' \in
%V$, we note $\{x,y\}_{G} \sim \{x',y'\}_{G}$ (or simply $\{x,y\}
%\sim \{x',y'\}$) if $\{x,y\} \in E$ if and only if $\{x',y'\} \in
%E$. Otherwise, we note $\{x,y\} \not\sim \{x',y'\}$.

The notions of isomorphism, subgraph and embedding are defined in
the following way. Let $G = (V, E)$ and $G' = (V', E')$ be
two graphs. First, a one-to-one correspondence $f$ from $V$ onto $V'$ is an
{\it isomorphism} from $G$ onto $G'$ provided that for $x, y \in V$,
$\{x, y\} \in E$ if and only if $\{f(x), f(y)\} \in E'$. The graphs
$G$ and $G'$ are said to be {\it isomorphic}, which is denoted
by $G \simeq G'$, if there is an isomorphism from $G$ onto $G'$. Second, for $X \subseteq V$, the graph $G[X] := (X, E
\cap \mathcal{P}_{2}(X))$ is an {\it induced subgraph} of $G$. For $X \subseteq V$
(resp. $x \in V $), the induced subgraph $G[V \setminus X]$ (resp. $G[ V
\setminus \{x\}]$) is denoted by $G-X$ (resp. $G-x$). If $G'$ is isomorphic to an induced subgraph of $G$, we say that
$G'$ {\it embeds} into $G$. Let $p$ be a partition of $V$; the graph $G$ is {\it
multipartite} w.r.t $p$ if for every $M \in p, \, G[M]$ is empty. It is {\it bipartite} when $\mid p\mid = 2$.\\ A nonempty subset $C$ of $V$ is \emph {a connected component of $G$} if for $x \in C$
and $y \in V \setminus C, \, \{x,y\} \notin E$ and if for $x \neq y
\in C,$ there is a sequence $x=x_{0}, \, \ldots, \, x_{n}=y$ of
elements of $C$ such that $\{x_{i},x_{i+1}\} \in E$ for $0\leq i
\leq n-1$. A vertex $x$ of $G$ is \emph{isolated} if $\{x\}$
constitutes a connected component of $G$. The graph $G$ is
\emph{connected} if it has at most one connected component of $G$. Otherwise, it is called \emph{non-connected}.

\subsection{Indecomposable graphs}
Given a graph $G = (V,E)$, a subset $I$ of $V$ is an \emph{interval} \cite{CIlle98, Fraisse, SchTro93} (clan \cite{EhR90}, module \cite{Spin92}) of $G$ provided that for every $x \in V \setminus I$, $x \sim I$. In other words, $I$ is an interval if every vertex outside $I$ has the same behavior to all elements of $I$. Clearly, $\emptyset $, $V$ and $\{x\}$, where $x \in V$,
are intervals of $G$, called {\it trivial} intervals. A
graph is then said to be {\it indecomposable} \cite{Ille97, SchTro93}
if all of its intervals are trivial. It is said to be {\it decomposable}
otherwise. Notice that the graphs $G$ and $\overline{G}$ share the same intervals. Thus, $G$ is indecomposable if and only if $\overline{G}$ is indecomposable. For example, all graphs of cardinality $3$ are decomposable and up to isomorphism, the graph $P_{4}$ is the unique indecomposable graph of cardinality $4$.\\

We review relevant properties of indecomposable graphs.\\
Given a graph $G = (V,E)$, consider a subset $X$ of $V$ such that $\mid X
\mid \geq 4$ and  $G[X]$ is indecomposable.  We use the following
subsets of $V\setminus X$.
\begin{itemize}
\item $\rm{Ext}(X)$ is the set of  $v \in V\setminus X$ such that $G[X \cup
\{v\}]$ is indecomposable;

\item $\langle X\rangle $ is the set of  $v \in V\setminus X$ such that $v \sim X$;

\item For each  $u \in X$, $X(u)$ is the set of  $v \in V\setminus X$ such
that  $ \{u, v\}$ is an interval of $G[X \cup \{v\}]$.
\end{itemize}

The family constituted by $\rm{Ext}(X)$, $\langle X\rangle$ and
$X(u)$, where $u\in X$, is denoted by $p_{X}$.\\

Besides, the family $p_{X}$ is divided as follows.
\begin{itemize}
  \item $X^{-}$ is the set of elements $v$ of $V \setminus X$
such that $v .\ldots X$.
  \item $X^{+}$ is the set of elements $v$ of $V \setminus X$
such that $v \rule[0.4mm]{0.6cm}{0.1mm}\,\,X$.
  \item $X^{-}(u)$ is the set of elements $v$ of $X(u)$
such that $\{u,v\} \notin E$.
 \item $X^{+}(u)$ is the set of elements $v$ of $X(u)$
such that $\{u,v\} \in E$.
\end{itemize}

We then introduce the three families below :
\begin{itemize}
  \item $q_{X} = \{Ext(X),X^{-},X^{+}\} \cup \{X^{-}(u),X^{+}(u)\}_{u \in X}$
  \item $q^{-}_{X}= \{X^{-}\} \cup \{X^{-}(u); u \in X\}$
  \item $q^{+}_{X}= \{X^{+}\} \cup \{X^{+}(u); u \in X\}$\\
\end{itemize}

\begin{theorem} \label{lemEhR}{\rm (\cite{EHR99})}
Let  $G = (V,E)$ be a graph, consider a subset $X$ of $V$ such that
$\mid X \mid \geq 4$ and $G[X]$ is indecomposable. The family
$p_{X}$ realizes a partition of $V\setminus X$. Moreover, the
following hold.

\begin{enumerate}
\item Given   $u \in X$, $v \in X(u)$ and  $w \in V\setminus(X \cup
X(u))$. If $G[X \cup \{v, w\}]$ is  decomposable, then  $\{u, v\}$
is an interval of  $G[X \cup \{v, w\}]$.

\item Let  $v \in \langle X\rangle$ and   $w \in V\setminus(X \cup \langle
X\rangle)$. If $G[X \cup \{v, w\}]$ is  decomposable, then  $X \cup
\{w\}$ is an interval of  $G[X \cup \{v, w\}]$.

\item Let  $v \neq w \in \rm{Ext}(X)$. If $G[X \cup \{v, w\}]$ is
decomposable, then  $\{v, w\}$ is an interval of  $G[X \cup \{v,
w\}]$.

\end{enumerate}
\end{theorem}

As a consequence of the above theorem, we obtain the following.
\begin{corollary} \label{corolEhR}{\rm (\cite{EhR90})} Let $G = (V,E)$ be an indecomposable graph. If $X$ is a subset of
$V$ such that  $\mid X \mid \geq 4$, $\mid V\setminus X \mid \geq 2$
and  $G[X]$ is  indecomposable, then there are two distinct elements $x$ and $y$ of $V\setminus X$ such that $G[X\cup\{x, y\}]$ is
indecomposable.
\end{corollary}

Given Corollary \ref{corolEhR}, we introduce the following graph.
Let $G=(V,E)$ be a graph, $X$ be a subset of $V$ such that $\mid X\mid
\geq 4$, $\mid V \setminus X\mid \geq 2$ and $G[X]$ is
indecomposable. The graph $G_{X}=(V \setminus X,E_{X})$ is defined
as follows. For each $x \neq y \in V \setminus X ,\,\, \{x,y\} \in
E_{X}$ if $G[X \cup \{x,y\}]$ is indecomposable.\\

We make the following remark:
\begin{remark}
Given a graph  $G=(V,E)$, let $X$ be a subset of $V$ such
that $\mid X\mid \geq 4$, $\mid V \setminus X\mid \geq 2$ and $G[X]$ is indecomposable. Consider distinct elements $x$ and $y$ of $V \setminus X$.
If $x,y \in \langle X\rangle$, then $X$ is an interval of $G[X \cup \{x,y\}]$. If $x,y \in X(u)$ where $u \in X$, then $\{u,x,y\}$ is an interval of
$G[X \cup \{x,y\}]$. Consequently, for each $M \in p_{X} \setminus \{Ext(X)\}, \, G_{X}(M)$ is empty. In other words, if $Ext(X) = \emptyset$, then
$G_{X}$ is multipartite by $p_{X}$. Moreover, if $x \in X^{-}$ (resp. $x \in X^{+}$) and $y \notin \langle X\rangle$, then $\{x,y\} \in E$ if and only
if $\{x,y\} \in E_{X}$ (resp. $\{x,y\} \notin E_{X}$). Finally, assume that $x \in X(u), y \in X(v)$ where $u \neq v \in X$ such that
$\{u,v\} \notin E$ (resp. $\{u,v\} \in E$). We also obtain that $\{x,y\} \in E$ if and only
if $\{x,y\} \in E_{X}$ (resp. $\{x,y\} \notin E_{X}$).
\end{remark}

\section{Critical and partially critical graphs}
In this section, we recall the characterization of critical and partially critical graphs which is used in our proof.\\
To begin with, we have to introduce the following definitions. Consider
an indecomposable graph $G = (V,E)$ with $v(G) \geq 2$. A vertex $x$
of $V$ is called a \emph{critical vertex} of $G$ if $G-x$ is
decomposable. The graph $G$ is \emph {critical} if all its vertices are critical. For example, for each integer $n \geq 2$, the graph $G_{2n}$ shown in Figure $3$ and defined below is critical. The vertex set of $G_{2n}$ is $\{0,
\ldots, 2n-1\}$ and for $i \neq j \in \{0, \ldots, 2n-1\},
\, \{i,j\}$ is an edge of $G_{2n}$ if there exist $k \leq l \in \{0,
\ldots, n-1\}$ such that $\{i,j\}=\{2k,2l+1\}$.
\begin{figure}[h!]
\centering
\includegraphics[width=2.5cm]{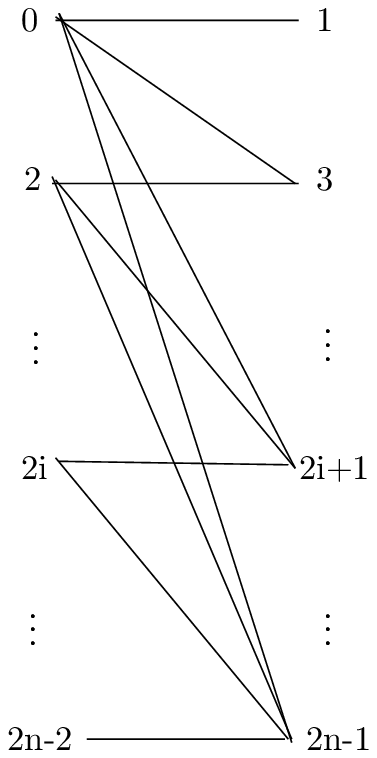}%\vspace{3.75 cm}
\caption{$G_{2n}$}
\end{figure}

\begin{theorem}{\rm (\cite{SchTro93})}
Let $G=(V,E)$ be an indecomposable graph. $G$ is critical if and
only if $G$ is isomorphic to $G_{2n}$ or $\overline{G_{2n}}$ where
$n\geq 2$.
\end{theorem}

Given an indecomposable graph $G = (V,E)$, consider a proper subset
$X$ of $V$ such that \\$\mid \! X\! \mid \geq 4$ and $G[X]$ is
indecomposable. The graph $G$ is \emph {critical according to
$G[X]$} if each element of $V \setminus X$ is critical. A graph is \emph {partially critical} if it is
critical according to one of its indecomposable induced subgraphs containing
at least $4$ vertices.\\

The partially critical graphs are characterized by the following two results.
\begin{theorem} \label{ThBDInonconnexe}{\rm (\cite{BDI08})}
Given a graph  $G=(V,E)$, let $X$ be a proper subset of $V$ such
that $\mid X\mid \geq 4$ and $G[X]$ is indecomposable. The graph $G$
is indecomposable and critical according to $G[X]$ if and only if
the three assertions below are satisfied.
\begin{itemize}
  \item $H1$: The partitions $p_{X}$ and $q_{X}$ coincide.
  \item $H2$: For each $M \in q^{-}_{X}$, $G[M]$ is
empty, and for each $N \in q^{+}_{X}$, $G[N]$ is complete.
  \item $H3$: For each connected component $C$ of $G_{X}$, $G[X \cup
C]$ is indecomposable and critical according to $G[X]$.
\end{itemize}
\end{theorem}

\begin{theorem} \label{ThBDIconnexe}{\rm (\cite{BDI08})}
Given a graph  $G=(V,E)$, let $X$ be a proper subset of $V$ such
that $\mid X\mid \geq 4$, $\mid V \setminus X\mid \geq 3$, $G[X]$ is
indecomposable and $G_{X}$ is connected. The graph $G$ is
indecomposable and critical according to $G[X]$ if and only if the
three assertions below are satisfied.
\begin{itemize}
  \item $K1$: $Ext(X) = \emptyset$.
  \item $K2$: The partitions $p_{X}$ and $q_{X}$ coincide.
  \item $K3$: For each $M \in q^{-}_{X}$, $G[M]$ is
empty, and for each $N \in q^{+}_{X}$, $G[N]$ is complete.
  \item $K4$: The graph $G_{X}$ is critical and bipartite by
$p_{X}$.\\
\end{itemize}
\end{theorem}

The following corollary is an immediate consequence of Theorems \ref{ThBDInonconnexe} and \ref{ThBDIconnexe}.
\begin{corollary} \label{corolBDI} {\rm (\cite{BDI08})}
If a graph $G$ is critical according to some induced subgraph $G[X]$, then $G_{X}$ has no isolated vertices.
\end{corollary}

\section{Description of the classes  $\mathcal{P}_{-1}$, $\mathcal{P}_{-3}$, $\mathcal{P}_{-5}$,
$\mathcal{Q}_{-1}$, $\mathcal{Q}_{-3}$ and $\mathcal{Q}_{-5}$}
First, notice the following. If $G$ is an $\{a,b\}$-minimal graph,
with $a \neq b$ and $v(G)\geq 6$, then $\mathbb I (G)$ is
$\{a,b\}$-covered. Indeed, if $\mathbb I (G)$ is
 not $\{a,b\}$-covered, then there exists $\{c,d\}\in E(\mathbb I (G))$
such that $\{c,d\}\cap\{a,b\}=\emptyset$. Thus, $G-\{c,d\}$ is
 indecomposable, which contradicts the minimality of $G$. \\

\subsection{Proof of Proposition \ref{theofond}}
Consider a graph $G=(V,E)$ with $v(G) \geq 6$ and assume that $\mathbb I (G)$ is $\{a,b\}$-covered where $a \neq b \in V$. Consider a minimal subset $X$ of $V$ under inclusion among the subsets $Y$ of $V$ satisfying $\mid Y \mid \geq 4$,
$\{a, b\} \subseteq Y$, and  $G[Y]$ is indecomposable. By minimality of $X$, $G[X]$
is $\{a,b\}$-minimal. Clearly, from Theorem \ref{thillemoins2},
$\mid V\setminus X\mid \, < 6$. It remains to verify that $\mid
V\setminus X\mid = 0,1,3$ or $5$. As $\mathbb I (G)$ is
$\{a,b\}$-covered, $\mid V\setminus X\mid \neq 2$. Moreover, Corollary \ref{corolEhR}
implies that $\mid V\setminus X\mid \neq 4$.\\
{\hspace*{\fill}$\Box$\medskip}

Now, we describe each of the classes  $\mathcal{P}_{-1}$, $\mathcal{P}_{-3}$, $\mathcal{P}_{-5}$,
$\mathcal{Q}_{-1}$, $\mathcal{Q}_{-3}$ and $\mathcal{Q}_{-5}$.

\subsection{The class  $\mathcal{P}_{-1}$}
The next proposition describes the class $\mathcal{P}_{-1}$.
\begin{proposition} \label{proppmoins1} Given a graph $G$ defined on $\mathbb{N}_{n}$,
where $n \geq 9$, $G \in \mathcal{P}_{-1}$\,  if and only if \,\,$G-n = P_{n-1}$ and either $N_{G}(n) = \{k\} \,where \,\,k \in \{3, \ldots,n-3\} \cup \{1,n-1\}$ or $N_{G}(n) = \{k,k+1\} \,where \,\,k \in \{2, \ldots,n-3\}$.

\end{proposition}

\begin{figure}[h!]
\centering
\includegraphics[width=9cm]{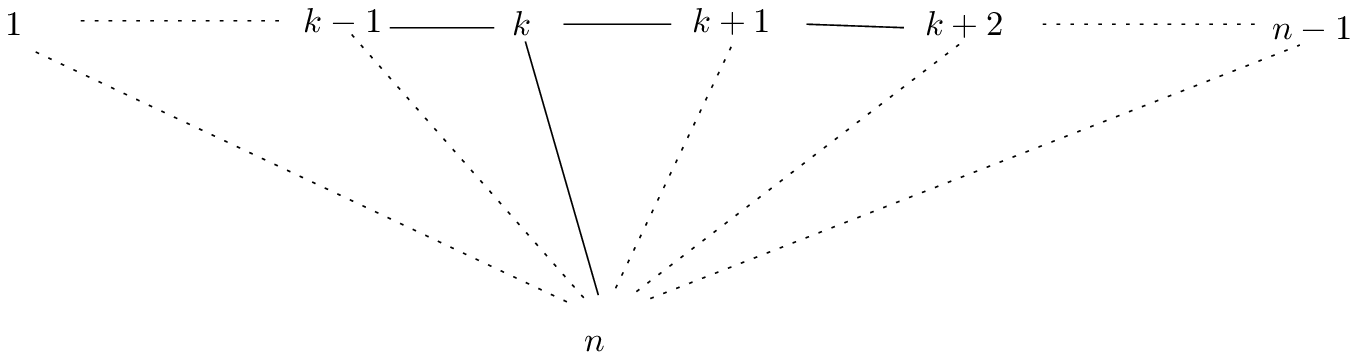}%\vspace{3.75 cm}
\caption{$N_{G}(n) = \{k\}$}
\end{figure}

\begin{figure}[h!]
\centering
\includegraphics[width=9cm]{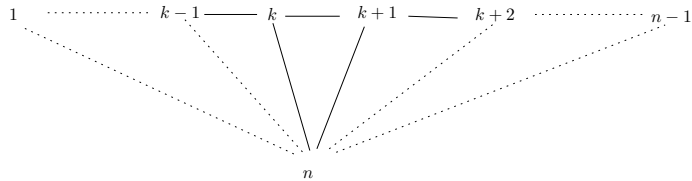}%\vspace{3.75 cm}
\caption{$N_{G}(n) = \{k, k+1\}$}
\end{figure}

\noindent{\em Proof\/.} Let $G$ be a graph defined on
$\mathbb{N}_{n}$, where $n \geq 9$ such that $G-n = P_{n-1}$. \\
First, assume that $N_{G}(n)=\{k\}$ where $k \in \{3, \ldots,n-3\} \cup \{1,n-1\}$. To start with, observe that if $N_{G}(n)=\{1\}$ (resp. $N_{G}(n)=\{n-1\}$), then $G
\simeq P_{n}$ and in particular, $G$ is indecomposable. Moreover, $\mathbb I(G)=(\mathbb{N}_{n}, \{ \{1,n\},
\{n-2,n-1\}, \{n,n-1\}\})$ (resp. $\mathbb I(G)=(\mathbb{N}_{n}, \{
\{1,2\}, \{n-1,n\}, \{1,n\}\})$). Hence, $\mathbb I(G)$ is
$\{1,n-1\}-$covered. \\
Assume that $k \in \{3, \ldots,n-3\}$. We show that $G$ is indecomposable. Since $G[\mathbb{N}_{n-1}] = P_{n-1}$ is indecomposable, we use the partition $p_{(\mathbb{N}_{n-1})}$ as follows.
We have $n \notin \langle \mathbb{N}_{n-1}\rangle$ because $n-1\, .\ldots\, n
\rule[0.4mm]{0.6cm}{0.1mm}\,\,k$. Moreover, $n \notin \mathbb{N}_{n-1}(n-1)$ because $n-1\, .\ldots\, k
\rule[0.4mm]{0.6cm}{0.1mm}\,\,n,$ and $n \notin \mathbb{N}_{n-1}(1)$ because $1\, .\ldots\, k
\rule[0.4mm]{0.6cm}{0.1mm}\,\,n$. Furthermore, $n \notin \mathbb{N}_{n-1}(i)$ for $2 \leq i \leq k$
because $n\, .\ldots\, i-1 \rule[0.4mm]{0.6cm}{0.1mm}\,\,i$. Similarly, as $n\, .\ldots\, i+1 \rule[0.4mm]{0.6cm}{0.1mm}\,\,i$, $n \notin \mathbb{N}_{n-1}(i)$ for
$k+1 \leq i \leq n-2$. Thus,$\, n \notin \bigcup\limits _{i=2}^{n-2} \mathbb{N}_{n-1}(i).$\\
Therefore, $n \notin \langle \mathbb{N}_{n-1}\rangle$ and $n \notin \mathbb{N}_{n-1}(i)$ for each $i \in \mathbb{N}_{n-1}$. Since $p_{(\mathbb{N}_{n-1})}$ is a
partition by Theorem \ref{lemEhR}, $n \in Ext(\mathbb{N}_{n-1})$ or equivalently, $G[\mathbb{N}_{n-1} \cup \{n\}]=G$ is indecomposable.\\
Now, we prove that $\mathbb I(G)$ is $\{1,n-1\}-$covered. Given $i<j \in \mathbb{N}_{n} \setminus \{1,n-1\}$, we have to verify that $G-\{i,j\}$ is
decomposable. If $i \geq k+1$, then $\mathbb{N}_{i-1} \cup \{n\} \setminus \, \{j\}$ is
a non-trivial interval of $G-\{i,j\}$. If $i \leq k$, then $G-\{i,j\}$ is decomposed into $(\{i+1, \ldots,n\}
\setminus \, \{j\}) \, \ldots \,  \{1, \ldots ,i-1 \}$. Consequently, $\{i+1, \ldots,n\} \setminus \, \{j\}$ is a non-trivial interval of $G-\{i,j\}$.\\

\noindent Second, assume that $N_{G}(n)=\{k,k+1\}$ where $k \in \{2, \ldots,n-3\}$. As $G[\mathbb{N}_{n-1}] = P_{n-1}$ and
$n-1 \geq 4$, then $G[\mathbb{N}_{n-1}]$ is indecomposable. We show that $n \in Ext(\mathbb{N}_{n-1})$. In fact, $n \notin
\langle\mathbb{N}_{n-1}\rangle$ because $n-1\, .\ldots\, n
\rule[0.4mm]{0.6cm}{0.1mm}\,\,k$. Besides, $n \notin
\mathbb{N}_{n-1}(1)$ because $1\, .\ldots\, k+1
\rule[0.4mm]{0.6cm}{0.1mm}\,\,n$. Moreover, $n \notin
\mathbb{N}_{n-1}(n-1)$ because $n-1 .\ldots k
\rule[0.4mm]{0.6cm}{0.1mm}\,\,n$. We also have for $2 \leq i \leq k,
\,n \notin \mathbb{N}_{n-1}(i)$ because $n\, .\ldots\, i-1
\rule[0.4mm]{0.6cm}{0.1mm}\,\,i $ and for $k+1 \leq i \leq n-2, \,
\,n \notin \mathbb{N}_{n-1}(i)$ because $n\, .\ldots\, i+1
\rule[0.4mm]{0.6cm}{0.1mm}\,\,i$. Thus, $\, n \notin \bigcup \limits_{i=1}^{n-1} \mathbb{N}_{n-1}(i).$\\
It follows by Theorem \ref{lemEhR}, that $G$ is indecomposable.\\
To prove that $\mathbb I(G)$ is $\{1,n-1\}-$covered, we proceed as previously.\\\\

Conversely, consider $G \in \mathcal{P}_{-1}$. We distinguish the following cases.
\begin{itemize}
\item \underline{Case 1: } $1 \rule[0.4mm]{0.6cm}{0.1mm}\,\,n$. As $\mathbb{N}_{n-1}$ is not an interval of $G$, there is $i \in \{2, \ldots,n-1\}$ such that $i .\ldots
n$. Set $\gamma=min (\{i \in \{2, \ldots,n-1\}\, : \, n \ldots i\})$. Clearly, $\gamma\, .\ldots\, n
\rule[0.4mm]{0.6cm}{0.1mm}\,\, \{1, \ldots,\gamma-1\}$. Set $Z=\{4, .\ldots , n-1\}$. Since $G[Z] \simeq P_{n-4}$ and $n \geq 9$, then $G[Z]$ is indecomposable. As $Z .\ldots 1, \, 1 \in Z^{-}$. Since $G-\{2,3\}$
is decomposable, it follows from Theorem \ref{lemEhR} that $n \in \langle Z\rangle$. $n \in Z^{+}$. We have $n .\ldots 2$ or $n .\ldots 3$ because $G$ is indecomposable; which implies $n \notin \langle \mathbb{N}_{n-4}\rangle $. Moreover, as $n-1 \in \mathbb{N}_{n-4}^{-}$ and $n \rule[0.4mm]{0.6cm}{0.1mm}\,\,n-1$, then by $2$ of
Theorem \ref{lemEhR}, $G-\{n-2,n-3\}$ is indecomposable, that contradicts $\mathbb I(G)$ is $\{1,n-1\}-$covered. Consequently, $n \in Z^{-}$ and so, $2\leq \gamma \leq 4$. If $\gamma=2$, we have $3 .\ldots n$ because $\{2,n\}$ is not an interval of $G$ and so $N_{G}(n)=\{1\}$. If $\gamma=3$ or $\gamma=4$, then $\{\gamma-2,n\}$ would be a non-trivial interval of $G$, that is impossible.\\
\item \underline{Case 2: } $1\, .\ldots\, n$. Since $\mathbb{N}_{n-1}$ is not an interval of $G$, there is $i \in \{2, \ldots,n-1\}$ such that $i \rule[0.4mm]{0.6cm}{0.1mm}\,\,n$. Set $\mu=min (\{i \in \{2, \ldots,n-1\}\, : \, n \rule[0.4mm]{0.6cm}{0.1mm}\,\, i\})$. Clearly, $\{1, \ldots,\mu -1\}\, .\ldots\, n \rule[0.4mm]{0.6cm}{0.1mm}\,\, \mu$. If $\mu=n-1$, then $N_{G}(n)=\{n-1\}$. Furthermore, $\{n-1,n\}$ would be a non-trivial interval of $G$ if $\mu=n-2$, that is impossible.
If $\mu=n-3$, then $n .\dots n-1$ because $\{n-2,n\}$ is not an interval of $G$, so we obtain for $k = n-3$, that either $N_{G}(n) =  \{k\} \, when \,\, n-2 .\ldots\, n $ or $N_{G}(n)=\{k, k+1\} \, when \,\,n-2 \rule[0.4mm]{0.6cm}{0.1mm}\,\, n$.\\

Thus, we may assume that $2 \leq \mu \leq n-4$. Observe that $n .\ldots n-1$. Otherwise, $n \rule[0.4mm]{0.6cm}{0.1mm}\,\, n-1, \, n-1 \in \mathbb{N}_{n-4}^{-}$ and $1\, .\ldots\, n \rule[0.4mm]{0.6cm}{0.1mm}\,\,\mu$, so that $n \notin \langle \mathbb{N}_{n-4} \rangle$. It would follow from $2$ of
Theorem \ref{lemEhR} that $G[\mathbb{N}_{n-4} \cup \{n,n-1\}]=G-\{n-3,n-2\}$ is indecomposable and $\mathbb I(G)$ would not be $\{1,n-1\}-$covered.\\

\begin{description}
\item [Case $2.1$ :] $\mu =2$. Since $\{1,n\}$ is not an interval of $G$, then there exists
 $j \in \{3, \ldots, n-2\}$ such that $n \rule[0.4mm]{0.6cm}{0.1mm}\, j$.
Let $m=max\{j\in\{3,\ldots,n-2\}: n \rule[0.4mm]{0.6cm}{0.1mm}\,
j\}$. That means $3 \leq m \leq n-2$ and $\{m+1,\ldots, n-1\} .\ldots\, n \rule[0.4mm]{0.6cm}{0.1mm}\,\,m$. We prove that $m=3$, so that $N_{G}(n) = \{2,3\}$. Clearly, $m\neq 4$, because $\{n,3\}$ is
not an interval of $G$. Suppose for a contradiction that $m\geq 5$. As $G[\{5,\ldots,n-1\}] \simeq P_{n-5}$ and
$n-5\geq 4$, then  $G [\{5,\ldots,n-1\}]$ is indecomposable.
Moreover, by $2$ of Theorem \ref{lemEhR}, $G[\{5,\ldots,n-1\}\cup\{2, n\}]$ is indecomposable
because $\{5,\ldots,n-1\} \, .\ldots \, 2$, $n\notin
\langle\{5,\ldots,n-1\}\rangle$ and $2 \rule[0.4mm]{0.6cm}{0.1mm}\,
n$. Set $W= \{5,\ldots,n-1\}\cup\{2, n\}$. We have $1\notin \langle
W\rangle$ because $n\, .\ldots\, 1
\rule[0.4mm]{0.6cm}{0.1mm}\,\,2$. Besides, $1\notin W(2)$ because $1\,
.\ldots\,n\rule[0.4mm]{0.6cm}{0.1mm}\,\, 2$ and $1\notin W(n)$ because $1\,
.\ldots\,m\rule[0.4mm]{0.6cm}{0.1mm}\,\, n$. Moreover, as $1 \in
\{5,\ldots n-1\}^{-}$, we obtain : $1 \notin \bigcup \limits_{i=5}^{i=n-1} W(i).$\\
Then, $1\notin \langle W\rangle$ and $n \notin W(i)$ for all $i \in W$. Thus,  by Theorem \ref{lemEhR}, $G-\{3,
4\}$ is indecomposable and $\mathbb{I}(G)$ would not be $\{1,n-1\}$-covered, that is impossible, so $m \leq 4$.

\item [Case $2.2$ :] $3 \leq \mu \leq n-4$. Set $X = \{1,\, \ldots,\mu\} \, \cup \{n\}$. It is clear that
$G[X] \simeq P_{\mu+1}$ with $\mu+1 \geq
4$. Thus, $G[X]$ is indecomposable. Since $\,3\leq \mu\leq
n-4$, $\{\mu+1,\mu+2\} \cap \{1,n-1\} = \emptyset$ then $G-\{\mu+1,\mu+2\}$ is
decomposable. Consider a non-trivial interval $I$ of
$G-\{\mu+1,\mu+2\}$.\\

\begin{itemize}
\item Assume that $I \cap X \neq \emptyset$. As $G[X]$ is
indecomposable, we have either $I \cap X = \{u\}, \,$ where $u \in X,$ or $I \cap X = X$.
Suppose that $I \cap X = \{u\}$. Since $\mid I\mid \geq 2$, there exists
$v \in \{\mu+3,\, \ldots, n-1\} \cap I$. As $v\, .\ldots\, 2 \rule[0.4mm]{0.6cm}{0.1mm}\,\,1$ and $v\, .\ldots\, \mu
\rule[0.4mm]{0.6cm}{0.1mm}\,\,n$, then $u \notin \{1,n\}$. Moreover, $u \notin \{2, \ldots,\mu\}$ because otherwise, $u-1 \notin I$ and $v\, .\ldots\, u-1
\rule[0.4mm]{0.6cm}{0.1mm}\,\,u$. It follows that $I \cap X = X$.\\
Now, we prove that $I \cap \{\mu+3,\ldots, n-1\} = \emptyset$. As $I$ is a non-trivial interval, there exists $k\in \{\mu+3,\ldots,
n-1\} \setminus I$. Let $\nu= min(\{\mu+3,\ldots, n-1\} \setminus I)$. If $\nu> \mu+3$, we obtain that $1\, .\ldots\, \nu
\rule[0.4mm]{0.6cm}{0.1mm}\,\,\nu-1$; which contradicts the fact that $\nu-1, 1\in I$
and $\nu \notin I$. Thus, $\nu = \mu +3$. Besides, since for each $k, \, \mu+4 \leq k\leq n-1$,
we have $1\, .\ldots\, k-1
\rule[0.4mm]{0.6cm}{0.1mm}\,\,k$, then we prove by induction that
$\,k\notin I$ for each $k, \, \mu+4 \leq k\leq
n-1$.\\
Thus, $I = X$ and then $\{\mu+3,\ldots, n-1\} .\ldots n$. Besides, $\{n, \mu+1\}$ is not
an interval of  $G$, then $\mu+2 .\ldots n$. We obtain that either $N_{G}(n) =  \{\mu\} \, when \,\, \mu+1 .\ldots\, n$ or $N_{G}(n)=\{\mu, \mu+1\} \, when \,\,\mu+1\rule[0.4mm]{0.6cm}{0.1mm}\,\, n$.\\

\item Assume that $I \cap X = \emptyset$. As $I$ is a non-trivial interval of  $G-\{\mu+1, \mu+2\}$, we have
$I \subset \{\mu+3,\ldots,n-1\}$ and $\mid I\mid \geq 2$, then $n\geq \mu+5$. We distinguish two cases.
\begin{itemize}

\item $n\geq \mu+7$ or $n=\mu+5$. In this case, since
$G[\{\mu+3,\ldots,n-1\}]$ is  indecomposable, then
$I=\{\mu+3,\ldots,n-1\}$. As $n\notin I$ and
$n-1 .\ldots n$, then  $\{\mu+3,\ldots,n-1\} .\ldots n$. Moreover, $\{n,
\mu+1\}$ is not an interval of  $G$, then $\mu+2 .\ldots n$ and we get that either $N_{G}(n) =  \{\mu\} \, when \,\, \mu+1 .\ldots\, n$ or $N_{G}(n)=\{\mu, \mu+1\} \, when \,\,\mu+1\rule[0.4mm]{0.6cm}{0.1mm}\,\, n$.\\

\item $n= \mu+6$. We have, $I$ is a
non-singleton interval of $G[\{n-3, n-2, n-1\}]$. It follows that either
$I= \{n-3, n-2, n-1\}$ or $I= \{n-3,n-1\}$. In the first case, we get as previously, that either $N_{G}(n) =  \{\mu\} \, when \,\, \mu+1 .\ldots\, n$ or $N_{G}(n)=\{\mu, \mu+1\} \, when \,\,\mu+1\rule[0.4mm]{0.6cm}{0.1mm}\,\, n$.\\

In the second case, as $n-1 .\ldots n, \, n \notin I$ and $n-3 \in I$, then
$n-3 .\ldots n$. Moreover, we show that $n-2 .\ldots n$.
Suppose for a contradiction that $n-2\rule[0.4mm]{0.6cm}{0.1mm}\,\, n$
and set $Y=\{1, \ldots, \mu+1\}$. Clearly, as $n \geq 9, \, \, G[Y]$ is indecomposable.
Besides, $n-2 \in \langle Y\rangle, \, n \notin \langle Y\rangle$ (because $1
.\ldots n \rule[0.4mm]{0.6cm}{0.1mm}\,\,\mu)$ and $n-2
\rule[0.4mm]{0.6cm}{0.1mm}\,\,n$. Thus, by $2$ of Theorem \ref{lemEhR}, $G[Y \cup \{n-2,n\}]$ is
indecomposable. Set $Z = Y \cup \{n-2,n\}$. We have $n-1 \notin Z(u)$
for $u \in Y$ because $n-1 \in Y^{-}$. Moreover, $n-1 \notin
Z(n-2)$, $n-1 \notin Z(n)$ and $n-1
\notin \langle Z\rangle$. Hence, by Theorem \ref{lemEhR}, $n-1 \in Ext(Z)$, that is, $G-\{n-3, n-4\}$ is
indecomposable; which contradicts the fact that $\mathbb{I}(G)$ is $\{1,n-1\}$-covered. Therefore, $n\, .\ldots
\{n-1,n-2,n-3\}$ and, as seen above, $\mu+2 .\ldots n$. So, we have either $N_{G}(n) =  \{\mu\} \, when \,\, \mu+1 .\ldots\, n$ or $N_{G}(n)=\{\mu, \mu+1\} \, when \,\,\mu+1\rule[0.4mm]{0.6cm}{0.1mm}\,\, n$.\\
\end{itemize}

\end{itemize}

\end{description}

\end{itemize}

{\hspace*{\fill}$\Box$\medskip}

\subsection{The class $\mathcal{Q}_{-1}$}
The next proposition describes the class $\mathcal{Q}_{-1}$.
\begin{proposition} \label{propqmoins1} Given a graph $G$ defined on $\mathbb{N}_{n}$,
where $n \geq 11$, $G \in \mathcal{Q}_{-1}$\, if and only if \,$G-n = Q_{n-1}$ and either $N_{G}(n) = \mathbb{N}_{n-3} \cup \{n-1\}$ or $N_{G}(n) = \{n-3\}$ or $N_{G}(n) =\{n-2,k\} \, where \,\,k \in \mathbb{N}_{n-3} \setminus \{2\}$ or $N_{G}(n) =\{n-2,k,k+1\} \, where \,\,k \in \{2,\ldots,n-4\}.$\\

\end{proposition}

\noindent{\em Proof\/.}
Given a graph $G$ defined on
$\mathbb{N}_{n}$, where $n \geq 11$ such that $G-n = Q_{n-1}$.
Suppose that $N_{G}(n) \in \{\mathbb{N}_{n-3} \cup \{n-1\}\} \cup
\{\{n-3\}\} \cup \{\{n-2,k\}: k \in \mathbb{N}_{n-3} \setminus
\{2\}\} \cup \{\{n-2,k,k+1\}: k \in \{2,\ldots,n-4\}\}$. To verify that $G$ is indecomposable and $\mathbb{I}(G)$ is $\{1,n-1\}$-covered, we proceed
as it is done at the beginning of the proof of Proposition \ref{proppmoins1}.\\

Conversely, assume that $G \in \mathcal{Q}_{-1}$. Set $X=\mathbb{N}_{n-3}$. We have $G[X] = P_{n-3}$ is indecomposable. Clearly, $n-1 \in \langle X \rangle$. Similarly, set $Y=\mathbb{N}_{n-4}$. We have $G[Y] = P_{n-4}$ is indecomposable and $n-1 \in \langle Y \rangle$. Also set $Z = \{4, \ldots, n-1\}$. We have $G[Z] \simeq Q_{n-4}$ is indecomposable. Observe that $1 \in Z(n-1)$.\\
Let $u \in Y$. For a contradiction, suppose that $n \in X(u)$. We have $n \in Y(u)$ as well. Since $\mathbb{I}(G)$ is $\{1,n-1\}$-covered, $G-\{n-2,n-3\}=G[Y \cup \{n-1,n\}]$ is decomposable. By
Theorem \ref{lemEhR}, $\{u,n\}$ is an interval of $G[Y \cup \{n-1,n\}]$. In particular, $n .\ldots n-1$. Now, we prove that $n \rule[0.4mm]{0.6cm}{0.1mm}\,\, n-2$, which implies that $\{u,n\}$ would be a non-trivial interval of $G$.\\
We distinguish the following two cases.
\begin{itemize}
  \item \underline{Case 1: } $u \neq 1$. Set $Y'=(Y \setminus \{u\}) \cup \{n\}$. As $\{u,n\}$ is an interval of $G[Y \cup \{n\}]$, $G[Y] \simeq G[Y']$ and hence $G[Y']$ is indecomposable. We have $n-1 \in \langle Y' \rangle$ because $n .\ldots n-1$. Since $\mathbb{I}(G)$ is $\{1,n-1\}$-covered, $G-\{u,n-3\}=G[Y' \cup \{n-1,n-2\}]$ is decomposable. It follows from
  Theorem \ref{lemEhR} that $n-2 \in \langle Y' \rangle$. In particular $n \rule[0.4mm]{0.6cm}{0.1mm}\,\, n-2$.
 \item \underline{Case 2: } $u = 1$. For a contradiction, suppose that $n-2 .\ldots n$. Set $Z'=\mathbb{N}_{n-5}$. The graph $G[Z']= P_{n-5}$ is indecomposable. Moreover, $n \in Z'(1), \, n-2 \in \langle Z'\rangle$ and $n .\ldots \,n-2 \rule[0.4mm]{0.6cm}{0.1mm}\,\,1$. It follows from
 Theorem \ref{lemEhR}, that $G[Z' \cup \{n-2,n\}]$ is indecomposable. Set $Z"=Z' \cup \{n-2,n\}$. We have $n-1 \notin \langle Z" \rangle$ because $n .\ldots \,n-1 \rule[0.4mm]{0.6cm}{0.1mm}\,\,n-2$. Since $n-1 .\ldots Z', \, n-2 \rule[0.4mm]{0.6cm}{0.1mm}\,\,Z'$, then $n-1 \notin Z"(n-2)$. Besides, $n-1 \notin Z"(n)$ because otherwise, as $n-1 .\ldots Z'$ we obtain that $n .\ldots Z'$ and then $n \in \langle Z' \rangle$; which contradicts the fact that $n \in Z'(1)$. Hence, $n-1 \notin Z"(n-2) \cup Z"(n)$. In addition,
 by Theorem \ref{lemEhR}, $n-1 \notin Z'(v)$ for $v \in Z'$ because $n-1 \in \langle Z'\rangle$. Thus, $n-1 \notin Z"(v)$ for $v \in Z"$. It follows from
 Theorem \ref{lemEhR}, that $n-1 \in Ext(Z")$. Thus, $G[Z" \cup \{n-1\}]=G-\{n-4,n-3\}$ is indecomposable which contradicts the fact that $\mathbb{I}(G)$ is $\{1,n-1\}$-covered.
\end{itemize}

Consequently, $n \notin \bigcup \limits_{u=1}^{u=n-4} X(u).$\\

\noindent Since $p_{X}$ is a partition of $\{n-2,n-1,n\}$, then by Theorem \ref{lemEhR}: $n \in X(n-3) \cup \langle X \rangle \cup Ext(X).$\\
We distinguish the following cases.
\begin{itemize}
\item \underline{Case 1: } $n \in X(n-3)$. As $n \in X(n-3)$, $n-4 \rule[0.4mm]{0.6cm}{0.1mm}\,\,n$. Thus, $n-1 .\ldots \,n-4 \rule[0.4mm]{0.6cm}{0.1mm}\,\,n$ and so, $n \notin Z(n-1)$. Furthermore, $G[Z \cup \{1,n\}]=G-\{2,3\}$ is decomposable because $\mathbb{I}(G)$ is $\{1,n-1\}$-covered. Since $1 \in Z(n-1)$, $\{1,n-1\}$ is an interval of $G[Z \cup \{1,n\}]$. In particular, $n .\ldots n-1$, so that $\{n-3,n\}$ is an interval of $G-\{n-2\}$. Therefore, $\{n-3,n\}$ is not an interval of $G[\{n-3,n-2,n\}]$ and $n \rule[0.4mm]{0.6cm}{0.1mm}\,\,n-2$. For $k=n-4$, we obtain that either $N_{G}(n) =  \{n-2,k\} \, when \,\, n-3 .\ldots\, n $ or $N_{G}(n)=\{n-2,k,k+1\} \, when \,\,n-3 \rule[0.4mm]{0.6cm}{0.1mm}\,\, n$.

\item \underline{Case 2: } $n \in \langle X \rangle$. Suppose for a contradiction that $n \in Z(n-1)$. We have $\{4, \ldots, n-3\} .\ldots n \rule[0.4mm]{0.6cm}{0.1mm}\,\, n-2$. As $n \in \langle X \rangle$, we obtain $\{1, \ldots, n-3\} .\ldots n \rule[0.4mm]{0.6cm}{0.1mm}\,\, n-2$ and
$\{n-1,n\}$ would be a non-trivial interval of $G$. Thus, $n \notin Z(n-1)$. Since $\mathbb{I}(G)$ is $\{1,n-1\}$-covered, $G-\{2,3\}=G[Z \cup \{1,n\}]$ is decomposable. As $1 \in Z(n-1)$ and $n \notin Z(n-1)$, $\{1,n-1\}$ is an interval of $G[Z \cup \{1,n\}]$. We get either $\{1,n-1\} \rule[0.4mm]{0.6cm}{0.1mm}\,\, n$ or $\{1,n-1\} .\ldots\, n$. Suppose for a contradiction that $\{1,n-1\} .\ldots\, n$. Since $n \in \langle X \rangle$, $\mathbb{N}_{n-3} .\ldots \, n$. If $n-2 .\ldots \, n$, then $\mathbb{N}_{n-1}$ would be a non-trivial interval of $G$, and if $n-2 \rule[0.4mm]{0.6cm}{0.1mm}\,\, n$, then $\{n-1,n\}$ would be a non-trivial interval of $G$. Therefore, $\{1,n-1\} \rule[0.4mm]{0.6cm}{0.1mm}\,\, n$. As $n \in \langle X \rangle$, $\mathbb{N}_{n-3} \,\rule[0.4mm]{0.6cm}{0.1mm}\,\, n$. Since $\mathbb{N}_{n-1}$ is not an interval of $G$, $n \,.\ldots \, n-2$ and : $N_{G}(n) = \mathbb{N}_{n-3} \cup \{n-1\}$.

\item \underline{Case 3: } $n \in Ext(X)$. For a contradiction, suppose that $n \rule[0.4mm]{0.6cm}{0.1mm}\,\, n-1$. As $\mathbb{I}(G)$ is $\{1,n-1\}$-covered, $G-\{n-2,n-3\}=G[Y \cup \{n-1,n\}]$ is decomposable. Since $Y .\ldots
 n-1 \rule[0.4mm]{0.6cm}{0.1mm}\,\, n$, it follows from Theorem \ref{lemEhR}, that $n \in \langle Y \rangle$. Furthermore, as $n \in Ext(X)$, $G[X \cup \{n\}]=G[Y \cup \{n-3,n\}]$ is indecomposable. Thus, we get either $Y .\ldots
 n \rule[0.4mm]{0.6cm}{0.1mm}\,\, n-3$ or $n-3 .\ldots
 n \rule[0.4mm]{0.6cm}{0.1mm}\,\, Y$. If $n-3 .\ldots
 n \rule[0.4mm]{0.6cm}{0.1mm}\,\, Y$, then $\{n-2,n\}$ would be a non-trivial interval of $G$. Suppose that $Y .\ldots
 n \rule[0.4mm]{0.6cm}{0.1mm}\,\, n-3$. Since $n-1 .\ldots
 n-3 \rule[0.4mm]{0.6cm}{0.1mm}\,\, n$, $n \notin Z(n-1)$. As $1 \in Z(n-1)$ and $1 .\ldots
 n \rule[0.4mm]{0.6cm}{0.1mm}\,\, n-1$, it would follow from Theorem \ref{lemEhR} that $G[Z \cup \{1,n\}]=G-\{2,3\}$ is indecomposable and $\mathbb{I}(G)$ would not be $\{1,n-1\}$-covered. Consequently, $n \,.\ldots \,n-1.$\\

 Lastly, consider $X'=X \cup \{n\}$. We have $G[X']$ is indecomposable because $n \in Ext(X)$. We verify that $\mathbb{I}(G[X'])$ is $\{1,n-3\}$-covered. Otherwise, there exist $x \neq y \in X' \setminus \{1,n-3\}$ such that $G[X']-\{x,y\}$ is indecomposable. Set $Y'=X' \setminus \{x,y\}$. We have $Y' \,.\ldots n-1$ because $n .\ldots n-1$. Therefore, $n-1 \in \langle Y' \rangle$. Moreover, $n-2 \notin \langle Y' \rangle$ because $n-3 .\ldots n-2 \rule[0.4mm]{0.6cm}{0.1mm}\,\, 1$. Since $Y' .\ldots
 n-1 \rule[0.4mm]{0.6cm}{0.1mm}\,\, n-2$, $G[Y' \cup \{n-2,n-1\}]$ is indecomposable by Theorem \ref{lemEhR}. As $x,y \in X' \setminus \{1,n-3\} \subseteq \mathbb{N}_{n} \setminus \{1,n-1\}$, $\mathbb{I}(G)$ would not be $\{1,n-1\}$-covered. It follows that $ \mathbb{I}(G[X'])\,\,  is \,\, \{1,n-3\}-covered.$\\
 By Proposition \ref{proppmoins1}, there is either $k \in \{3, \ldots ,
n-5\} \cup \{1,n-3\}$ such that $N_{G[X']}(n) = \{k\}$ or
$k \in \{2, \ldots , n-5\}$ such that $N_{G[X']}(n) = \{k,k+1\}$.

\begin{itemize}
  \item Assume that $N_{G[X']}(n) = \{k\}$ where $k \in \{3, \ldots ,
n-5\} \cup \{1,n-3\}$. We prove that if $k \neq n-3$, then $n
\rule[0.4mm]{0.6cm}{0.1mm}\,\, n-2$. Suppose, by contradiction, that $k \leq n-5$ and
$n .\ldots n-2$. Clearly, $G[\mathbb{N}_{n-5}] = P_{n-5}$ is indecomposable. Besides, $\mathbb{N}_{n-5}  \rule[0.4mm]{0.6cm}{0.1mm}\,\, n-2,
\, n \notin \langle \mathbb{N}_{n-5} \rangle$ (because $2 .\ldots n
\rule[0.4mm]{0.6cm}{0.1mm}\,\, k$) and $n-2 .\ldots n$. It follows
that $G[\mathbb{N}_{n-5} \cup \{n-2,n\}]$ is indecomposable. Set $W
= \mathbb{N}_{n-5} \cup \{n-2,n\}$. We have $ n-1 \notin \langle W
\rangle$ (because $n .\ldots n-1 \rule[0.4mm]{0.6cm}{0.1mm}\,\,
n-2$) and $n-1 \notin W(v)$ where $v \in \{n-2,n\}$ (because $n-1
.\ldots k \rule[0.4mm]{0.6cm}{0.1mm}\,\, \{n-2,n\}$). Moreover, $n-1
\notin W(v), v \in \{1, \ldots ,n-5\}$ because $n-1 .\ldots
\mathbb{N}_{n-5}$. It results that $n-1 \in Ext(W)$ and then $G-\{n-4,n-3\}$ is indecomposable; which contradicts the fact that $\mathbb{I}(G)$ is $\{1,n-1\}$-covered. Thus, we get that either $N_{G}(n) =  \{n-3\} \, when \,\, n-2 .\ldots\, n $ or $N_{G}(n)= \{k,n-2\}; \,k \in  \mathbb{N}_{n-3} \setminus \{2,n-4\} \, when \,\,n-2 \rule[0.4mm]{0.6cm}{0.1mm}\,\, n$.\\

  \item Assume that $N_{G[X']}(n) = \{k,k+1\}$ where $k \in \{2, \ldots , n-5\}$. Similarly, we show that
$n \rule[0.4mm]{0.6cm}{0.1mm}\,\, n-2$ and thus $N_{G}(n) =
\{k,k+1,n-2\}$ where $k \in \{2, \ldots , n-5\}$.
\end{itemize}

\end{itemize}
{\hspace*{\fill}$\Box$\medskip}

%%%%%%%%%%%%%%%%%%%%%%%%%%%%%%%%%%%%%%%%%%%%%%%%%%%%%%%%%%%%%%%%%%%%%%%%%%%%%%%%%%%%%%%%%%%%%%%%%%%%%%%%%%%%%%%%%%%%%%%%%%%%%%%%%%%%%%%%%%
\subsection{The class $\mathcal{P}_{-3}$}
With the aim  to describe  the class $\mathcal{P}_{-3}$, we
introduce  the  class $\mathcal{G}$  of graphs $G$ defined on
$\mathbb{N}_{n}$, where $n\geq 12$, such that $G[X]=P_{n-3}$, $n-2
\in X^{-}$, where $X=\mathbb{N}_{n-3}$ and satisfying one and only
one of the following assertions :

\begin{itemize}
\item $n-1 \in X^{-}, X(1) = \{n\}, E(G_{X}) = \{\{n,n-1\}\}$ and $n-1 \rule[0.4mm]{0.6cm}{0.1mm}\,\,
n-2$.
\item $X(2) = \{n-1\}, X(1) = \{n\}, \mid E(G_{X})\mid \geq 2$ and $1 .\ldots n$ if $\{n,n-1\} \in E(G_{X})$.
\item $X(u) = \{n-1,n\}$ where $u \in \{1,2\}$ and either $\mid E(G_{X})\mid = 2$ with $u \not\sim \{n-1,n\}$, or $E(G_{X})=\{\{n,n-2\}\}$ with $n \not\sim \{n-1,u\}.$
\end{itemize}

\begin{proposition}
The indecomposability graph of graphs of the  class $\mathcal{G}$
are $\{1, n-3\}$-covered.
\end{proposition}

\noindent{\em Proof\/. } Let  $G$ be a graph of the class
$\mathcal{G}$. Given  $i< j\in \mathbb{N}_{n} \setminus \{1, n-3\}$,
then $G-\{i, j\}$ is decomposable. Indeed, if $3 \leq j\leq n-5$,
then $\{j+1, \ldots ,n-3\}$ would be a non trivial interval of
$G-\{i,j\}$. As more, if $j=n-4$, $\mathbb{N}_{n} \setminus
\{i,j,n-3\}$ is a non trivial interval of $G-\{i,j\}$.\\
If $j \geq n-2$ and $i \geq n-2$, then $G-\{i, j\}$ is
decomposable because $\rm{Ext}(\mathbb{N}_{n-3})=\emptyset$. \\
If $j\geq n-2$ and $ 2 \leq i \leq n-4$, we have to examine the
following cases.
\begin{itemize}
  \item $3 \leq i \leq n-4$. Clearly, $(\mathbb{N}_{i-1} \cup \{n,n-1,n-2\}) \setminus \{j\}$ would be a non trivial interval of
$G-\{i,j\}$.
  \item $i=2$. We have to distinguish the following three cases
according to $\mathbb{N}_{n-3}(2)$.
\begin{itemize}
  \item $\mathbb{N}_{n-3}(2) = \emptyset$. In this case, $\{1,n,n-1,n-2\} \setminus
\{j\}$ would be a non trivial interval of $G-\{i,j\}$.
  \item $\mid \mathbb{N}_{n-3}(2) \mid = 1.$ If $j=n-1$, it is clear
that $\{1,n,n-2\}$ is a non trivial interval of $G-\{i,j\}$. If
$j=n-2$ (resp. $j=n$), then $\{3, .\ldots, n-3\} \cup \{1,n-1\}$ is
a non trivial interval of $G-\{i,j\}$ if $n .\ldots n-1$ (resp. $n-2
.\ldots n-1$) and $\{1,n\}$ (resp. $\{1,n-2\}$) is a non trivial
interval of $G-\{i,j\}$ if $n \rule[0.4mm]{0.6cm}{0.1mm}\,\, n-1$
(resp. $n-2 \rule[0.4mm]{0.6cm}{0.1mm}\,\, n-1$).
  \item $\mid \mathbb{N}_{n-3}(2) \mid = 2.$ If $j=n-2, \,
\{n,n-1\}$ is a non trivial interval of $G-\{i,j\}$. If $j \neq n-2$
and $n-1 \rule[0.4mm]{0.6cm}{0.1mm}\,\, n-2$, then $\{1,n-2\}$ would
be a non trivial interval of $G-\{i,j\}$. If $j=n$ (resp. $j=n-1$)
and $n-1 .\ldots n-2$, then $\{3, .\ldots, n-3\} \cup \{1,n-1\}$
(resp. $\{1,n-2\} $) is a non trivial interval of $G-\{i,j\}$.
\end{itemize}
\end{itemize}
It remains to prove that the graph $G$ is indecomposable. We pose
$X=\mathbb{N}_{n-3}$. First, assume that $\mid X(u)\mid = 2$ where
$u \in \{1,2\}$. If $n-2 \,\rule[0.4mm]{0.6cm}{0.1mm}\,\,\{n-1,n\}$
(resp. $n-1 .\ldots n-2 \,\rule[0.4mm]{0.6cm}{0.1mm}\, n$), then
Theorem \ref{lemEhR} claims that $G[X \cup \{n-2,n\}]$ is
indecomposable. Set $W = X \cup \{n-2,n\}$. As $n-1 \notin \langle
X\rangle$, then $n-1 \notin \langle W\rangle$. Besides, %$n-1 \notin
%W(u)$ because $u .\ldots n-2 \,\rule[0.4mm]{0.6cm}{0.1mm}\, n-1$
%(resp. because $n \not\sim \{n-1,u\}$) and $n-1 \notin W(n)$ because
%$n-1 .\ldots u \,\rule[0.4mm]{0.6cm}{0.1mm}\, n$ (resp. because $n-1
%.\ldots n-2 \,\rule[0.4mm]{0.6cm}{0.1mm}\, n$).
we verify that there is no $i \in W$ such that $n-1 \in W(i)$ and so
$G$ is indecomposable. Finally, assume that $\mid X(1)\mid = 1$. If
$n \,\rule[0.4mm]{0.6cm}{0.1mm}\,\, n-1$ (resp. $n .\ldots n-1$),
then by Theorem \ref{lemEhR}, $G[\mathbb{N}_{n-1}]$ is indecomposable
(resp. $G[X \cup \{n,n-1\}]$ is indecomposable). It is clear that $n
\notin \langle \mathbb{N}_{n-1} \rangle$ (resp. $n-2 \notin \langle
Z \rangle$ where $Z = X \cup \{n,n-1\}$). Moreover, we verify as
seen above that there is no $i \in \mathbb{N}_{n-1}$ (resp. no $i
\in Z$) such that $n \in \mathbb{N}_{n-1}(i)$ (resp. such that $n-2
\in Z(i)$). Hence, $n \in Ext(\mathbb{N}_{n-1})$(resp. $n-2 \in
Ext(Z)$) which allows us to conclude.
{\hspace*{\fill}$\Box$\medskip}

\begin{proposition} \label{proppmoins3}
Up to isomorphism, the graphs of the class $\mathcal{P}_{-3}$ are
those of cardinality  $\geq 12$ of the class $\mathcal{G}$.
\end{proposition}

For the proof of this proposition, we need the following lemma.
\begin{lemma}\label{lempmoins3}
Let $G$ be an indecomposable graph defined on $\mathbb{N}_{n-1}$
where $n\geq 12$, verifying : $G[\mathbb{N}_{n-3}]=P_{n-3}$ and for
each vertex $i$ of $\mathbb{N}_{n-1} \setminus \{1, n-3\}$, $i$ is critical.
Then one and only one of the following assertions holds:

\begin{enumerate}
\item $n-1 \in \mathbb{N}_{n-3}^{-}$ and $n-2\in
\mathbb{N}_{n-3}(u)$, where  $u\in\{1, 2, n-4, n-3\}$.

\item $n-1\in \mathbb{N}_{n-3}(u)$ and $n-2\in \mathbb{N}_{n-3}(v)$, where $\{u,v\}=\{n-4, n-3\}$ or $\{u,v\}=\{1, 2\}$. Moreover, if  $u=n-4$
and  $v=n-3$, (resp. $u=2$ and $v=1$), then $n-2\, .\ldots\, \{n-1,
n-3\}$ (resp. $n-2\, .\ldots\, \{1, n-1\}$).

\end{enumerate}
\end{lemma}

\noindent{\em Proof\/. }% Firstly, notice that since
%$G[\mathbb{N}_{n-3}]=P_{n-3}$ is indecomposable and for each vertex
%$i$ of $\mathbb{N}_{n-1}-\{1, n-3\}$, $i$ is critical, then
%$Ext(\mathbb{N}_{n-3})=\emptyset$.
Assume that $n-1\, .\ldots\,\, \mathbb{N}_{n-3}$. As $n-1$ is a
critical vertex of $G$, then $n-2\notin Ext(\mathbb{N}_{n-3})$.
Moreover, since  $G$ is indecomposable, $n-2\notin
\langle\mathbb{N}_{n-3}\rangle$. Thus, $n-2\in\mathbb{N}_{n-3}(u)$
where  $u\in \mathbb{N}_{n-3}$ and
$n-2\rule[0.4mm]{0.6cm}{0.1mm}\,\, n-1$. Suppose by contradiction
that $\, u \in \{3,\, \ldots ,n-5\}$. We prove that $u$ is not a
critical vertex of $G$ which is impossible. In fact, consider the
bijection $f\, : \, \mathbb{N}_{n-1}-\{u\}\, \longrightarrow \,
\mathbb{N}_{n-2}$ defined for each $i \in \mathbb{N}_{u-1} \cup
\{u+1, \ldots, n-3\}$, $f(i) = i, \, f(n-1) = n-1, \, f(n-2) = u$ is
an isomorphism from $G-\{u\}$ onto an element of $\mathcal{P}_{-1}$.
So, $G-\{u\}$ is indecomposable; contradiction. Hence, $n-2\in
\mathbb{N}_{n-3}(u)$ where  $u\in\{1, 2, n-4, n-3\}$.\\

It may be assumed now that $u\neq v\in\mathbb{N}_{n-3}$ such that
$n-1\in\mathbb{N}_{n-3}(u)$ and $n-2\in\mathbb{N}_{n-3}(v)$. First,
suppose that $\{u, v\}\neq\{1,n-3\}$. Without loss of generality, we
assume that $u\notin \{1, n-3\}$. Set $Y=(\mathbb{N}_{n-3} \setminus
\{u\})\cup\{n-1\}$. Clearly, $G[Y]\simeq P_{n-3}$. As  $u$ is a
critical vertex, then $n-2\notin Ext(Y)$. Suppose that $n-2\in
\langle Y\rangle$. If $n-2\rule[0.4mm]{0.6cm}{0.1mm}\,\, Y$. In this
case, $v \rule[0.4mm]{0.6cm}{0.1mm}\,\, \mathbb{N}_{n-3}-\{u,v\}$
and $u\, .\ldots \,v$. Thus, $\mid N_{G[\mathbb{N}_{n-3}]} (v)\mid
\, > 2$ which is impossible because $G[\mathbb{N}_{n-3}] = P_{n-3}$.
If $n-2 \, .\ldots Y$. In this case, $v\, .\ldots
\,\mathbb{N}_{n-3}-\{u,v\}$ and $u\rule[0.4mm]{0.6cm}{0.1mm}\,\, v$;
which implies that $v=1$ and $u=2$ or $v=n-3$ and $u=n-4$. If $v=1$
and $u=2$ (resp. $v=n-3$ and $u=n-4$), as $n-2\, .\ldots \, Y$, then
$n-2\, \ldots \, \{n-1,1\}$
(resp. $n-2\, .\ldots \, \{n-1,n-3\}$).\\
Assume in the sequel that $n-2\notin \langle Y\rangle$. So, there is
$w\in Y$ such that $n-2\in Y(w)$. $v\neq w$, because $\{v,n-2\}$ is
not a interval of $G[\{v,n-2,n-1\}]$. As more
$n-2\in\mathbb{N}_{n-3}(v)$, then $\{v,w\}$ is an interval of
$G[\mathbb{N}_{n-3}-\{u\}]$. Necessarily, $\{v, w\}=\{1, 2\}$ and
$u=3$ or $\{v, w\}=\{n-3, n-4\}$ and $u=n-5$ or $\{v, w\}=\{1, 3\}$
and $u=4$ or $\{v, w\}=\{n-3, n-5\}$ and $u=n-6$. By isomorphism, it
suffices to study the case where $\{v, w\}=\{n-3, n-4\}$ and
$u=n-5$, and the case where $\{v, w\}=\{n-3, n-5\}$ and $u=n-6$.\\
Now suppose that $\{v, w\}=\{n-4, n-3\}$. Using Theorem
\ref{lemEhR},we demonstrate, that $G-\{n-4\}$ is indecomposable
which is impossible.
%In the first case, if $n-5\rule[0.4mm]{0.6cm}{0.1mm}\,\, n-1$ (resp.
%$n-1\, .\ldots \, n-5$), then $F=G-\{n-4\} \in
%\mathcal{P}_{-1}$($F-\{n-1\} \simeq P_{n-1}$ and $N_{F}(n-1) =
%\{n-6,n-5\}$ (resp. $N_{F}(n-1) = \{n-6\})$). In the second case,
%since $\{n-4,n-3\}$ is an interval of $G-\{n-5\}$, then
%$n-2\rule[0.4mm]{0.6cm}{0.1mm}\,\, n-3$. Besides,
%$G[\mathbb{N}_{n-5}]$ is indecomposable, $n-1\in
%\mathbb{N}_{n-5}(n-5)$, $n-2\, .\ldots \,\mathbb{N}_{n-5}$ and
%$n-1\rule[0.4mm]{0.6cm}{0.1mm}\,\, n-2$. Hence,
%$G[\mathbb{N}_{n-5}\cup\{n-1, n-2\}]$ is indecomposable. Posing $Z=
%\mathbb{N}_{n-5}\cup\{n-1, n-2\}$. We have $n-3\notin \langle
%Z\rangle$ because $n-5 \, .\ldots
%\,n-3\rule[0.4mm]{0.6cm}{0.1mm}\,\, n-2$, $n-3\notin Z(n-1)$ (resp.
%$n-3\notin Z(n-2)$) because $n-3\, .\ldots \,
%n-6\rule[0.4mm]{0.6cm}{0.1mm}\,\, n-1$ (resp. $n-3\, .\ldots \,
%n-1\rule[0.4mm]{0.6cm}{0.1mm}\,\, n-2$). Moreover, as $n-3\in
%\langle\mathbb{N}_{n-5}\rangle$, there is no $i\in\mathbb{N}_{n-5}$
%such that $n-3\in Z(i)$. By Lemma \ref{lemEhR}, $G-\{n-4\}$ is
%indecomposable.\\
Presently suppose that $\{v, w\}=\{n-3, n-5\}$. Similarly, we
demonstrate by Theorem \ref{lemEhR}, that $G-\{n-5\}$ is indecomposable.%In the first case,
%if $n-1\rule[0.4mm]{0.6cm}{0.1mm}\,\, n-6$ (resp. $n-1\, .\ldots \,
%n-6$), then $F'=G-\{n-5\} \in \mathcal{P}_{-1}$($F'-\{n-1\} \simeq
%P_{n-1}$ and $N_{F'}(n-1) = \{n-7,n-6\}$ (resp. $N_{F'}(n-1) =
%\{n-7\})$). In the second case and as seen above, it is easy to
%verify that $G-\{n-5\}$ is indecomposable. Indeed,
%$G[\mathbb{N}_{n-6}]$ is indecomposable, $n-1 \notin \langle
%\mathbb{N}_{n-6}\rangle$, $n-2\, .\ldots \, \mathbb{N}_{n-6}$ and
%$n-1 \rule[0.4mm]{0.6cm}{0.1mm}\,\, n-2$. Hence,
%$G[\mathbb{N}_{n-6}\cup\{n-1, n-2\}]$ is indecomposable. Set $W=
%\mathbb{N}_{n-6}\cup\{n-1, n-2\}$. As $\{n-3,n-5\}$ is an interval
%of $G-\{n-6\}$, we have $n-3\, .\ldots \, n-2$. Consequently, $n-3
%\in W^{-}$ and $n-4 \notin \langle W\rangle$ (because $n-1 .\ldots
%n-4 \rule[0.4mm]{0.6cm}{0.1mm}\,\, n-2$), and since $n-3
%\rule[0.4mm]{0.6cm}{0.1mm}\,\, n-4$, then $G[W \cup\{n-3, n-4\}]
%= G-\{n-5\}$ is indecomposable.\\
Finally, assume that  $u=1$ and $v=n-3$. Clearly,
$G[\{6,\ldots,n-3\}]$ is indecomposable because $n \geq 12$. As
$n-2\in\{6,\ldots, n-3\}(n-3)$, $\{6,\ldots,n-3\}\, .\ldots \, n-1$
and $n-1\rule[0.4mm]{0.6cm}{0.1mm}\,\, n-2$, then $G[\{6,\ldots,
n-3\}\cup\{n-1, n-2\}]$ is indecomposable. We pose $Z= \{6,\ldots,
n-3\}\cup\{n-1, n-2\}$. We have $3 \, .\ldots \, Z$, $2 \notin
\langle Z\rangle$ (because $n-3\, .\ldots \,
2\rule[0.4mm]{0.6cm}{0.1mm}\,\, n-1$) and
$2\rule[0.4mm]{0.6cm}{0.1mm}\,\, 3$. So, $G[Z \cup \{2,3\}]$ is
indecomposable. Set $Z' = Z \cup \{2,3\}$. We shall examine the two
cases.
\begin{itemize}
\item \underline{Case 1: }$ n-1 \, .\ldots 1$. We have $1 \in Z'(3), \, 4 \notin
Z'(3)$ (because $4\, .\ldots \, 2\rule[0.4mm]{0.6cm}{0.1mm}\,\, 3$)
and $1\, .\ldots \, 4\rule[0.4mm]{0.6cm}{0.1mm}\,\, 3$. Thus, $G[Z'
\cup \{3,4\}] = G-\{5\}$ is indecomposable; contradiction.
\item \underline{Case 2: }$ n-1 \rule[0.4mm]{0.6cm}{0.1mm}\,\, 1$. As
$1 \notin \langle Z'\rangle$ (because $3\, .\ldots \,
1\rule[0.4mm]{0.6cm}{0.1mm}\,\, 2$), $1 \notin Z'(n-1)$ (because
$1\, .\ldots \, n-2\rule[0.4mm]{0.6cm}{0.1mm}\,\, n-1$), $1 \notin
Z'(n-2)$ (because $n-2\, .\ldots \, 2\rule[0.4mm]{0.6cm}{0.1mm}\,\,
1$), $1 \notin Z'(2)$ (because $1\, .\ldots \,
3\rule[0.4mm]{0.6cm}{0.1mm}\,\, 2$), $1 \notin Z'(3)$ (because $3\,
.\ldots \, n-1\rule[0.4mm]{0.6cm}{0.1mm}\,\, 1$). So, there is no
$x\in Z'$ such that $1\in Z'(x)$ (because $1 .\ldots \{6, \ldots,
n-3\}$). Hence, $G[Z' \cup \{1\}]$ is indecomposable. Set $Z"=Z'
\cup \{1\}$. As $G[Z" \cup \{4\}]= G-\{5\}$, it suffices to
demonstrate that $G[Z" \cup \{4\}]$ is indecomposable. In fact, $4
\notin \langle Z"\rangle$(because $1\, .\ldots \,
4\rule[0.4mm]{0.6cm}{0.1mm}\,\, 3$), $4 \notin Z"(n-1)$ (because
$4\, .\ldots \, n-2\rule[0.4mm]{0.6cm}{0.1mm}\,\, n-1$), $4 \notin
Z"(n-2)$ (because $4\, .\ldots \, n-1\rule[0.4mm]{0.6cm}{0.1mm}\,\,
n-2$), $4 \notin Z"(2)$ (because $4\, .\ldots \,
1\rule[0.4mm]{0.6cm}{0.1mm}\,\, 2$), $4 \notin Z"(1)$ (resp. $4
\notin Z"(3)$) (because $4\, .\ldots \,
2\rule[0.4mm]{0.6cm}{0.1mm}\,\, 1$) (resp. (because $4\, .\ldots \,
2\rule[0.4mm]{0.6cm}{0.1mm}\,\, 3$)) and there is no $x\in Z"$ such
that $4\in Z"(x)$ (because $4 .\ldots \{6, \ldots, n-3\}$). By Theorem
\ref{lemEhR}, $G-\{5\}$ is indecomposable; impossible.

\end{itemize}
{\hspace*{\fill}$\Box$\medskip} \vspace*{0.5cm}

\noindent{\em Proof of  Proposition \ref{proppmoins3} \/. } Let $G$
be a graph defined on $\mathbb{N}_{n}$ where  $n\geq 12$. Assume
that $\mathbb{I}(G)$ is $\{1,n-3\}$-covered and
$G[\mathbb{N}_{n-3}]=P_{n-3}$. Set $X=\mathbb{N}_{n-3}$. Notice that
for $z\neq t\in\{n-2,n-1,n\}$, if $G[X \cup\{z,t\}]$ is
indecomposable, then  for each vertex
$i\in(X\cup\{z,t\}) \setminus \{1,n-3\}$, $i$ is a critical vertex of
$G[X\cup\{z,t\}]$. Even more notice that from Corollary 2.2, there
is  $x\neq y\in\{n-2, n-1, n\}$ such that $G[X\cup\{x,y\}]$ is
indecomposable. Without loss of generality, we may assume that
$\{x,y\}=\{n-1,n-2\}$. Let us distinguish the two cases.
\begin{itemize}

\item \underline{Case 1: }$n-1\in \langle X \rangle$
and $n-2\in X(u)$ where  $u\in X$. From what precedes, for each
vertex $i\in\mathbb{N}_{n-1} \setminus \{1,n-3\}$, $i$ is a critical vertex of
$G[\mathbb{N}_{n-1}]$. So, we may assume that $n-1\, .\ldots \, X$
and by Lemma \ref{lempmoins3}, $u\in \{1,2,n-4,n-3\}$. Consider the application
$f\, : \, \mathbb{N}_{n-1}\, \longrightarrow \, \mathbb{N}_{n-1}$
defined for each $i \in \mathbb{N}_{n-3}$, $f(i) = n-i+1, \, f(n-1)
= n-1$ and $f(n-2) = n-2$, we may assume that $u\in \{1,2\}$.

\begin{itemize}

\item If $n \rule[0.4mm]{0.6cm}{0.1mm}\,\,X$, then either
$n \, .\ldots \, n-1$ or $n \, .\ldots \, n-2$. Otherwise, $n
\rule[0.4mm]{0.6cm}{0.1mm}\,\,\mathbb{N}_{n-1}$, that is, $G$ is
decomposable; contradiction. Suppose that $n \, .\ldots \, n-1$
(resp. $n \, .\ldots \, n-2$), we show that $G-\{n-4,n-5\}$ is
indecomposable which is impossible. Indeed, we pose $Y = \{1,
\ldots, n-6\}$. Clearly, $G[Y]$ is indecomposable ($n \geq 12$).
$n-2 \in Y(u)$, $n-1 \in \langle Y \rangle$ and $n-1
\rule[0.4mm]{0.6cm}{0.1mm}\,\, n-2$. Hence, $G[Y \cup \{n-1,n-2\}]$
is indecomposable. Set $Z = Y \cup \{n-1,n-2\}$. Since $G[Z]$ is
indecomposable, $n \notin \langle Z\rangle$ because $n-1\, .\ldots
\, n\rule[0.4mm]{0.6cm}{0.1mm}\,\, n-6$ (resp. because $n-2\,
.\ldots \, n\rule[0.4mm]{0.6cm}{0.1mm}\,\, n-6$), $n-3 \in \langle
Z\rangle$(because $n-3\, .\ldots \, Z$) and
$n\rule[0.4mm]{0.6cm}{0.1mm}\,\, n-3$, this implies that $G[Z \cup
\{n-3,n\}] = G-\{n-4,n-5\}$ is indecomposable.

\item If $n \, .\ldots \,X$, then $n\, .\ldots \,
n-2$. Otherwise, $\{n,n-1\}$ is an interval of $G$; impossible.
Moreover, as $G$ is indecomposable, $n
\rule[0.4mm]{0.6cm}{0.1mm}\,\, n-1$. If $u=1$ (resp. $u=2$), then
$G$ is isomorphic to one of the elements of $\mathcal{G}$. It
suffices to consider the application $f\, : \, \mathbb{N}_{n}\,
\longrightarrow \, \mathbb{N}_{n}$ defined for each $i \in
\mathbb{N}_{n-3}$, $f(i) = i, \, f(n) = n-2, \, f(n-1) = n-1$ and
$f(n-2) = n$ (resp. $f\, : \, \mathbb{N}_{n}\, \longrightarrow \,
\mathbb{N}_{n}$ defined for each $i \in \mathbb{N}_{n-3} \setminus
\{2\}$, $f(i) = i, \, f(n-2) = 2, \, f(n-1) = n, \, f(2) = n-1$ and
$f(n) = n-2$).

\item $n\in X(v)$ where  $v\in X$. As
$\{n,v\}$ is not an interval of $G$, then either $v\, .\ldots
\,n-1 \rule[0.4mm]{0.6cm}{0.1mm}\,\, n$ or $n-2 \not\sim \{n,v\}$.\\

First, assume that $v\, .\ldots \,n-1 \rule[0.4mm]{0.6cm}{0.1mm}\,\,
n$. In this case, $G-\{n-2\}$ is indecomposable. Using Lemma \ref{lempmoins3}, $v
\in \{1,2,n-4,n-3\}$. If $v=u \in \{1,2\}$, we have $u \not\sim
\{n,n-2\}$ (because $\{n,n-2\}$ is not an interval of $G$). We may
assume that $n\, .\ldots \,u \rule[0.4mm]{0.6cm}{0.1mm}\,\, n-2$ and
$G$ is isomorphic to one of the elements of $\mathcal{G}$ by
permuting $n-1$ and $n-2$. If $v \neq u \in \{1,2\}$. Assume, for
instance, that $u=1, \, v=2$. If $n \rule[0.4mm]{0.6cm}{0.1mm}\,\,
n-2$, then $G$ is isomorphic to one of the elements of $
\mathcal{G}$. If $n .\ldots n-2$, then necessarily $n-2 .\ldots 1$
(otherwise, we verify by Theorem \ref{lemEhR}, that $G-\{2,n-1\}$ is
indecomposable; impossible). Hence, $G$ is isomorphic to one of the
elements of $ \mathcal{G}$. Suppose that $v
\in \{n-4,n-3\}$. We distinguish the two cases.\\
\hspace*{0.2cm} - If $u=1$, then $G-\{2,3\}$ is indecomposable if
$1\rule[0.4mm]{0.6cm}{0.1mm}\,\, n-2$ and
%Indeed, set $Z= \{4, \ldots, n-3\} \cup \{n,n-1\}$. Clearly, $G[Z]$
%is indecomposable. $1 \in Z^{-}, \, n-2 \notin \langle Z \rangle$
%and $1 \rule[0.4mm]{0.6cm}{0.1mm}\,\, n-2$.
$G-\{3,4\}$ is indecomposable if $1 .\ldots n-2$; impossible.\\
% In fact, set $W= \{5, \ldots, n-3\}
%\cup \{n,n-1\}$. It is clear that $G[W]$ is indecomposable and so
 %$G[W \cup \{n,n-2\}]$ is indecomposable. We pose $W'=W \cup \{n,n-2\}$.
%$1 \notin \langle W' \rangle$, $1 \notin W'(2)$, $1 \notin W'(n-2)$ and $1 \in W^{-}$, so $1 \in Ext(W')$.\\
\hspace*{0.2cm} - If $u=2$, we verify as previously, that
$G-\{2,3\}$ is indecomposable; impossible.\\

Now, assume that $n-2 \not\sim \{n,v\}$. We can assume that $n\,
.\ldots \,n-1$. If $u=1$, then since $G[X \cup \{n-2,n\}]$ is
indecomposable, Lemma \ref{lempmoins3} implies that $v \in \{1,2\}$. If $u=v=1$,
$G$ is isomorphic to one of the elements of $\mathcal{G}$. If $u=1$
and $v=2$, then $n-2\, .\ldots \,n$ and $G-\{n-1\}$ is
indecomposable. By Lemma \ref{lempmoins3}, $1\, .\ldots \,n-2$ and $G$ is
isomorphic to one of the elements of $\mathcal{G}$. If $u=2$, then
we obtain in the same manner that $v \in \{1,2\}$. If $u=v=2$, $G$
is isomorphic to one of the elements of $\mathcal{G}$. If $u=2$ and
$v=1$, then $n\, .\ldots \,n-2\rule[0.4mm]{0.6cm}{0.1mm}\,\, 1$ and
$G[X \cup \{n-2,n\}]$ is indecomposable. Lemma \ref{lempmoins3} implies that $n\,
.\ldots \,1$ and $G$ is isomorphic to one of the elements of
$\mathcal{G}$.
\end{itemize}

\item \underline{Case 2.}$\langle X \rangle = \emptyset$. We may assume that
$n-2\in X(u)$, $n-1\in X(v)$, $n\in X(w)$ where $\{u,v,w\}\subset X$
and $G[X \cup \{n-2,n-1\}]$ is indecomposable. Using Lemma \ref{lempmoins3}, we
may assume that $n-2\in X(1)$, $n-1\in X(2)$
 and $n-2\, .\ldots \, \{1,n-1\}$. Set  $Y=(X \setminus \{2\}) \cup \{n-1\}$. In this case,
$G[Y]\simeq P_{n-3}$ $n-2 \, .\ldots \, Y$ and $2\in Y(n-1)$. We may
then return to the first case.

\end{itemize}
{\hspace*{\fill}$\Box$\medskip}

\subsection{The class $\mathcal{Q}_{-3}$}
The next proposition describes the class $\mathcal{Q}_{-3}$.\\
We first introduce the class $\mathcal{G'}$ of graphs $G$ defined on
$\mathbb{N}_{n}$, where $n\geq 10$, such that $G[X]=Q_{n-3}$, $n-2
\in X^{-}(n-3)$ where $X=\mathbb{N}_{n-3}$, and satisfying one and
only one of the following assertions :
\begin{enumerate}
  \item $n-1 \in X^{-}, \, n \in X^{+}, \, \mid E(G_{X})\mid \geq 1$, and either $n-1 \not\sim \{n,n-2\}$ or
 $n \not\sim \{n-1,n-2\}$.

  \item $n-1,n \in X^{-}$ (resp. $n-1,n \in X^{+}), \, E(G_{X})= \{\{n-1,n-2\}\}$, and
$n-1 \sim \{n,n-2\}$.

  \item $n-1 \in X^{-}, \, n \in X(n-4)$ and $\mid E(G_{X}) \mid \geq 2$ and if $\{n,n-2\} \notin E(G_{X})$, then $\mid E(G_{X}) \mid = 2$.

 \item $n-1 \in X^{+}, \, n \in X(n-4)$ and $E(G_{X})= \{\{n-1,n-2\}, \{n-2,n\}\}$.

 \item $n-1 \in \langle X \rangle$, $n \in X^{-}(n-3), \, E(G_{X})= \{\{n-1,n-2\}\}$ and
 $n-2 \rule[0.4mm]{0.6cm}{0.1mm}\,\,n$.

 \item $n-1 \in X(n-4), \,n \in X^{+}(n-4)$ and either $n-1 .\ldots n$ and $E(G_{X})= \{\{n,n-2\}\}$, or $n-1 .\ldots n-4$ and $E(G_{X})= \{\{n-1,n-2\},\{n-2,n\}\}.$

 \item $n-1 \in X(1), \,n \in X(2)$ and either $E(G_{X})= \{\{n-1,n-2\},\{n-2,n\}\}$, or $\{\{n-1,n\},\{n-2,n\}\} \subseteq E(G_{X})$ with $n-1 \sim \{1,n\}.$

\item $n-1 \in X(u)$ where $u \in \{1,2\}, \,n \in X^{-}(n-3)$, $\,
E(G_{X})= \{\{n-1,n-2\}\}$ and
 $n-2 \rule[0.4mm]{0.6cm}{0.1mm}\,\,n$.

 \item $n-1 \in X^{+}(u), n \in X(u)$ where $u \in \{1,2\}$ and either $E(G_{X})= \{\{n-1,n-2\}\}$ and $n-1 .\ldots n$, or $E(G_{X})= \{\{n-1,n-2\},\{n-2,n\}\}$ and $n .\ldots u$.
\end{enumerate}

\vskip 0.9cm
\begin{proposition}
The indecomposability graph of graphs of the  class $\mathcal{G}'$
are $\{1, n-3\}$-covered.
\end{proposition}

\noindent{\em Proof\/. } Let  $G$ be a graph of the class
$\mathcal{G}'$. Set $X=\mathbb{N}_{n-3}$. Given  $i< j\in
\mathbb{N}_{n} \setminus \{1, n-3\}$, then $G-\{i, j\}$ is decomposable.
Indeed, if $3 \leq j\leq n-4$ and $n-1 \in \langle X \rangle \cup
X(n-4)$, then $\mathbb{N}_{i-1} \cup \{n-3\}$ is a non trivial
interval of $G-\{i,j\}$. If $3 \leq j\leq n-5$ (resp. $j=n-4$) and
$n-1 \in X(1) \cup X(2)$, then $(\mathbb{N}_{j-1} \cup
\{n-3,n-2,n-1,n\}) \setminus \{i\}$ (resp. $\mathbb{N}_{n} \setminus
\{i,j,n-3\}$) would be a non trivial interval of $G-\{i,j\}$.\\
If $j\geq n-2$ and $i \geq n-2$, then $G-\{i, j\}$ is
decomposable because $\rm{Ext}(\mathbb{N}_{n-3})=\emptyset$. \\
If $j\geq n-2$ and $ 2 \leq i \leq n-4$, we have to examine the
following cases.
\begin{itemize}
  \item $ 2 \leq i \leq n-4$ and $n-1 \in \langle X \rangle \cup
X(n-4)$.  If $i \neq n-4$, we get $\mathbb{N}_{i-1} \cup \{n-3\}$ is
a non trivial interval of $G-\{i,j\}$. If $i=n-4$, we distinguish
the following cases.
\begin{itemize}
  \item $n-1 \in \langle X \rangle$ and $n \notin X(n-4)$. Clearly, $\mathbb{N}_{i-1} \cup \{n-3\}$ is a non trivial
interval of $G-\{i,j\}$.
  \item $n \in X(n-4)$. If $n-1 .\ldots n-2$, then $\mathbb{N}_{n} \setminus
\{i,j,n-2\}$ would be a non trivial interval of $G-\{i,j\}$. If $n-1
\rule[0.4mm]{0.6cm}{0.1mm}\,\,n-2$ and $n-1 \in X^{-}$, then we have
to distinguish the following three cases according to $j$.
\begin{itemize}
  \item $j=n$. We verify that $\mathbb{N}_{n} \setminus
\{i,j,n-3\}$ would be a non trivial interval of $G-\{i,j\}$.
  \item $j=n-1$. If $n .\ldots n-2$, then $\mathbb{N}_{n} \setminus
\{i,j,n-2\}$ would be a non trivial interval of $G-\{i,j\}$. If $n
\rule[0.4mm]{0.6cm}{0.1mm}\,\,n-2$, then $\{n-3,n-2\}$ a non trivial
interval of $G-\{i,j\}$.
  \item $j=n-2$. We have either $\{n-3,n-1\}$ or $\mathbb{N}_{n} \setminus
\{i,j,n-1\}$ is a non trivial interval of $G-\{i,j\}$.
\end{itemize}

If $n-1 \rule[0.4mm]{0.6cm}{0.1mm}\,\,n-2$ and $n-1 \in X(n-4)$,
then we have to examine the following three cases according to $j$.
\begin{itemize}
  \item $j=n$. We verify that $\{n-3,n-2\}$ would be a non trivial interval of $G-\{i,j\}$.\\
  \item $j=n-1$. It is clear that $\mathbb{N}_{n} \setminus
\{i,j,n-2\}$ would be a non trivial interval of $G-\{i,j\}$.
  \item $j=n-2$. We get $\{n,n-1\}$ is a non trivial interval of $G-\{i,j\}$.\\
\end{itemize}
\end{itemize}
  \item $ 2 \leq i \leq n-4$ and $n-1 \in X(1) \cup X(2)$. If $i \neq 2$ then $(\mathbb{N}_{i-1} \cup
\{n-3,n-2,n-1,n\}) \setminus \{j\}$ would be a non trivial interval
of $G-\{i,j\}$. If $i=2$, we have to distinguish the following three
cases according to $X(2)$.
\begin{itemize}
  \item $X(2) = \emptyset$, we verify that $\{1,n-3,n-2,n-1,n\}
\setminus \{j\}$ is a non trivial interval of $G-\{i,j\}$.
  \item $\mid X(2) \mid = 1$. If $n \in X^{-}(n-3)$, then $\{1,n-2\}$ would be a
non trivial interval of $G-\{i,j\}$ if $j=n$. Moreover, $\{1,n-3\}$
would be a non trivial interval of $G-\{i,j\}$ if $j=n-1$. Finally,
$\{n,n-3\}$ is a non trivial interval of $G-\{i,j\}$ if $j=n-2$.\\

If $n \in X(2)$ and $n \rule[0.4mm]{0.6cm}{0.1mm}\,\,n-1$ (resp. $n
.\ldots n-1$), then $\{1,n-3,n-2,n-1\})$ would be a non trivial
interval of $G-\{i,j\}$ if $j=n$. Moreover, $\{1,n-2\}$ would be a
non trivial interval of $G-\{i,j\}$ if $j=n-1$. Finally, $\{1,n-1\}$
(resp. $\{n-1,n-3\}$) would be a non trivial interval of $G-\{i,j\}$
if $j=n-2$.

  \item $\mid X(2) \mid = 2$. If $n \rule[0.4mm]{0.6cm}{0.1mm}\,\,n-2$ (resp. $n
.\ldots n-2$), then $\{1,n-2\}$ (resp. $\{n-3,n-2\}$) would be a non
trivial interval of $G-\{i,j\}$ if $j=n-1$. Moreover, $\{1,n-2\}$
would be a non trivial interval of $G-\{i,j\}$ if $j=n$. Finally,
$\{n,n-1\}$ is a non trivial interval of $G-\{i,j\}$ if $j=n-2$.
\end{itemize}
\end{itemize}

It remains to verify that the graph $G$ is indecomposable. Since we
have either $\{n-1,n-2\} \in E(G_{X})$ or $\{n,n-2\} \in E(G_{X})$,
we prove by Theorem \ref{lemEhR}, that $G$ is indecomposable.\\
{\hspace*{\fill}$\Box$\medskip} \vskip 0.9cm

\begin{proposition} \label{propqmoins3}
Up to isomorphism, the graphs of the class $\mathcal{Q}_{-3}$ are
those of cardinality $\geq 10$ of the class $\mathcal{G'}$.
\end{proposition}

\vskip 0.3cm For the proof of this proposition, we need the next
results.

\begin{remark} \label{rqqmoins3}For $n \geq 10$ and $2 \leq i \leq n-4$, we introduce
the set $\mathbb{I}_{\mathbb{N},i}$ of non trivial intervals of
$Q_{n-3}-\{i\}$. We have :
\begin{enumerate}
  \item $\mathbb{I}_{\mathbb{N},2} = \{\{1,n-3\}\}$.
  \item $\mathbb{I}_{\mathbb{N},3} = \{\{1,2\}, \{1,2,n-3\}\}$.
  \item $\mathbb{I}_{\mathbb{N},4} = \{\{1,3\}, \{1,2,3\},
\{1,2,3,n-3\}\}$ for $n \geq 11$ and \\
 $\mathbb{I}_{\mathbb{N},4} = \{\{1,3\}, \{1,2,3\},
\{1,2,3,7\}, \{1,2,3,6,7\}\}$ for $n=10$.

 \item $\mathbb{I}_{\mathbb{N},5} = \{\{1,2,3,4\}, \{1,2,3,4,n-3\}$ for $n \neq 11$ and \\
 $\mathbb{I}_{\mathbb{N},5} = \{\{1,2,3,4\}, \{1,2,3,4,8\}, \{1,2,3,4,7,8\}\}$ for $n=11$.

 \item $\mathbb{I}_{\mathbb{N},n-4} = \{\mathbb{N}_{n-5}\}$.

 \item $\mathbb{I}_{\mathbb{N},i} = \{ \mathbb{N}_{i-1}, \mathbb{N}_{i-1} \cup
\{n-3\}\}$ for $6 \leq i \leq n-7$ or $i=n-5$.

 \item $\mathbb{I}_{\mathbb{N},n-6} = \{\mathbb{N}_{n-7}, \mathbb{N}_{n-7} \cup
\{n-3\}, \mathbb{N}_{n-7} \cup \{n-3,n-4\} \}$ for $n \geq 11$.

\end{enumerate}
\end{remark}

\begin{lemma}\label{lemqmoins3}
Let $G$ be an indecomposable graph defined on $\mathbb{N}_{n-1}$
where $n\geq 10$, verifying : $G[\mathbb{N}_{n-3}]=Q_{n-3}$ and for
each vertex $i$ of $\mathbb{N}_{n-1} \setminus \{1, n-3\}$, $i$ is critical.
Then one and only one of the following assertions holds:
\begin{enumerate}

\item $n-1\, .\ldots\, \mathbb{N}_{n-3}$ and $n-2\in
\mathbb{N}_{n-3}(u)$, where  $u\in\{n-4, n-3\}$. Moreover, if
$u=n-3$, then $n-2\, .\ldots\, n-3$.

\item $n-1 \rule[0.4mm]{0.6cm}{0.1mm}\,\, \mathbb{N}_{n-3}$ and $n-2\in
\mathbb{N}_{n-3}^{-}(n-3)$.

\item $n-1\in\mathbb{N}_{n-3}(u)$ and $n-2\in\mathbb{N}_{n-3}^{-}(v)$, where either $u \in \{1,2\}$ and $v=n-3$ or $u=2$ and $v=1$.

\end{enumerate}
\end{lemma}

\noindent{\em Proof\/. } Assume that $n-1 \, \ldots
\mathbb{N}_{n-3}$ (resp. $n-1 \rule[0.4mm]{0.6cm}{0.1mm}\,\,
\mathbb{N}_{n-3}$). As $n-1$ is a critical vertex, then $n-2 \notin
Ext(\mathbb{N}_{n-3})$. Besides, $n-2 \notin \langle
\mathbb{N}_{n-3} \rangle$ because $G$ is indecomposable. Thus, $n-2
\in \mathbb{N}_{n-3}(u)$ where $u \in \mathbb{N}_{n-3}$ and $n-2
\rule[0.4mm]{0.6cm}{0.1mm}\,\, n-1$ (resp. $n-1 .\ldots n-2$).
Therefore, we distinguish the following cases.
\begin{description}
  \item[Case $1$:] $n-1 \, \ldots \mathbb{N}_{n-3}$. \\
First, assume that $u=1$. We prove that $G-\{n-5\}$ is
indecomposable which is impossible. In fact, set $W=\{1, \ldots,
n-6\} \cup \{n-2,n-1\}$. As $G[\{1, \ldots, n-6\}]$ is
indecomposable, Theorem \ref{lemEhR} implies that $G[W]$ is
indecomposable. Moreover, $n-3 .\ldots W, \, n-4 \notin \langle
W\rangle$ (because $n-1 .\ldots n-4
\rule[0.4mm]{0.6cm}{0.1mm}\,\,1$) and $n-3
\rule[0.4mm]{0.6cm}{0.1mm}\,\, n-4$. So, $G[W \cup \{n-4,n-3\}]$ is
indecomposable.\\

Now, assume that $u=n-3$. We prove that $n-2 \ldots n-3$. Suppose
that it is not so. Clearly, $G[\{3,\ldots, n-3\}]$ is
indecomposable, and using Theorem \ref{lemEhR}, $G[\{3,\ldots, n-1\}]$
is indecomposable. Indeed, Set $W=\{3,\ldots,n-1\}$. It is easy to
verify that $1 \in Ext(W)$ so $G-\{2\}$ is
indecomposable; contradiction.\\

Finally, assume that $u \notin \{1,n-3\}$. Clearly,
$G[\mathbb{N}_{n-3} \cup \{n-2\} \setminus \{u\}] \, \simeq \,
Q_{n-3}$ and it is indecomposable. Set $Z = \mathbb{N}_{n-2}
\setminus \{u\}$. As $G-\{u\}$ is decomposable and $1 \,.\ldots n-1
\rule[0.4mm]{0.6cm}{0.1mm}\,\, n-2$, then $n-1 \in Z(v)$ where $v
\in Z$. Moreover, $v \neq n-2$ because otherwise $\{u,n-1,n-2\}$ is
a non trivial interval of $G$; contradiction. So, $N_{G[Z \cup
\{n-1\}]}(n-1) = \{n-2\}$ and $N_{G[Z]}(v) = \{n-2\}$. As, for each
$j \in \mathbb{N}_{n-3} \setminus \{n-5,n-3\}, \, \mid N_{Q_{n-3}}(j)\mid
\geq 2$, then $v \in \{n-5,n-3\}$. As more, $n-2
\rule[0.4mm]{0.6cm}{0.1mm}\,\, v, \,\, n-2 \in \mathbb{N}_{n-3}(u)$
then $u \rule[0.4mm]{0.6cm}{0.1mm}\,\, v$. Thus, $u=n-4$ when
$v=n-3$ and $u=n-6$ when $v=n-5$. Suppose
that $u=n-6$, we get as seen above that $G-\{n-5\}$ is indecomposable; impossible.\\

In summary, $u \in \{n-3,n-4\}$. Besides, if $u = n-3$, then $n-2
\ldots n-3$.

  \item[Case $2$:] $n-1 \rule[0.4mm]{0.6cm}{0.1mm}\,\,
\mathbb{N}_{n-3}$. \\
First assume that $u=1$. We prove that $G-\{n-4\}$ is
indecomposable; impossible. Indeed, by Theorem \ref{lemEhR}, $G[\mathbb{N}_{n-5} \cup \{n-2,n-1\}]$ is indecomposable.\\

Now assume that $u=n-3$. Suppose that $n-2
\rule[0.4mm]{0.6cm}{0.1mm}\,\, n-3$. Necessarily, $1 \notin \{3,
\ldots, n-1\}(n-3)$ and we obtain , as previously, that $G-\{2\}$ is
indecomposable; contradiction. Thus, if $u=n-3$ then $n-2 .\ldots n-3$.\\

Finally, assume that $u \notin \{1,n-3\}$. Clearly,
$G[\mathbb{N}_{n-2}] \setminus \{u\} \, \simeq Q_{n-3}$ and it is
indecomposable. Set $Z=\mathbb{N}_{n-2}-\{u\}$. As
$G[\mathbb{N}_{n-1}]-\{u\}$ is decomposable and $n-2 .\ldots n-1
\rule[0.4mm]{0.6cm}{0.1mm}\,\, 1$, then $n-1 \in Z(v)$ where $v \in
Z$. Moreover, $v \neq n-2$ because $G$ is indecomposable.
Consequently, $\overline{N}_{G[Z \cup \{n-1\}]}(n-1) = \{n-2\}$ and
$\overline{N}_{G[Z]}(v) = \{n-2\}$. Since we have for each $j \in
\mathbb{N}_{n-3} \setminus \{n-4\},\,\mid \overline{N}_{Q_{n-3}}(j)\mid \geq
2$, then $v=n-4$. As more, $n-2 .\ldots\, v, \,\, n-2 \in
\mathbb{N}_{n-3}(u)$ then $u .\ldots v$ and $u=n-5$. Suppose that
$u=n-5$, we prove in the same manner that $G-\{n-4\}$ is
indecomposable which is impossible.\\

 In summary, $u = n-3$, and $n-2 \ldots n-3$.

  \item[Case $3$:] $n-1 \in \mathbb{N}_{n-3}(u), \, n-2 \in \mathbb{N}_{n-3}(v)$ where $u \neq v \in \mathbb{N}_{n-3}$.\\
 Notice that if $v \notin \{1,n-3\}$ (resp. $u \notin \{1,n-3\}$),
then by interchanging $n-2$ with $v$ (resp. $n-1$ with $u$) and
building on the preceding cases, we may assume that $n-1 \in
(\mathbb{N}_{n-2} \setminus \{v\})(w)$ (resp. $n-2 \in ((\mathbb{N}_{n-3} \setminus \{u\}) \cup \{n-1\})(t)$) where $w \in \mathbb{N}_{n-2} \setminus \{u,v\}$ (resp. $t\in (\mathbb{N}_{n-3} \cup \{n-1\}) \setminus \{u,v\}$). So,
$\{u,w\}$ (resp. $\{v,t\}$) is an interval of $Q_{n-3}-v$ (resp.
$Q_{n-3}-u$. Then, using Remark \ref{rqqmoins3}, we can suppose that
$\{u,v\} \cap \{1,n-3\} \neq \emptyset$.\\

Assume first that $v=n-3$. It results from what precedes and Remark
\ref{rqqmoins3}, that $u \in \{1,2\}$. At present, showing that $n-3
.\ldots n-2$. Suppose that it is not so. As $n-1
\rule[0.4mm]{0.6cm}{0.1mm}\,\,n-2, \, n-2 \in \mathbb{N}_{n-5}^{-}$
and $n-1 \in \mathbb{N}_{n-5}(u)$, then Theorem \ref{lemEhR} claims
that $G[\mathbb{N}_{n-5} \cup \{n-1,n-2\}]$ is indecomposable. Set
$W=\mathbb{N}_{n-5} \cup \{n-1,n-2\}$. We have $n-3 \notin \langle W
\rangle$ (because $n-1 .\ldots n-3
\rule[0.4mm]{0.6cm}{0.1mm}\,\,n-2$) and for $x \in W, \, n-3 \notin
W(x)$, that means $n-3 \in Ext(W)$. So, $G-\{n-4\}$ is
indecomposable; impossible. Hence, $u \in \{1,2\}$, $v=n-3$ and $n-2 .\ldots n-3$.\\

Notice that by Remark \ref{rqqmoins3}, if $v=1$, then $u \in
\{n-3,2,3,4\}$. Suppose that $u \in \{3,4\}$. Necessarily, $1
\rule[0.4mm]{0.6cm}{0.1mm}\,\, n-2$. Otherwise, we verify that
$G-\{3\}$ is indecomposable; contradiction. But, $1
\rule[0.4mm]{0.6cm}{0.1mm}\,\, n-2$ gives that $G-\{2\}$ is
indecomposable which is also impossible. So, $u \in \{n-3,2\}$.\\
\hspace*{0.5cm}- If $u=2$, then we prove that $1 .\ldots n-2$.
Otherwise, we verify that $G[\{3, \ldots,n-3\} \cup \{1,n-1\}]$ is
indecomposable ($1 \in \{3, \ldots,n-3\}(n-3), \, n-1 \notin \{3,
\ldots,n-3\}(n-3)$ and $n-3 .\ldots n-1
\rule[0.4mm]{0.6cm}{0.1mm}\,\, 1$). Consequently, we get $n-2 \in
Ext(\{3, \ldots,n-3\} \cup \{1,n-1\})$, that means, $G-\{2\}$ is indecomposable; contradiction.\\
\hspace*{0.5cm}- If $u=n-3$, we may return to the previous case
where $u=1$ and $v=n-3$. \\

Hence, $u=2$, $v=1$ and  $n-2 .\ldots 1$.
\end{description}
{\hspace*{\fill}$\Box$\medskip} \vspace*{0.2cm}

\noindent{\em Proof of  Proposition \ref{propqmoins3}\/. } Let $G$
be a graph defined on $\mathbb{N}_{n}$ where  $n\geq 10$. Assume
that $\mathbb{I}(G)$ is $\{1,n-3\}$-covered and
$G[\mathbb{N}_{n-3}]=Q_{n-3}$. Notice that for $z\neq
t\in\{n-2,n-1,n\}$, if $G[\mathbb{N}_{n-3}\cup\{z,t\}]$ is
indecomposable, then  for each vertex
$i\in(\mathbb{N}_{n-3}\cup\{z,t\}) \setminus \{1,n-3\}$, $i$ is a
critical vertex of $G[\mathbb{N}_{n-3}\cup\{z,t\}]$. Even more
notice that from Corollary \ref{corolEhR}, there is  $x\neq
y\in\{n-2, n-1, n\}$ such that $G[\mathbb{N}_{n-3}\cup\{x,y\}]$ is
indecomposable. Without loss of generality, we may assume that
$\{x,y\}=\{n-1,n-2\}$. Using Lemma \ref{lemqmoins3}, we have to
distinguish the three following cases.
\begin{itemize}

\item \underline{Case 1: }$n-1 .\ldots \mathbb{N}_{n-3}$
and $n-2\in \mathbb{N}_{n-3}(u)$ where  $u\in \mathbb{N}_{n-3}$. In
this case, $n-1 \rule[0.4mm]{0.6cm}{0.1mm}\,\,n-2$. By Lemma
\ref{lemqmoins3}, we have $u\in \{n-3,n-4\}$. Moreover, if $n-2
.\ldots n-3$ then $u=n-3$. As $G-\{n-2,n-1\}$ is decomposable, we
examine the next cases.

\begin{itemize}
\item If $n \in \langle \mathbb{N}_{n-3}\rangle$. Assume that $n\, .\ldots \,
\mathbb{N}_{n-3}$. As $\{n,n-1\}$ is not an interval of $G$, $n
.\ldots n-2$. Moreover, since $\mathbb{N}_{n-1}$ is not an interval
of $G$, then $n \rule[0.4mm]{0.6cm}{0.1mm}\,\,n-1$. If $u=n-3$
(resp. $u=n-4$), $G \in \mathcal{G'}$ (resp. $G$ is isomorphic to
one of the elements of $\mathcal{G'}$ by permuting the vertices
$n-2$ and $n-4$). Assume now that $n
\rule[0.4mm]{0.6cm}{0.1mm}\,\,\mathbb{N}_{n-3}$. As
$\mathbb{N}_{n-1}$ is not an interval of $G$, then $n .\ldots n-1$
or $n .\ldots n-2$. If $u=n-3$, then $G \in \mathcal{G'}$. If
$u=n-4$, we verify that $n \rule[0.4mm]{0.6cm}{0.1mm}\,\,n-2$
(because $n .\ldots n-2$ implies that $G-\{n-4,n-1\}$ is
indecomposable). Besides, as $G$ is indecomposable, then $n .\ldots
n-1$ and by permuting $n-2$ and $n-4$, we get that $G$ is isomorphic
to one of the elements of $\mathcal{G'}$.

\item $n\in \mathbb{N}_{n-3}(v)$ where  $v\in\mathbb{N}_{n-3}$.
First, assume that $n \rule[0.4mm]{0.6cm}{0.1mm}\,\,n-1$. So,
$G[\mathbb{N}_{n-3} \cup \{n-1,n\}]$ is indecomposable, that is,
$G-\{n-2\}$ is indecomposable. By exchanging $n-2$ and $n$, and
applying Lemma \ref{lemqmoins3}, we obtain either $v=n-4$ or $v=n-3$
and $n-3 .\ldots n$.
\begin{itemize}
  \item If $u=n-3$ then $n-2 .\ldots n-3$. If $v=n-4$, $G \in
\mathcal{G'}$. If $v=n-3$ then $n-3 .\ldots n$, so $\{n,n-2\}$ is an
interval of $G$; impossible.
  \item If $u=n-4$. If $v=n-4$, we may assume that $n .\ldots n-4 \rule[0.4mm]{0.6cm}{0.1mm}\,\,n-2$
(because $\{n,n-2\}$ is not an interval of $G$), and by permuting
$n-2$ and $n-4$, $G$ is isomorphic to one of the elements of
$\mathcal{G'}$. If $v=n-3$ and $n-3 .\ldots n$, then $G$ is
isomorphic to one of the elements of $\mathcal{G'}$.
\end{itemize}

Then, assume that $n \ldots n-1$. We prove that $v \in \{n-3,n-4\}$.
Indeed, we have to distinguish the two cases.\\
$\hspace*{0.2cm} \triangleright \,\,\,u=n-3$ and $n-2 .\ldots n-3$.
Suppose that $v \notin \{n-3,n-4\}$. As $\{n,v\}$ is not an interval
of $G$, then $n \rule[0.4mm]{0.6cm}{0.1mm}\,\, n-2$. If $v=n-5$, we
prove that $G-\{n-5,n-1\}$ is indecomposable. If $v < n-5$, we prove
that $G-\{n-5,n-6\}$ is indecomposable.\\
$\hspace*{0.2cm} \triangleright \,\,\,u=n-4$. If $v=n-5$ (resp. $v <
n-5$), then $n \rule[0.4mm]{0.6cm}{0.1mm}\,\,n-2$ (resp. $n .\ldots
n-2$) because $\{n,v\}$ is not an interval of $G$. In this case,we
get $G-\{n-4,n-1\}$ is indecomposable; impossible.\\

Hence, $v \in \{n-3,n-4\}$.
\begin{itemize}
  \item If $u=n-3$. If $v=n-4$, then $n-2 .\ldots n$ (otherwise
$\{n,n-4\}$ is an interval of $G$; impossible), and $G \in
\mathcal{G'}$. If $v=n-3$, then $n-2 \rule[0.4mm]{0.6cm}{0.1mm}\,\,
n$ (otherwise $\{n,n-3\}$ is an interval of $G$; impossible).
Besides, we prove that $n-3 .\ldots n$ (otherwise $G-\{2,3\}$ is
indecomposable; contradiction). So, $G \in \mathcal{G'}$.
  \item If $u=n-4$. If $v=n-4$, we can assume
that $n .\ldots n-2 \rule[0.4mm]{0.6cm}{0.1mm}\,\,n-4$ (because
$\{n,n-4\}$ is not an interval of $G$). By permuting $n-2$ and
$n-4$, $G$ is isomorphic to one of the elements of $\mathcal{G'}$.
If $v=n-3$, then $n-2 .\ldots n$ (otherwise $\{n,n-3\}$ is an
interval of $G$;impossible). Besides, by Lemma \ref{lemqmoins3},
$n-3 .\ldots n$. Clearly, $G$ is isomorphic to one of the elements
of $\mathcal{G'}$.
\end{itemize}

\end{itemize}

\item \underline{Case 2: }$n-1 \rule[0.4mm]{0.6cm}{0.1mm}\,\, \mathbb{N}_{n-3}$
and $n-2\in \mathbb{N}_{n-3}(u)$ where  $u\in \mathbb{N}_{n-3}$.
Lemma \ref{lemqmoins3} implies that $u=n-3$ where $n-2 .\ldots n-3$.
As $G-\{n-2,n-1\}$ is decomposable, we examine the following cases.

\begin{itemize}
\item If $n \in \langle \mathbb{N}_{n-3}\rangle$. Assume that $n \rule[0.4mm]{0.6cm}{0.1mm}\,\,
\mathbb{N}_{n-3}$. As $\{n,n-1\}$ is not an interval of $G$, $n
\rule[0.4mm]{0.6cm}{0.1mm}\,\, n-2$. Moreover, since
$\mathbb{N}_{n-1}$ is not an interval of $G$, then $n .\ldots n-1$
and $G \in \mathcal{G'}$. Assume now that $n .\ldots
\,\,\mathbb{N}_{n-3}$. Necessarily, $n
\rule[0.4mm]{0.6cm}{0.1mm}\,\, n-1$ or $n-2
\rule[0.4mm]{0.6cm}{0.1mm}\,\,n$. In both cases, $G$ is isomorphic
to one of the elements of $\mathcal{G'}$.

\item $n\in \mathbb{N}_{n-3}(v)$ where  $v \in \mathbb{N}_{n-3}$.
First, we prove that $n \rule[0.4mm]{0.6cm}{0.1mm}\,\, n-1$. Suppose
that it is not so. We have $G[\mathbb{N}_{n-3} \cup \{n-1,n\}]$ is
indecomposable and using Lemma \ref{lemqmoins3}, we obtain that
$v=n-3$ and $n-3 .\ldots n$. Hence, $\{n,n-2\}$ is an interval of
$G$ which is impossible. Then, we show that $v \in \{n-4,n-3\}$.
Suppose that $v \leq n-5$, we verify that $G-\{n-4,n-5\}$ is
indecomposable. In fact, as $\{n,v\}$ is not an interval of $G$, $n
\rule[0.4mm]{0.6cm}{0.1mm}\,\, n-2$. Besides, $n-2 \in
\mathbb{N}_{n-6}^{-}, \, n \notin \langle \mathbb{N}_{n-6} \rangle$
and $n \rule[0.4mm]{0.6cm}{0.1mm}\,\, n-2$. So, $G[W]$ is
indecomposable where $W=\mathbb{N}_{n-6} \cup \{n-2,n\}$. Since $n-3
\in W^{-}, \, n-1 \notin \langle W\rangle$ and $n-3
\rule[0.4mm]{0.6cm}{0.1mm}\,\, n-1$, then $G[W \cup \{n-3,n-1\}]$ is
indecomposable. Therefore, we examine the two following cases:
\begin{itemize}
  \item If $v=n-4$, then $n .\ldots n-2$ (because $\{n,n-4\}$ is not an
interval of $G$), and $G \in \mathcal{G'}$.
  \item If $v=n-3$, then $n-2
\rule[0.4mm]{0.6cm}{0.1mm}\,\, n$ (otherwise $\{n,n-3\}$ is an
interval of $G$; impossible). Besides, we prove that $n-3 .\ldots n$
(otherwise $G-\{2,3\}$ is indecomposable; contradiction). Thus, $G
\in \mathcal{G'}$.
\end{itemize}

\end{itemize}

\item \underline{Case 3: }$ \{n-1,n-2\} \cap \langle\mathbb{N}_{n-3}\rangle=\emptyset$.
In this case, $n-1\in \mathbb{N}_{n-3}(u)$ and $n-2\in
\mathbb{N}_{n-3}(v)$ where $u \neq v\in\mathbb{N}_{n-3}$. By Lemma
\ref{lemqmoins3}, we have either $u \in \{1,2\}$ and $v=n-3$ with
$n-3 .\ldots n-2$, or $u=2$ and $v=1$ with $1 .\ldots n-2$. \\
Moreover, if $n \in \langle \mathbb{N}_{n-3}\rangle$, then either
$G-\{n-1\}$ or $G-\{n-2\}$ is indecomposable, which refers to one of
the first two cases. Thus, we can assume that $n \in
\mathbb{N}_{n-3}(\gamma)$ where $\gamma \in \mathbb{N}_{n-3}$.
\begin{itemize}
  \item $v=n-3$ and $u=1$. First, assume that $\gamma \notin
\{u,v\}$. As $\{n,\gamma\}$ is not an interval of $G$, then either
$G[\mathbb{N}_{n-3} \cup \{n,n-2\}]$ or $G[\mathbb{N}_{n-3} \cup
\{n,n-1\}]$ is indecomposable. Using Lemma \ref{lemqmoins3}, we may
assume that $\gamma=2$. In the first case, $n-2
\rule[0.4mm]{0.6cm}{0.1mm}\,\, \{n-1,n\}$. If $n
\rule[0.4mm]{0.6cm}{0.1mm}\,\,n-1, \, G  \in \mathcal{G'}$. If $n
.\ldots n-1$, we prove that $1 .\ldots n-1$. Suppose that it is not
so, it is easy to verify that $G-\{2,n-2\}$ is indecomposable;
impossible. In this case, $n-1 .\ldots \{1,n\}$ and $G  \in
\mathcal{G'}$. In the second case, we may assume that $n .\ldots
n-1$ and $n .\ldots n-2$. In the same manner, as $G-\{2,n-2\}$ is
decomposable, we obtain $1 .\ldots n-1$. By permuting $n$ and $2$,
$G$ is
isomorphic to one of the elements of $\mathcal{G'}$.\\

Presently, assume that $\gamma = v = n-3$. As $\{n-2,n\}$ is not an
interval of $G$, we have either $n .\ldots n-1$ or $n
\rule[0.4mm]{0.6cm}{0.1mm}\,\, n-3$. In the first case, $n-2
\rule[0.4mm]{0.6cm}{0.1mm}\,\, n$ (because $\{n-3,n\}$ is not an
interval of $G$). Besides, as $G-\{n-4,n-5\}$ is decomposable, $n
.\ldots n-3$ and $ G  \in \mathcal{G'}$. In the second case, we can
suppose that $n \rule[0.4mm]{0.6cm}{0.1mm}\,\, n-1$ and we verify
that $G-\{n-4,n-2\}$ is indecomposable; impossible. \\

Now, assume that $\gamma = u = 1$. As $\{n-1,n\}$ is not an interval
of $G$, we have either $n .\ldots n-2$ or $1 \not\sim \{n-1,n\}$ .
In the first case, as $\{1,n\}$ is not an interval of $G$, then by
interchanging $1$ with $n$, we may assume that $ n .\ldots \, n-1
\rule[0.4mm]{0.6cm}{0.1mm}\,\, 1$, and $G  \in \mathcal{G'}$. In the
second case, we may assume that $n-2 \rule[0.4mm]{0.6cm}{0.1mm}\,\,
n$, $ n .\ldots \, 1 \rule[0.4mm]{0.6cm}{0.1mm}\,\, n-1$ and $G  \in
\mathcal{G'}$.

  \item $v=n-3$ and $u=2$. First assume that $\gamma \neq \{u,v\}$.
As $\{n,\gamma\}$ is not an interval of $G$ then either
$G[\mathbb{N}_{n-3} \cup \{n,n-2\}]$ or $G[\mathbb{N}_{n-3} \cup
\{n,n-1\}]$ is indecomposable, and as previously seen, we may assume
that $\gamma = 1$. In the first case, $n
\rule[0.4mm]{0.6cm}{0.1mm}\,\, n-2$. As seen above, we have either
$n \rule[0.4mm]{0.6cm}{0.1mm}\,\, n-1$ and $G$ is isomorphic to one
of the elements of $\mathcal{G'}$, or $n .\ldots \{1,n-1\}$ and $G$
is isomorphic to one of the elements of $\mathcal{G'}$. In the
second case, we have $n .\ldots n-1$ and we can assume that $n
.\ldots n-2$. Besides, as $G-\{2,n-2\}$ is decomposable, we have $1
.\ldots n$. So, $n .\ldots \{1,n-1,n-2\}$ and $G$ is isomorphic to
one of the elements of $\mathcal{G'}$.\\

Presently, assume that $\gamma = n-3$. Since $\{n-3,n\}$ is not an
interval of $G$, we have either $n \rule[0.4mm]{0.6cm}{0.1mm}\,\,
n-1$ or $n \rule[0.4mm]{0.6cm}{0.1mm}\,\, n-2$. In the first case,
$G-\{n-4,n-2\}$ is decomposable implies that $n .\ldots n-3$ and
$\{n,n-2\}$ is an interval of $G$; contradiction. Assume now that
$n-1 .\ldots n$ and $n \rule[0.4mm]{0.6cm}{0.1mm}\,\, n-2$. As
$G-\{n-4,n-5\}$ is decomposable, then $n-3 .\ldots n$, and $G  \in
\mathcal{G'}$.\\

Consider the case where $\gamma = 2$. If $n
\rule[0.4mm]{0.6cm}{0.1mm}\,\, n-2$, then as $\{n-1,n\}$ is not an
interval of $G$, we may assume that $n .\ldots 2
\rule[0.4mm]{0.6cm}{0.1mm}\,\, n-1$. So, $ G  \in \mathcal{G'}$. If
$n-2 \ldots n$, by exchanging $2$ and $n$, we may assume that $n
.\ldots n-1 \rule[0.4mm]{0.6cm}{0.1mm}\,\, 2$ (because $\{2,n\}$ is
not an interval of $G$) and $G  \in \mathcal{G'}$.

  \item $v=1$ and $u=2$. In this case, $n-2 \in \mathbb{N}_{n-3}^{-}(1)$.
Set  $Y=(\mathbb{N}_{n-3} \setminus \{2\}) \cup \{n-1\}$. In this
case, $G[Y]\simeq Q_{n-3}$, $n-2 \in Y^{-}(n-3)$, $2 \in Y(n-1)$ and
$n-2 \rule[0.4mm]{0.6cm}{0.1mm}\,\, 2$. Then, by interchanging $n-1$
with $2$, we may then return to the previous case where $v=n-3$ and
$u=2$.

\end{itemize}
\end{itemize}

{\hspace*{\fill}$\Box$\medskip}

\subsection{The class $\mathcal{P}_{-5}$}
We describe the class $\mathcal{P}_{-5}$ by the following proposition.
\begin{proposition} \label{proppmoins5} Up to isomorphism, the elements of $\mathcal{P}_{-5}$ are the graphs $G$ defined on $\mathbb{N}_{n}$ where
$n\geq 14$, such that $G[X]=P_{n-5}$, $X=\mathbb{N}_{n-5}$,
satisfying that for each $Y \in q_{X}^{-}$ (resp. $Y
\in q_{X}^{+}$), $G[Y]$ is empty (resp. complete) and the graph $G_{X}$ is isomorphic to the bipartite graph $P_{5}$, with the bipartition $\{X^{-}, X^{-}(1)\}$ or $\{X^{-}, X^{+}(1)\}$.\\
\end{proposition}

\noindent{\em Proof\/. } Consider a graph $G=(V,E)$ satisfying the hypotheses of Proposition \ref{proppmoins5}. Assume, for instance, that the graph $G_{X}=(\{n-4,n-3,n-2,n-1,n\},E_{X})$ where either $E_{X} = \{ \{n-3,n-4\}, \{n-4,n-1\}, \{n-1,n-2\}, \{n-2,n\} \}\,\, when \,\,X(1) = \{n-3,n-1,n\}$ or $E_{X}=\{ \{n-4,n-3\}, \{n-3,n-2\}, \{n-2,n-1\}, \{n-1,n\}  \}\,\, when \,\, X(1) = \{n-3,n-1\}$.\\

\noindent We verify that $\mathbb{I}(G)$ is $\{1,n-5\}$-covered. Given  $i< j\in \mathbb{N}_{n} \setminus \{1, n-5\}$,
 we show that $G-\{i, j\}$ is decomposable.
% If $3 \leq j\leq n-7$, then $\{j+1, \ldots ,n-5\}$ is a non-trivial interval of
%$G-\{i,j\}$. In addition, if $j=n-6$, $\mathbb{N}_{n} \setminus
%\{i,j,n-5\}$ is a non-trivial interval of $G-\{i,j\}$. If $j\geq n-4$
If $2 \leq i \leq n-6$, it is clear that $(\mathbb{N}_{i-1} \cup
\{n-4,n-3,n-2,n-1,n\}) \setminus \{j\}$ is a non-trivial interval of
$G-\{i,j\}$. Assume that $i \geq n-4$. In this case, if $\{i,j\} \in q_{X}$, then $V \setminus (X \cup \{i,j\})$ is a non-trivial interval of $G-\{i,j\}$. Otherwise, $\{i,j\} \cap Y \neq \emptyset$ where $Y$ is the element of $q_{X}$ such that $\mid  Y\mid = 3$. Clearly, there is $k \in \{i,j\}$ such that $G_{X}-k$ is isomorphic to $P_{4}$. Thus, by Theorem \ref{ThBDIconnexe}, $G-k$ is critical according to $G[X]$ and then, $G-\{i,j\}$ is decomposable. It remains to prove that the graph $G$ is
indecomposable. Since $G[X]=P_{n-5}$ is indecomposable, then by Theorem \ref{ThBDIconnexe}, $G[\mathbb{N}_{n-1}]$ is critical according to $G[X]$. Consequently, $G[\mathbb{N}_{n-1}]$ is indecomposable, and
using Theorem \ref{lemEhR}, it is easy to
show that $n \in Ext(\mathbb{N}_{n-1})$.\\

Conversely, let $G$ be a graph defined on
$\mathbb{N}_{n}$ where  $n\geq 14$. Assume that $\mathbb{I}(G)$ is
$\{1,n-5\}$-covered and $G[\mathbb{N}_{n-5}]=P_{n-5}$. Set
$X=\mathbb{N}_{n-5}$. Observe that from Corollary \ref{corolEhR},
there are $x\neq y \in \{n-4,n-3,n-2,n-1,n\}$ such that $G[X \cup
\{x,y\}]$ is indecomposable. We may assume that $x=n-4$ and $y=n-3$.
Similarly, observe that there are $z\neq t \in \{n-2,n-1,n\}$ such
that $G[\mathbb{N}_{n-3} \cup \{z,t\}]$ is indecomposable. Assume,
for instance, that $z=n-2$ and $t=n-1$. Thus, $G[\mathbb{N}_{n-1}]$
is indecomposable. Set $H = G[\mathbb{N}_{n-1}]$. Since the
graph $\mathbb{I}(G)$ is $\{1,n-5\}$-covered, then for each $i \in
\mathbb{N}_{n-1} \setminus \{1,n-5\}$, $i$ is a critical vertex of
$H$. Thus, $H$ is critical according to $G[X]$. Consider the
graph $H_{X}$ defined on $\{n-4,n-3,n-2,n-1\}$. Corollary
\ref{corolBDI} claims that $H_{X}$ has no isolated vertex.
Moreover, we have either $H_{X}$ is connected and by Theorem
\ref{ThBDIconnexe}, $H_{X} \simeq P_{4}$ or $H_{X}$ has two
connected components and each of them is isomorphic to $P_{2}$. Without loss of
generality, we may assume that either $E(H_{X}) =
\{\{n-2,n-1\},\{n-1,n-4\},\{n-4,n-3\}\}$ or $E(H_{X}) =
\{\{n-2,n-1\},\{n-4,n-3\}\}$. In particular, $H[X \cup \{n-4,n-3\}]$ and
$H[X \cup \{n-2,n-1\}]$ are indecomposable.
%Set $W_{1} = H[X \cup \{n-4,n-3\}]$ and $W_{2} = H[X \cup \{n-2,n-1\}]$.
Notice that if $\{\alpha,\beta\} \in E(H_{X})$, then for
each vertex $i\in (X \cup \{\alpha,\beta\}) \setminus \{1,n-5\}$, $i$ is a critical vertex of $H[X \cup \{\alpha,\beta\}]$. Otherwise, there exists a non critical vertex $i\in (X \cup \{\alpha,\beta\}) \setminus \{1,n-5\}$ of $H[X \cup
\{\alpha,\beta\}]$, that is, $H[(X \cup \{\alpha,\beta\}) \setminus
\{i\}]$ is indecomposable. Set $Z = (X \cup \{\alpha,\beta\})
\setminus \{i\}$. So, $G[Z]$ is indecomposable and $\mid V \setminus
Z\mid = 4$. It follows from Corollary \ref{corolEhR} that there is $\{\gamma,\delta\}
\subseteq V \setminus Z$ such that $G-\{\gamma,\delta\}$ is indecomposable
which contradicts the fact that $\{\gamma,\delta\} \cap \{1,n-5\} =
\emptyset$. \\
%Consequently, we can apply successively Lemma \ref{lempmoins3} to both graphs $W_{1}$ and $W_{2}$. \\
Consequently, we may assume using Lemma \ref{lempmoins3}, that $\alpha \in X^{-}$ and $\beta \in X(\mu)$ where
$\mu \in \{1,2\}$. Indeed, if $\alpha \in X^{-}$ and $\beta \in X(n-6) \cup X(n-5)$, then by considering the bijection $f$ defined on $\mathbb{N}_{n-5} \cup \{\alpha,\beta\}$ by: for each $i \in \mathbb{N}_{n-5}$, $f(i) = n-i-4, \, f(\alpha)
= \alpha$ and $f(\beta) = \beta$, we can assume that $\alpha \in X^{-}$ and $\beta \in X(\mu)$ where $\mu \in \{1,2\}$. If $\alpha \in X(1)$ and $\beta \in X(2)$, it suffices to permute the vertices $\beta$ and $2$. Finally, if $\alpha \in X(n-6)$ and $\beta \in X(n-5)$, we may easily return to the first case by permuting the vertices $n-6$ and $\alpha$. \\
Now, we pose $Y = X \cup \{n-4,n-3\}$ and we distinguish the following cases according to $E(H_{X})$.
\begin{itemize}
\item \underline{Case 1: } $E(H_{X}) =
\{\{n-2,n-1\},\{n-1,n-4\},\{n-4,n-3\}\}$. We may assume, up to isomorphism, that $n-2 \in X^{-}$ and $n-1 \in X(\mu)$ where $\mu \in \{1,2\}$. By Theorem \ref{ThBDIconnexe}, we have to examine the two cases.
\begin{description}
\item [Case $1.1$ :] If $\{n-2,n-4\} \subseteq X^{-}$ and $\{n-1,n-3\} \subseteq
X(1)$. By Theorem \ref{ThBDIconnexe}, $n-2 .\ldots \{n-4,n-3\}$. Clearly, as $G-\{n-2,n-1\}$ is decomposable, then $n \notin
Ext(Y)$. We distinguish the two following cases according to $n$.

\begin{description}
\item [$\hookrightarrow$] $ n \in \langle Y \rangle$. If $n .\ldots Y$. Since $n-1 \rule[0.4mm]{0.6cm}{0.1mm}\,\, n-2
.\ldots Y$ and $\{n,n-2\}$ is not an interval of $G$, then $n
.\ldots n-1$. Moreover, $\mathbb{N}_{n-1}$ is not an interval of
$G$, then $n \rule[0.4mm]{0.6cm}{0.1mm}\, n-2$, so that $G-\{n-3,n-4\}$ is indecomposable and $\mathbb{I}(G)$ would not be $\{1,n-5\}$-covered.\\
If $n \rule[0.4mm]{0.6cm}{0.1mm}\,\,Y$, then either $n .\ldots n-2$
or $n .\ldots n-1$ and we verify that $G-\{n-3,n-4\}$ is
indecomposable.

\item [$\hookrightarrow$] $n \in Y(u)$ where $u \in Y$. Since
$\{n,u\}$ is not an interval of $G$,
we have either $n-2 \not\sim \{n,u\}$ or $n-1 \not\sim \{n,u\}$.\\

$\ast \, u \in X.$ If $n-2 \not\sim \{n,u\}$, then
$n-2 \rule[0.4mm]{0.6cm}{0.1mm}\,\,n$. Moreover, as $G-\{n-3,n-4\}$ is decomposable, we get $u = 1$. In this case, we verify that $1
\sim \{n,n-1\}$ (because $G-\{n-3,n-4\}$ is decomposable) and $n-1
\sim \{1,n\}$ (because $G-\{n-3,n-2\}$ is decomposable). Thus, $G$ is one of the graphs defined in Proposition \ref{proppmoins5}.\\
If $n-1 \not\sim \{n,u\}$, we may assume that $n-2 \sim
\{n,u\}$, which means $n-2
.\ldots n$. Similarly, we prove that $G-\{n-3,n-4\}$ is indecomposable.\\

$\ast \, u = n-4.$  We have $n \rule[0.4mm]{0.6cm}{0.1mm}\,\, n-3$.
Moreover, $n-2 \not\sim \{n,n-4\}$ or $n-1 \not\sim \{n,n-4\}$. If
$n-2 \not\sim \{n,n-4\}$, that is $n-2
\rule[0.4mm]{0.6cm}{0.1mm}\,\, n$, then we prove that
$G-\{n-1,n-4\}$ is indecomposable. If $n-1 \not\sim
\{n,n-4\}$, that is $n .\ldots n-1$ and we can assume that $n
.\ldots n-2$. As $G-\{n-3,n-2\}$ is decomposable, we verify
that $n.\ldots n-4$ and so $G$ is isomorphic to one of the graphs defined in Proposition \ref{proppmoins5}.\\

$\ast \, u = n-3.$ We have $n-4 \rule[0.4mm]{0.6cm}{0.1mm}\,\, n$
and $1 \sim \{n-3,n\}$. If $n-2 \not\sim \{n,n-3\}$, that is $n
\rule[0.4mm]{0.6cm}{0.1mm}\,\,n-2$. Suppose that $n-3 \sim
\{n,n-1\}$, then $\{n-1,n\}$ is a non-trivial interval of $G$. Besides, suppose that $n-3 \not\sim \{n,n-1\}$, we prove
that $G-\{n-1,n-4\}$ is indecomposable. If $n-1
\not\sim \{n,n-3\}$, then we prove
that $G-\{n-3,n-4\}$ is indecomposable.
\end{description}

\item [Case $1.2$ :] If $\{n-2,n-4\} \subseteq X^{-}$ and $\{n-1,n-3\} \subseteq
X(2)$. We proceed in the same way as previously seen if $n \in
\langle Y \rangle$. If $n \in Y(u), \, u \in Y$, we have to examine the following cases.\\
$\ast \, u \in X.$ If $n-2 \not\sim \{n,u\}$ and $u \neq 2$ (resp.
$u = 2$), then we prove that $G-\{n-3,n-4\}$ is indecomposable
(resp. $G-\{2,n-1\}$ is indecomposable). If $n-1 \not\sim \{n,u\}$, we may assume that $n-2 .\ldots n$ and $G-\{n-3,n-4\}$ would be indecomposable.\\

$\ast \, u = n-4.$ If $n-2 \not\sim \{n,n-4\}$, then $n-2
\rule[0.4mm]{0.6cm}{0.1mm}\,\,n$ and $G-\{n-1,n-4\}$
is indecomposable. If
$n-1 \not\sim \{n,n-4\}$, we may assume that $n-2 .\ldots n$ and we verify that $G-\{2,n-4\}$ is indecomposable.\\

$\ast \, u = n-3.$ If $n-2 \not\sim \{n,n-3\}$, then necessarily
$n-3 \not\sim \{n,n-1\}$. Otherwise, $\{n-1,n\}$ would be a non-trivial
interval of $G$. In this case, we prove that
$G-\{n-1,n-4\}$ is indecomposable. If $n-1 \not\sim
\{n,n-3\}$, then we prove
that $G-\{n-3,n-4\}$ is indecomposable.

\end{description}

\item \underline{Case 2: } $E(H_{X}) =
\{\{n-2,n-1\}, \{n-4,n-3\}\}$. We can assume that $n-2 \in X^{-}$ and $n-1 \in X(\mu)$ where $\mu \in \{1,2\}$. We shall examine the following cases.
\begin{description}
  \item [Case $2.1$ :] $n-4 \in X^{-}, \, n-3 \in X(\nu)$ where $\nu \in
\{1,2,n-6,n-5\}$.
\begin{itemize}
  \item If $\mu=\nu=1$, then $\{n-2,n-4\} \subseteq X^{-}$ and $\{n-1,n-3\} \subseteq
X(1)$. As seen in Case $1$, $n \in Y(u)$ where $u \in Y$ and we have to distinguish the following cases.\\

$\ast \, $ If $u \in X$ and $n-2 \not\sim \{n,u\}$, we prove that
$G-\{n-3,n-4\}$ is indecomposable if $u \neq 1$ or $u = 1$ and $1
\not\sim \{n,n-1\}$. Therefore, $u = 1$ and $1 \sim \{n,n-1\}$, that
means $\{n-1,n\}$ would be a non-trivial interval of $G$. If $u \in X$ and $n-1 \not\sim \{n,u\}$, then we prove that $G-\{n-3,n-4\}$ is indecomposable.\\

$\ast \, u = n-4.$ If $n-2 \not\sim \{n,n-4\}$, then $n-2
\rule[0.4mm]{0.6cm}{0.1mm}\,\,n$ and $G-\{n-1,n-4\}$ would be indecomposable. If $n-1 \not\sim \{n,n-4\}$, then
$n-1 \rule[0.4mm]{0.6cm}{0.1mm}\,\,n$ and we can assume that $n-2
.\ldots n$. Suppose that $n \rule[0.4mm]{0.6cm}{0.1mm}\,\, n-4$, we
prove that $G-\{n-2,n-3\}$ is indecomposable. Hence,
$n .\ldots n-4$ and $G$ is isomorphic to one of the graphs defined in Proposition \ref{proppmoins5}.\\

$\ast \, u = n-3.$ We have $1 \sim \{n,n-3\}$. If $n-2 \not\sim
\{n,n-3\}$, then $n-2 \rule[0.4mm]{0.6cm}{0.1mm}\,\,n$. Suppose that
$n \not\sim \{1,n-3\}$, we verify that $G-\{n-1,n-4\}$ is
indecomposable. Moreover, as $G-\{n-3,n-2\}$ is
decomposable, then $n \sim \{1,n-1\}$. Clearly, $G$ is isomorphic to
one of the graphs defined in Proposition \ref{proppmoins5}. If $n-1 \not\sim
\{n,n-3\}$, then we prove that $G-\{n-3,n-4\}$ is indecomposable.\\

 \item If $\mu=\nu=2$, then we demonstrate easily that $2$ is a non-critical vertex
of $H$ which contradicts the fact that $\mathbb{I}(G)$ is $\{1,n-5\}$-covered. \\
%In fact, $H[(X \setminus \{2\}) \cup
%\{n-3\}] \simeq P_{n-5}$ is indecomposable. We pose $Z = (X
%\setminus \{2\}) \cup \{n-3\}$. Since $n-1 \in Z(n-3), \, n-4 \in
%Z^{-}(1)$ and $n-4 .\ldots \,n-1$, then $H[Z \cup \{n-1,n-4\}]$ is
%indecomposable. Set $W = Z \cup \{n-1,n-4\}$. It is easy to verify
%that $n-2 \in Ext(W)$, so $H[W \cup \{n-2\}]$ is
%indecomposable.\\

 \item If $\mu=1$ and $\nu =2$, then $\{n-2,n-4\} \subseteq X^{-}$, $n-3 \in X(2)$ and $n-1 \in X(1)$. Since $2$ is a critical
vertex of $H$, then necessarily, $n-1 \in X^{-}(1)$. Similarly, we have $n \in Y(u)$ where $u \in Y$ and we distinguish the following cases.\\

$\ast \, $ If $u \in X$ and $n-2 \not\sim \{n,u\}$, we prove that
$G-\{n-3,n-4\}$ is indecomposable if $u \neq 1$ or $u = 1$ and $1
\rule[0.4mm]{0.6cm}{0.1mm}\,\, n$. Therefore, $u = 1$ and $1 .\ldots
n$, that means $\{n-1,n\}$ would be a non-trivial interval of $G$. If $u \in X$ and $n-1 \not\sim \{n,u\}$, then we may assume that $n-2 .\ldots n$ and we prove that $G-\{n-3,n-4\}$ is indecomposable.\\

$\ast \, u = n-4.$ If $n-2 \not\sim \{n,n-4\}$, then we prove that
$G-\{n-1,n-4\}$ is indecomposable. If
$n-1 \not\sim \{n,n-4\}$, then we prove that $G-\{n-2,n-4\}$ is indecomposable.\\

$\ast \, u = n-3.$ We have either $n \rule[0.4mm]{0.6cm}{0.1mm}\,\, n-2$ or $n .\ldots n-1$. In both cases, we obtain that $G-\{n-3,n-4\}$ is indecomposable.\\

 \item If $\mu=1$ and $\nu=n-5 \,\, ($resp. $\mu=2$ and $\nu \in \{n-6,n-5\}$), then we
have $\{n-2,n-4\} \subseteq X^{-}, \, n-3 \in X(\nu)$ and $n-1 \in
X(\mu)$. Similarly, $n \in Y(u), \, u \in Y$ and we distinguish the following cases.\\
$\ast \, u \in X \setminus \{1\} \, $ (resp. $u \in X \setminus
\{2\}$), we verify that $G-\{n-3,n-4\}$ is indecomposable.\\
$\ast \, u=1 \,$ (resp. $u=2$). If $n-2 \not\sim \{n,u\}$, then $n-2
\rule[0.4mm]{0.6cm}{0.1mm}\,\, n$ and as $G-\{n-3,n-4\}$ is
decomposable, we get $u \sim \{n-1,n\}$. Therefore,
$\{n,n-1\}$ would be a non-trivial interval of $G$.\\
If $n-1 \not\sim \{n,u\}$, we can assume that
$n .\ldots n-2$ and then, $G-\{n-3,n-4\}$ would be indecomposable.\\

$\ast \, u = n-4.$ If $n-2 \not\sim \{n,n-4\}$, then $n-2
\rule[0.4mm]{0.6cm}{0.1mm}\,\, n$ and we prove that $G-\{n-1,n-4\}$
is indecomposable. If $n-1 \not\sim \{n,n-4\}$, then
$n-1 \rule[0.4mm]{0.6cm}{0.1mm}\,\, n$ and we prove
that $G-\{n-2,n-4\}$ is indecomposable.\\

$\ast \, u = n-3.$ We have $n \rule[0.4mm]{0.6cm}{0.1mm}\,\, n-2$ or
$n \rule[0.4mm]{0.6cm}{0.1mm}\,\, n-1$ and we get $G-\{n-3,n-4\}$ is
indecomposable.

\item If $\mu=1$ and $\nu =n-6$, then we may return to the case where $\mu=2$ and $\nu=n-5$.

\end{itemize}

  \item [Case $2.2$ :] $n-4 \in X^{-}(1)$ and $n-3 \in X(2)$.
\begin{itemize}
  \item If $\mu=1$, then using Theorem \ref{ThBDInonconnexe}, we get
$n-1 \in X^{-}(1)$ and $G[\{1,n-1,n-4\}]$ is empty. Moreover, by
permuting $2$ and $n-3$, we may return to one of the previous cases.
  \item If $\mu=2$, it suffices to permute the vertices $2$ and $n-3$.
\end{itemize}

 \item [Case $2.3$ :] $n-4 \in X^{-}(n-5)$ and $n-3 \in X(n-6)$. In both cases where $\mu \in \{1,2\}$, it suffices to permute the vertices
 $n-6$ and $n-3$.

\end{description}
\end{itemize}
{\hspace*{\fill}$\Box$\medskip}

\subsection{The class $\mathcal{Q}_{-5}$}
\noindent The next proposition describes the class $\mathcal{Q}_{-5}$.

\begin{proposition} \label{propqmoins5} Up to isomorphism, the elements of $\mathcal{Q}_{-5}$ are the graphs $G$ defined on
$\mathbb{N}_{n}$ where $n \geq 12$, such that $G[X]=Q_{n-5}$,
$X=\mathbb{N}_{n-5}$, satisfying that for each $Y \in q_{X}^{-}$
(resp. $Y \in q_{X}^{+}$), $G[Y]$ is empty (resp. complete) and the graph $G_{X}$ is isomorphic to the bipartite graph $P_{5}$, with the bipartition $\{X^{-}(n-5), Z\}$ where $Z \in \{X^{-}, X^{+}, X^{-}(1), X^{+}(1)\}$.\\
\end{proposition}

\noindent{\em Proof\/. } Consider a graph $G=(V,E)$ satisfying the hypotheses of Proposition \ref{propqmoins5}. Assume, for instance, that the graph $G_{X}=(\{n-4,n-3,n-2,n-1,n\},E_{X})$ where either $E_{X}= \{ \{n-3,n-4\}, \{n-4,n-1\}, \{n-1,n-2\}, \{n-2,n\} \}$ if $X^{-}(n-5) = \{n-4,n-2\}$ or $E_{X}=\{ \{n-4,n-3\}, \{n-3,n-2\}, \{n-2,n-1\}, \{n-1,n\}  \}$ if $X^{-}(n-5) = \{n-4,n-2,n\}$.\\

\noindent We verify that $\mathbb{I}(G)$ is $\{1,n-5\}$-covered. Indeed, given  $i< j\in
\mathbb{N}_{n} \setminus \{1, n-5\}$, we prove that $G-\{i, j\}$ is
decomposable. If $3 \leq j\leq n-6$ and $\langle X \rangle
\neq \emptyset$, then $\mathbb{N}_{i-1} \cup \{n-5\}$ is a non
trivial interval of $G-\{i,j\}$. If $3 \leq j\leq n-7$ (resp.
$j=n-6$) and $\langle X \rangle = \emptyset$, then
$(\mathbb{N}_{j-1} \cup \{n-5,n-4,n-3,n-2,n-1,n\}) \setminus \{i\}$
(resp. $\mathbb{N}_{n} \setminus \{i,j,n-5\}$) would be a non
trivial interval of $G-\{i,j\}$. It is obviously the case if $j\geq
n-4$ and $i \geq n-4$. If $j\geq n-4$ and $ 2 \leq i \leq n-6$, we
have to examine the two following cases.
\begin{itemize}
  \item If $\langle X \rangle \neq \emptyset$, then $\mathbb{N}_{i-1} \cup \{n-5\}$ is a non
trivial interval of $G-\{i,j\}$.
  \item If $\langle X \rangle = \emptyset$ and $i \neq 2$ then
$(\mathbb{N}_{i-1} \cup \{n-4,n-3,n-2,n-1,n\}) \setminus \{j\}$
would be a non trivial interval of $G-\{i,j\}$. Otherwise,
 $\{1,n-5,n-4,n-3,n-2,n-1,n\}) \setminus \{j\}$ is a non trivial interval of
$G-\{i,j\}$.
\end{itemize}

Now, we show that $G$ is an indecomposable graph.
As previously seen, we get by Theorem \ref{ThBDIconnexe} that
$G[\mathbb{N}_{n-1}]$ is indecomposable and we verify using Theorem
\ref{lemEhR}, that $n \in Ext(\mathbb{N}_{n-1})$.\\

Conversely, let  $G$ be a graph defined on
$\mathbb{N}_{n}$ where  $n\geq 12$. Assume that $\mathbb{I}(G)$ is
$\{1,n-5\}$-covered and $G[\mathbb{N}_{n-5}]=Q_{n-5}$. We pose
$X=\mathbb{N}_{n-5}$. As previously seen, notice that from Corollary \ref{corolEhR},
there are $x\neq y \in \{n-4,n-3,n-2,n-1,n\}$ such that $G[X \cup
\{x,y\}]$ is indecomposable. We may assume that $x=n-4$ and $y=n-3$.
Similarly, notice that there are $z\neq t \in \{n-2,n-1,n\}$ such
that $G[\mathbb{N}_{n-3} \cup \{z,t\}]$ is indecomposable. Assume,
for instance, that $z=n-2$ and $t=n-1$. Thus, $G[\mathbb{N}_{n-1}]$
is indecomposable. We pose $H = G[\mathbb{N}_{n-1}]$. Since the
graph $\mathbb{I}(G)$ is $\{1,n-5\}$-covered, then for each $i \in
\mathbb{N}_{n-1} \setminus \{1,n-5\}$, $i$ is a critical vertex of
$H$. So, $H$ is critical according to $G[X]$. Now, consider the
graph $H_{X}$ defined on $\{n-4,n-3,n-2,n-1\}$. Corollary
\ref{corolBDI} claims that $H_{X}$ has no isolated vertex.
Moreover, we have either $H_{X}$ is connected and by Theorem
\ref{ThBDIconnexe}, $H_{X} \simeq P_{4}$ or $H_{X}$ has two
connected component and each of them is isomorphic to $P_{2}$. Without loss of
generality, we may assume that either $E(H_{X}) =
\{\{n-2,n-1\},\{n-1,n-4\},\{n-4,n-3\}\}$ or $E(H_{X}) =
\{\{n-2,n-1\},\{n-4,n-3\}\}$. We have $H[X \cup \{n-4,n-3\}]$ and
$H[X \cup \{n-2,n-1\}]$ are indecomposable. Set $W_{1} = H[X \cup
\{n-4,n-3\}]$ and $W_{2} = H[X \cup \{n-2,n-1\}]$. Even more, notice
that for $\alpha \neq \beta \, \in \{n-4,n-3,n-2,n-1\}$, if $H[X \cup \{\alpha,\beta\}]$ is indecomposable, then for
each vertex $i\in (X \cup \{\alpha,\beta\}) \setminus \{1,n-5\}$,
$i$ is a critical vertex of $H[X \cup \{\alpha,\beta\}]$. Similarly to the
proof of Proposition \ref{proppmoins5}, we apply successively
Lemma \ref{lemqmoins3} to both of the graphs $W_{1}$ and $W_{2}$. We
pose $Y = X \cup \{n-4,n-3\}$ and we distinguish the following cases
according to $E(H_{X})$.
\begin{itemize}
\item \underline{Case 1: } $E(H_{X}) =
\{\{n-2,n-1\},\{n-1,n-4\},\{n-4,n-3\}\}$. Using Theorem
\ref{ThBDIconnexe} and Lemma \ref{lemqmoins3}, and up to
isomorphism, we may examine the three cases.
\begin{description}
\item [Case $1.1$ :] $n-1 \in X(\mu)$ where $\mu \in \{1,2\}$ and $n-2 \in X^{-}(n-5)$.
\begin{itemize}
  \item If $\mu =1$, then $\{n-2,n-4\}  \subseteq  X^{-}(n-5)$ and $\{n-1,n-3\}  \subseteq
X(1)$ and by Theorem \ref{ThBDIconnexe}, $G[\{n-1,n-3,1\}]$ is empty
or complete. Clearly, as $G-\{n-2,n-1\}$ is decomposable, then $n
\notin Ext(Y)$. We shall examine the two following cases according
to $n$.

\begin{description}
\item [$\hookrightarrow$] $ n \in \langle Y \rangle$. We simply verify that $G-\{n-3,n-4\}$ is indecomposable; contradiction.

\item [$\hookrightarrow$] $n \in Y(u)$ where $u \in Y$. Notice that since
$\{n,u\}$ is not an interval of $G$, we have either $n-2 \not\sim
\{n,u\}$ or $n-1 \not\sim \{n,u\}$. We distinguish the following
cases.
  $\hspace*{0.2cm} \checkmark \, u=n-5$. If $n-1 \not\sim \{n-5,n\}$, then $G-\{n-3,n-4\}$ is
decomposable gives that $n .\ldots n-5$. Consequently, $\{n,n-2\}$
is a non trivial interval of $G$; impossible. If $n-2 \not\sim
\{n-5,n\}$, then $n \rule[0.4mm]{0.6cm}{0.1mm}\,\,n-2$ and we may
assume that $n .\ldots \,n-1$. We demonstrate that $G-\{n-3,n-4\}$
is indecomposable; contradiction.

  $\hspace*{0.2cm} \checkmark \, u=1$. If $n-1 \not\sim \{1,n\}$, necessarily $n \rule[0.4mm]{0.6cm}{0.1mm}\,\,
n-2$ and $1 \sim \{n-1,n\}$ (otherwise $G-\{n-3,n-4\}$ is
indecomposable; contradiction). But in this case, $G-\{n-3,n-2\}$ is
indecomposable; impossible. If $n-2 \not\sim \{1,n\}$, then $n
\rule[0.4mm]{0.6cm}{0.1mm}\,\, n-2$ and we can assume that $n-1 \sim
\{1,n\}$. Moreover, as $G-\{n-3,n-4\}$ is decomposable, then $1 \sim
\{n-1,n\}$. Thus, $G$ is one of the graphs defined in Proposition \ref{propqmoins5}.

  $\hspace*{0.2cm} \checkmark \, u \in X \setminus \{1,n-5\}$. Clearly, $G-\{n-3,n-4\}$ is indecomposable; contradiction.

  $\hspace*{0.2cm} \checkmark \, u=n-3$. If $n-1 \not\sim \{n-3,n\}$, then $n-1 \not\sim
\{n,1\}$. Besides, $G-\{n-3,n-4\}$ is decomposable gives that $n
\rule[0.4mm]{0.6cm}{0.1mm}\,\,n-2$ and $G-\{n-1,n-4\}$ is
decomposable implies that $n \sim \{n-3,1\}$. So, $n-3 \sim
\{n,n-1\}$ and $\{n,n-1\}$ is a non trivial interval of $G$;
contradiction. If $n-2 \not\sim \{n-3,n\}$, then , $n
\rule[0.4mm]{0.6cm}{0.1mm}\,\,n-2$ and we can assume that $n-1 \sim
\{n,n-3\}$. We get as previously, $n-3 \sim \{n,n-1\}$ and thus,
$\{n,n-1\}$ is a non trivial interval of $G$ which is impossible.

  $\hspace*{0.2cm} \checkmark \, u=n-4$. If $n-1 \not\sim \{n-4,n\}$, then $n-1 .\ldots n$. As $G-\{n-1,n-4\}$ is decomposable, we have $n
.\ldots n-2$. Besides, $G-\{n-2,n-3\}$ is decomposable implies that
$n .\ldots n-4$. Hence, $n .\ldots \{n-2,n-4\}$ and $G$ is isomorphic to one of the graphs defined in Proposition \ref{propqmoins5}. If $n-2
\not\sim \{n-4,n\}$, then $n \rule[0.4mm]{0.6cm}{0.1mm}\,\, n-2$ and
we may assume that $n-1 \sim \{n-4,n\}$. We obtain that
$G-\{n-1,n-4\}$ is indecomposable; contradiction.\\

\end{description}

  \item If $\mu =2$, then $\{n-2,n-4\}  \subseteq  X^{-}(n-5)$ and $\{n-1,n-3\}  \subseteq
X(2)$ and by Theorem \ref{ThBDIconnexe}, $G[\{n-1,n-3,2\}]$ is empty
or complete. Similarly, we examine the two following cases according
to $n$.

\begin{description}
\item [$\hookrightarrow$] $ n \in \langle Y \rangle$. We verify that $G-\{n-3,n-4\}$ is indecomposable; contradiction.

\item [$\hookrightarrow$] $n \in Y(u)$ where $u \in Y$. We distinguish the following cases.

  $\hspace*{0.2cm} \checkmark \, u=2$. If $n-1 \not\sim \{2,n\}$, then $G-\{n-3,n-2\}$ is indecomposable; contradiction.
If $n-2 \not\sim \{2,n\}$, then $n \rule[0.4mm]{0.6cm}{0.1mm}\,\,
n-2$ and $n-1 \sim \{2,n\}$. Moreover, as $G-\{n-3,n-4\}$ is
decomposable, then $2 \sim \{n-1,n\}$. But, in this case, we prove
that $G-\{2,n-1\}$ is indecomposable; impossible.

  $\hspace*{0.2cm} \checkmark \, u=n-5$. If $n-1 \not\sim \{n-5,n\}$, then $n-1 \rule[0.4mm]{0.6cm}{0.1mm}\,\,
n$. As $G-\{n-3,n-4\}$ is decomposable, $n .\ldots n-5$. But, in
this case, $\{n,n-2\}$ is a non trivial interval of $G$; impossible.
If $n-2 \not\sim \{n-5,n\}$, we prove that $G-\{n-3,n-4\}$ is
indecomposable; contradiction.

  $\hspace*{0.2cm} \checkmark \, u \in X \setminus \{2,n-5\}$. It is clear that $G-\{n-3,n-4\}$ is
indecomposable which is impossible.

  $\hspace*{0.2cm} \checkmark \, u=n-3$. If $n-1 \not\sim \{n-3,n\}$, then as $G-\{n-3,n-4\}$ is decomposable,
we get $n \rule[0.4mm]{0.6cm}{0.1mm}\,\, n-2$. Besides, as
$G-\{n-1,n-4\}$ is decomposable, then $n-3 \sim \{n,n-1\}$ and so
$\{n,n-1\}$ is a non trivial interval of $G$; contradiction. If $n-2
\not\sim \{n-3,n\}$, then $n \rule[0.4mm]{0.6cm}{0.1mm}\,\, n-2$ and
as previously seen, we verify that $\{n,n-1\}$ is an interval of
$G$; impossible.

  $\hspace*{0.2cm} \checkmark \, u=n-4$. If $n-1 \not\sim \{n-4,n\}$, then $n .\ldots n-1$. As
$G-\{n-3,n-2\}$ is decomposable, we have $n .\ldots n-4$. In this
case, we prove that $G-\{2,n-2\}$ is indecomposable; contradiction.
If $n-2 \not\sim \{n-4,n\}$, then $n \rule[0.4mm]{0.6cm}{0.1mm}\,\,
n-2$. In this case, we prove that $G-\{n-1,n-4\}$ is indecomposable
which is impossible.

\end{description}
\end{itemize}

\item [Case $1.2$ :] $n-1 \in X^{+}$, $n-2 \in X^{-}(n-5)$, we get $\{n-2,n-4\}  \subseteq  X^{-}(n-5)$ and $\{n-1,n-3\}  \subseteq
X^{+}$. As previously seen, we have to distinguish the following
cases.
\begin{description}
\item [$\hookrightarrow$] If $n .\ldots Y$. It is easy to verify that $G-\{n-3,n-4\}$ is
indecomposable; contradiction. If $n
\rule[0.4mm]{0.6cm}{0.1mm}\,\,Y$ and $n
\rule[0.4mm]{0.6cm}{0.1mm}\,\,n-1$, then $n .\ldots n-2$. So, $G$ is isomorphic to one of the graphs defined in Proposition \ref{propqmoins5}. If $n \rule[0.4mm]{0.6cm}{0.1mm}\,\,Y$ and $n .\ldots n-1$, it suffices to
prove that $G-\{n-3,n-2\}$ is indecomposable; contradiction.

\item [$\hookrightarrow$] $n \in Y(u), u \in Y$. We distinguish the following cases.

  $\hspace*{0.2cm} \checkmark \, u \in X$. We demonstrate that $G-\{n-3,n-4\}$ is
indecomposable; impossible.

  $\hspace*{0.2cm} \checkmark \, u=n-3$. If $n-1 \not\sim \{n-3,n\}$, then $n .\ldots n-1$. As $G-\{n-3,n-4\}$ is
decomposable, then $n .\ldots n-2$. Besides, $G-\{n-4,n-1\}$ is
decomposable implies that $n \rule[0.4mm]{0.6cm}{0.1mm}\,\, n-3$.
But in this case, $\{n,n-1\}$ is a non trivial interval of $G$;
impossible. If $n-2 \not\sim \{n-3,n\}$, then , $n-2 .\ldots n$ and
we may assume that $n \rule[0.4mm]{0.6cm}{0.1mm}\,\, n-1$. Since
$G-\{n-1,n-4\}$ is decomposable, then $n
\rule[0.4mm]{0.6cm}{0.1mm}\,\, n-3$ and so $\{n,n-1\}$ is an
interval of $G$; impossible.

  $\hspace*{0.2cm} \checkmark \, u=n-4$. If $n-1 \not\sim \{n-4,n\}$, then as $G-\{n-1,n-4\}$ and $G-\{n-2,n-3\}$
are decomposable, we have $n .\ldots \{n-2,n-4\}$. Therefore, $G$ is isomorphic to one of the graphs defined in Proposition \ref{propqmoins5}. If $n-2
\not\sim \{n-4,n\}$, then $n \rule[0.4mm]{0.6cm}{0.1mm}\,\,n-2$ and
we obtain that $G-\{n-1,n-4\}$ is indecomposable which is
impossible.\\

\end{description}

\item [Case $1.3$ :] $n-1 \in X^{-}$, $n-2 \in X(\mu)$ where $\mu \in \{n-6,n-5\}$.
\begin{itemize}
  \item If $\mu =n-6$, then $\{n-2,n-4\}  \subseteq  X(n-6)$ and $\{n-1,n-3\}  \subseteq  X^{-}$.
% then by interchanging $n-4$ and $n-6$, we can assume that $H$ verifies
%$\{n-2,n-4\}  \subseteq  X^{-}(n-5), \,\{n-3,n-1\} \subseteq
%X(n-6)$, $H[\{n-3,n-1,n-6\}]$ is complete or empty and $E(H_{X})=
%\{\{n-1,n-4\},\{n-4,n-3\},\{n-3,n-2\}\}$. \\
Similarly, we shall examine the two following cases.

\begin{description}
\item [$\hookrightarrow$] $ n \in \langle Y \rangle$. If $n \rule[0.4mm]{0.6cm}{0.1mm}\,\, Y$, we prove that $G-\{n-3,n-4\}$ is indecomposable;
contradiction. If $n .\ldots Y$, then as $G-\{n-3,n-4\}$ is
decomposable we obtain that $n \rule[0.4mm]{0.6cm}{0.1mm}\,\, n-2$.
Moreover, we show that $G-\{n-2,n-3\}$ is indecomposable if $n
\rule[0.4mm]{0.6cm}{0.1mm}\,\, n-1$ and $G-\{n-6,n-1\}$ is
indecomposable if $n .\ldots\,\, n-1$ which is impossible.

\item [$\hookrightarrow$] $n \in Y(u)$ where $u \in Y$. We distinguish the following cases.

  $\hspace*{0.2cm} \checkmark \, u = n-5$. We have either $n-1 \rule[0.4mm]{0.6cm}{0.1mm}\,\, n$  or $n-2 .\ldots n$. In both cases, we prove that
$G-\{n-4,n-3\}$ is indecomposable; impossible.

  $\hspace*{0.2cm} \checkmark \, u=n-6$. If $n-1 \not\sim \{n-6,n\}$, then as $G-\{n-3,n-2\}$ is
decomposable, $n-6 \sim \{n,n-4\}$ that is, $n-6 \sim \{n,n-2\}$. It
follows that $\{n,n-2\}$ is a non trivial interval of $G$;
impossible. If $n-2 \not\sim \{n-6,n\}$, then we can assume that
$n-1 .\ldots n$ which implies that $G-\{n-4,n-3\}$ is
indecomposable; contradiction.

  $\hspace*{0.2cm} \checkmark \, u \in X \setminus \{n-5,n-6\}$. We
verify that $G-\{n-4,n-3\}$ is indecomposable which is impossible.

  $\hspace*{0.2cm} \checkmark \, u = n-3$. If $n-1 \not\sim
\{n-3,n\}$, then $n \rule[0.4mm]{0.6cm}{0.1mm}\,\, n-1$. Since
$G-\{n-3,n-4\}$ is decomposable, $n-2 \rule[0.4mm]{0.6cm}{0.1mm}\,\,
n$. Moreover, $G-\{n-1,n-4\}$ is decomposable implies $n-3 .\ldots
n$. Consequently, $\{n,n-1\}$ would be a non trivial interval of
$G$; contradiction. If $n-2 \not\sim \{n-3,n\}$, then $n
\rule[0.4mm]{0.6cm}{0.1mm}\,\,n-2$ and we verify in the same way,
that $\{n,n-1\}$ is a non trivial interval of $G$; contradiction.

  $\hspace*{0.2cm} \checkmark \, u = n-4$. If $n-1 \not\sim \{n-4,n\}$, then $n-1 .\ldots n$ and $G-\{n-4,n-6\}$ is
indecomposable; impossible. If $n-2 \not\sim \{n-4,n\}$, then
$G-\{n-4,n-1\}$ is indecomposable; contradiction.
\end{description}

  \item If $\mu =n-5$, then $\{n-2,n-4\}  \subseteq  X^{-}(n-5)$ and $\{n-1,n-3\}  \subseteq  X^{-}$. We distinguish the two cases.
\begin{description}
\item [$\hookrightarrow$] $ n \in \langle Y \rangle$. If $n .\ldots Y$ and $n .\ldots n-1$ (resp. $n \rule[0.4mm]{0.6cm}{0.1mm}\,\,n-1$),
then $n \rule[0.4mm]{0.6cm}{0.1mm}\,\,n-2$ and $G$ is one of the graphs defined in Proposition \ref{propqmoins5} (resp. we prove that $G-\{n-3,n-2\}$ is
indecomposable which is impossible). If $n
\rule[0.4mm]{0.6cm}{0.1mm}\,\,Y$, we verify that $G-\{n-3,n-4\}$ is
indecomposable; contradiction.

\item [$\hookrightarrow$] $n \in Y(u)$. Let us distinguish the
following cases.

  $\hspace*{0.2cm} \checkmark \, u \in X \setminus \{n-5\}$. It is easy to verify that
$G-\{n-3,n-4\}$ is indecomposable; impossible.

  $\hspace*{0.2cm} \checkmark \, u=n-5$. If $n-1 \not\sim \{n-5,n\}$, then $n .\ldots n-5$ (because $G-\{n-3,n-4\}$ is decomposable). We prove that $\{n,n-2\}$ is a non trivial interval
of $G$; impossible. If $n-2 \not\sim \{n-5,n\}$, then
$G-\{n-3,n-4\}$ is indecomposable; impossible.

  $\hspace*{0.2cm} \checkmark \, u=n-3$. If $n-1 \not\sim \{n-3,n\}$, then $n
\rule[0.4mm]{0.6cm}{0.1mm}\,\,n-1$. As $G-\{n-3,n-4\}$ is
decomposable, then $n \rule[0.4mm]{0.6cm}{0.1mm}\,\,n-2$. Besides,
$G-\{n-4,n-1\}$ is decomposable implies that $n .\ldots n-3$. But in
this case, $\{n,n-1\}$ is a non trivial interval of $G$;
contradiction. If $n-2 \not\sim \{n-3,n\}$, then $n
\rule[0.4mm]{0.6cm}{0.1mm}\,\,n-2$ and we may assume that $n .\ldots
n-1$. As $G-\{n-1,n-4\}$ is decomposable, then necessarily $n
.\ldots n-3$ and $\{n,n-1\}$ is a non trivial interval of $G$;
impossible.

  $\hspace*{0.2cm} \checkmark \, u=n-4$. If $n-1 \not\sim \{n-4,n\}$, then $n .\ldots n-1$.
As $G-\{n-4,n-1\}$ is decomposable, then $n .\ldots n-2$. Besides,
$G-\{n-3,n-2\}$ is decomposable implies that $n .\ldots n-4$.
Therefore, $n .\ldots \{n-1,n-2,n-4\}$ and $G$ is isomorphic to one of the graphs defined in Proposition \ref{propqmoins5}. If $n-2 \not\sim \{n-4,n\}$, we
obtain $G-\{n-1,n-4\}$ is indecomposable which is impossible.

\end{description}
\end{itemize}

\item \underline{Case 2: } $E(H_{X}) =
\{\{n-2,n-1\}, \{n-4,n-3\}\}$. We proceed as previously and we apply
Theorem \ref{ThBDInonconnexe} and Lemma \ref{lemqmoins3}. Besides,
up to isomorphism, we can examine the three cases.
\begin{description}
  \item [Case $2.1$ :]  $n-1 \in X(\mu)$ where $\mu \in \{1,2\}$ and $n-2 \in
X^{-}(n-5)$. We have to distinguish the following cases.
\begin{itemize}
  \item If $n-3 \in X(\nu)$ where $\nu \in \{1,2\}$ and $n-4 \in X^{-}(n-5)$, then we discuss the three cases.
   \begin{itemize}
    \item If $\mu=\nu=1$, we have $\{n-2,n-4\}  \subseteq  X^{-}(n-5)$ and $\{n-1,n-3\}  \subseteq  X(1)$. By theorem 2.7, $G[\{n-1,n-3,1\}]$
is empty or complete. As seen in Case $1$, we distinguish the two
cases.
\begin{description}
\item [$\hookrightarrow $] $\, n \in \langle Y \rangle$. We simply verify that $G-\{n-3,n-4\}$ is indecomposable; contradiction.

\item [$\hookrightarrow $] $\, n \in Y(u), u \in Y$. We examine the following cases.

  $\hspace*{0.2cm} \checkmark \, u=n-5$. If $n-1 \not\sim \{n-5,n\}$, then $G-\{n-3,n-4\}$ is
decomposable gives that $n .\ldots n-5$. Consequently, $\{n,n-2\}$
is a non trivial interval of $G$; impossible. If $n-2 \not\sim
\{n-5,n\}$, then $n \rule[0.4mm]{0.6cm}{0.1mm}\,\,n-2$ and we may
assume that $n .\ldots \,n-1$. We demonstrate that $G-\{n-3,n-4\}$
is indecomposable; contradiction.

 $\hspace*{0.2cm} \checkmark \, u=1$. If $n-1 \not\sim \{1,n\}$, necessarily $n
\rule[0.4mm]{0.6cm}{0.1mm}\,\, n-2$ and $1 \sim \{n-1,n\}$
(otherwise $G-\{n-3,n-4\}$ is indecomposable; contradiction). But in
this case, $\{n,n-1\}$ is a non trivial interval of $G$; impossible.
If $n-2 \not\sim \{1,n\}$, then $n \rule[0.4mm]{0.6cm}{0.1mm}\,\,
n-2$ and we can assume that $n-1 \sim \{1,n\}$. Moreover, as
$G-\{n-3,n-4\}$ is decomposable, then $1 \sim \{n-1,n\}$. Thus,
$\{n,n-1\}$ is a non trivial interval of $G$; impossible.

 $\hspace*{0.2cm} \checkmark \, u \in X \setminus \{1,n-5\}$. Clearly,
$G-\{n-3,n-4\}$ is indecomposable; contradiction.

 $\hspace*{0.2cm} \checkmark \, u=n-3$. If $n-1 \not\sim \{n-3,n\}$, then $n-1 \not\sim
\{n,1\}$. We prove that $G-\{n-2,n-3\}$ is indecomposable;
contradiction. If $n-2 \not\sim \{n-3,n\}$, then , $n
\rule[0.4mm]{0.6cm}{0.1mm}\,\,n-2$ and we can assume that $n-1 \sim
\{n,n-3\}$. We get as previously, $n-3 \sim \{n,n-1\}$ and thus, $G$ is isomorphic to one of the graphs defined in Proposition \ref{propqmoins5}.

 $\hspace*{0.2cm} \checkmark \, u=n-4$. If $n-1 \not\sim \{n-4,n\}$, then $n-1
\rule[0.4mm]{0.6cm}{0.1mm}\,\,n$. As $G-\{n-1,n-4\}$ is
decomposable, we have $n .\ldots n-2$. Besides, $G-\{n-2,n-3\}$ is
decomposable implies that $n .\ldots n-4$. Hence, $n .\ldots
\{n-2,n-4\}$ and $G$ is isomorphic to one of the graphs defined in Proposition \ref{propqmoins5}. If $n-2 \not\sim \{n-4,n\}$, then $n
\rule[0.4mm]{0.6cm}{0.1mm}\,\, n-2$ and we may assume that $n-1 \sim
\{n-4,n\}$. We obtain that $G-\{n-1,n-4\}$ is indecomposable;
contradiction.\\
\end{description}

   \item If $\mu=\nu=2$, we get $\{n-2,n-4\}  \subseteq  X^{-}(n-5)$ and $\{n-1,n-3\}  \subseteq
X(2)$, then we prove that $2$ is a non critical vertex of $G$ which
is impossible.\\

   \item If $\mu =2$ and $\nu=1$, then $\{n-2,n-4\}  \subseteq  X^{-}(n-5)$, $n-3 \in X(1)$ and
$n-1 \in  X(2)$, then necessarily $1 .\ldots n-3$ (otherwise we
prove that $2$ is a non critical vertex of $G$; contradiction).
Thus, $n-3 \in X^{-}(1)$. Let us distinguish the two cases.
\begin{description}
\item [$\hookrightarrow $]$n \in \langle Y \rangle$, it suffices to prove that $G-\{n-4,n-3\}$ is indecomposable which
is impossible.

\item [$\hookrightarrow $]$n \in Y(u), u \in Y$. We distinguish the following
cases.

$\hspace*{0.2cm} \checkmark \, u = n-5$. We have $n .\ldots
\{n-3,n-4\}$. If $n-1 \not\sim \{n-5,n\}$, then $n-1
\rule[0.4mm]{0.6cm}{0.1mm}\,\,n$. Since $G-\{n-4,n-3\}$ is
decomposable, then $n .\ldots n-5$  and so $\{n,n-2\}$
 is a non trivial interval of $G$; impossible. If $n-2 \not\sim \{n-5,n\}$, then $n \rule[0.4mm]{0.6cm}{0.1mm}\,\,
n-2$ and we may assume that $n .\ldots n-1$. So, $G-\{n-4,n-3\}$ is
indecomposable; impossible.

$\hspace*{0.2cm} \checkmark \, u=2$. If $n-1 \not\sim \{2,n\}$, we
prove that $n \rule[0.4mm]{0.6cm}{0.1mm}\,\,n-2$ (otherwise
$G-\{n-3,n-4\}$ is indecomposable; contradiction). In this case, we
demonstrate that $\{n,n-1\}$ is a non trivial interval of $G$ (if $2
\sim \{n,n-1\}$) and that $G-\{n-4,n-3\}$ is indecomposable (if $2
\not\sim \{n,n-1\}$); contradiction. If $n-2 \not\sim \{2,n\}$, then
$n \rule[0.4mm]{0.6cm}{0.1mm}\,\,n-2$ and we may assume that $n-1
\sim \{2,n\}$. In this case, we demonstrate that $\{n,n-1\}$ is a
non trivial interval of $G$ (if $2 \sim \{n,n-1\}$) and that
$G-\{n-3,n-4\}$ is indecomposable (if $2 \not\sim \{n,n-1\}$);
contradiction.

$\hspace*{0.2cm} \checkmark \, u \in X \setminus \{2,n-5\}$. We get
$G-\{n-3,n-4\}$ is indecomposable; contradiction.

$\hspace*{0.2cm} \checkmark \, u=n-3$. We prove that $G-\{n-3,n-4\}$
is indecomposable; contradiction.

$\hspace*{0.2cm} \checkmark \, u= n-4$. If $n-1 \not\sim \{n-4,n\}$,
then $n \rule[0.4mm]{0.6cm}{0.1mm}\,\,n-1$. We prove that
$G-\{n-2,n-4\}$ is indecomposable; impossible. If $n-2 \not\sim
\{n-4,n\}$, then $n \rule[0.4mm]{0.6cm}{0.1mm}\,\,n-2$ and we may
assume that $n-1 .\ldots n$. So, $G-\{n-4,n-3\}$ is indecomposable;
contradiction.
\end{description}
   \end{itemize}

\item If $n-3 \in X^{+}$ and $n-4 \in X^{-}(n-5)$, we have $\{n-2,n-4\}  \subseteq  X^{-}(n-5)$, $n-1 \in  X(\mu)$ where $\mu \in \{1,2\}$ and
$n-3 \in X^{+}$. Let us distinguish the two cases.
\begin{description}
\item [$\hookrightarrow $]$n \in \langle Y \rangle$. We simply verify that $G-\{n-3,n-4\}$ is indecomposable; contradiction.

\item [$\hookrightarrow $]$n \in Y(u), u \in Y$. We examine the following cases.

$\hspace*{0.2cm} \checkmark \, u=n-5 \,(resp. u=\mu)$. If $n-1
\not\sim \{n,u\}$, then as $G-\{n-3,n-4\}$ is decomposable, we get
$n .\ldots n-5$ (resp. $n \rule[0.4mm]{0.6cm}{0.1mm}\,\, n-2$ and $1
\sim \{n-1,n\}$). Consequently, $\{n,n-2\}$ (resp. $\{n,n-1\}$) is a
non trivial interval of $G$; impossible. If $n-2 \not\sim \{n,u\}$,
necessarily, $n \rule[0.4mm]{0.6cm}{0.1mm}\,\,n-2$ and we may assume
that $n-1 \sim \{n,u\}$. We demonstrate that $G-\{n-3,n-4\}$ is
indecomposable ($\{n,n-1\}$ is a non trivial interval of $G$);
contradiction.

$\hspace*{0.2cm} \checkmark \, u \in X \setminus \{\mu,n-5\}$. We
verify that $G-\{n-3,n-4\}$ is indecomposable; impossible.

$\hspace*{0.2cm} \checkmark \, u=n-3$. If $n-1 \not\sim \{n-3,n\}$,
we demonstrate that $G-\{n-2,n-3\}$ is indecomposable; impossible.
If $n-2 \not\sim \{n-3,n\}$, we verify that $G-\{n-3,n-4\}$ is
indecomposable; impossible.

$\hspace*{0.2cm} \checkmark \, u=n-4$. If $n-1 \not\sim \{n-4,n\}$,
then $n-1 \rule[0.4mm]{0.6cm}{0.1mm}\,\,n$ and we prove that
$G-\{n-2,n-4\}$ is indecomposable which is impossible. If $n-2
\not\sim \{n-4,n\}$, then $n \rule[0.4mm]{0.6cm}{0.1mm}\,\,n-2$. We
may assume that $n .\ldots n-1$ and we obtain that $G-\{n-3,n-4\}$
is indecomposable; contradiction.
\end{description}

 \item If $n-3 \in X^{-}$ and $n-4 \in X(\nu)$ where $\nu \in
\{n-6,n-5\}$. We have to distinguish the following cases.
\begin{itemize}
  \item If $\mu \in \{1,2\}$ and $\nu=n-6$, then we have $n-2 \in X^{-}(n-5)$, $n-4 \in X(n-6)$,
$n-3 \in X^{-}$ and $n-1 \in X(1)$ (resp. $n-1 \in X(2)$). Let us
distinguish the two cases.
\begin{description}
\item [$\hookrightarrow $]$n \in \langle Y \rangle$. As previously seen, we prove that $G-\{n-3,n-4\}$ is indecomposable; impossible.

\item [$\hookrightarrow $]$n \in Y(u), u \in Y$. We distinguish the following cases.

$\hspace*{0.2cm} \checkmark \, u = n-5$. If $n-1 \not\sim
\{n-5,n\}$, then $n \rule[0.4mm]{0.6cm}{0.1mm}\,\, n-1$. As
$G-\{n-6,n-3\}$ is decomposable, we obtain $n .\ldots n-5$. So,
$\{n,n-2\}$ is a non trivial interval of $G$; impossible. If $n-2
\not\sim \{n-5,n\}$, then $n \rule[0.4mm]{0.6cm}{0.1mm}\,\, n-2$ and
we may assume that $n-1 .\ldots n$. It follows that $G-\{n-3,n-4\}$
is indecomposable; impossible.

$\hspace*{0.2cm} \checkmark \, u=n-6$. In this case, we have either
$n-1 .\ldots n$ or $n-2 .\ldots n$. Consequently, $G-\{n-3,n-4\}$ is
indecomposable; impossible.

$\hspace*{0.2cm} \checkmark \, u=1 \,$(resp $u=2$). If $n-1 \not\sim
\{n,u\}$, then as $G-\{n-3,n-4\}$ is decomposable, $u \sim
\{n-1,n\}$. Besides, $G-\{n-6,n-3\}$ is decomposable implies that $n
\rule[0.4mm]{0.6cm}{0.1mm}\,\, n-2$. Therefore, $\{n,n-1\}$ would be
an interval of $G$; contradiction. If $n-2 \not\sim \{n,u\}$, then
$n \rule[0.4mm]{0.6cm}{0.1mm}\,\, n-2$ and we may assume that $n-1
\sim \{n,u\}$. As previously, as $G-\{n-3,n-4\}$ is decomposable
gives that $u \sim \{n-1,n\}$ and so $\{n,n-1\}$ is an interval of
$G$; impossible.

$\hspace*{0.2cm} \checkmark \, u \in X \setminus \{n-5,n-6,1\}$
(resp. $u \in X \setminus \{n-5,n-6,2\}$). It is easy to prove that
$G-\{n-3,n-4\}$ is indecomposable which is impossible.

$\hspace*{0.2cm} \checkmark \, u = n-3$. In this case, $n-1
\rule[0.4mm]{0.6cm}{0.1mm}\,\, n$ or $n-2
\rule[0.4mm]{0.6cm}{0.1mm}\,\, n$ and $G-\{n-3,n-4\}$ is
indecomposable; contradiction.

$\hspace*{0.2cm} \checkmark \, u = n-4$. We have either $n .\ldots
n-1$ or $n .\ldots n-2$ and $G-\{n-3,n-4\}$ is indecomposable;
impossible.
\end{description}

  \item If $\mu \in \{1,2\}$ and $\nu=n-5$, we get $\{n-2,n-4\} \subseteq  X^{-}(n-5)$, $n-1 \in  X(\mu)$ where $\mu \in \{1,2\}$ and
$n-3 \in  X^{-}$. Let us distinguish the two cases.
\begin{description}
\item [$\hookrightarrow $]$n \in \langle Y \rangle$. We simply verify that $G-\{n-3,n-4\}$ is indecomposable; contradiction.

\item [$\hookrightarrow $]$n \in Y(u), u \in Y$. We distinguish the following cases.

$\hspace*{0.2cm} \checkmark \, u=n-5 \,(resp. u=\mu)$. If $n-1
\not\sim \{n,u\}$, then as $G-\{n-3,n-4\}$ is decomposable, we get
$n .\ldots n-5$ (resp. $n \rule[0.4mm]{0.6cm}{0.1mm}\,\, n-2$ and $1
\sim \{n-1,n\}$). Consequently, $\{n,n-2\}$ (resp. $\{n,n-1\}$) is a
non trivial interval of $G$; impossible. If $n-2 \not\sim \{n,u\}$,
necessarily, $n \rule[0.4mm]{0.6cm}{0.1mm}\,\,n-2$ and we may assume
that $n-1 \sim \{n,u\}$. We demonstrate that $G-\{n-3,n-4\}$ is
indecomposable ($\{n,n-1\}$ is a non trivial interval of $G$);
contradiction.

$\hspace*{0.2cm} \checkmark \, u \in X \setminus \{\mu,n-5\}$, we
verify that $G-\{n-3,n-4\}$ is indecomposable; impossible.

$\hspace*{0.2cm} \checkmark \, u=n-3$. If $n-1 \not\sim \{n-3,n\}$,
we demonstrate that $G-\{n-2,n-3\}$ is indecomposable; impossible.
If $n-2 \not\sim \{n-3,n\}$, we verify that $G-\{n-3,n-4\}$ is
indecomposable; impossible.

$\hspace*{0.2cm} \checkmark \, u=n-4$. If $n-1 \not\sim \{n-4,n\}$,
then $n-1 \rule[0.4mm]{0.6cm}{0.1mm}\,\,n$ and we prove that
$G-\{n-2,n-4\}$ is indecomposable which is impossible. If $n-2
\not\sim \{n-4,n\}$, then $n \rule[0.4mm]{0.6cm}{0.1mm}\,\,n-2$. We
may assume that $n .\ldots n-1$ and we obtain that $G-\{n-3,n-4\}$
is indecomposable; contradiction.
\end{description}
\end{itemize}

\item If $n-3 \in X(2)$ and $n-4 \in X^{-}(1)$, then $\mu \in \{1,2\}$. In both cases, we permute $n-3$ and $2$ and we may return to one of the previous
cases.

\end{itemize}

  \item [Case $2.2$ :] $n-1 \in X^{+}$ and $n-2 \in X^{-}(n-5)$. By
Theorem \ref{ThBDInonconnexe}, we have to examine the three cases.
\begin{itemize}
  \item If $n-3 \in X(\nu)$ where $\nu \in \{1,2\}$ and $n-4 \in X^{-}(n-5)$, then it is clear that we may return to one of the previous cases.
  \item If $n-3 \in X^{+}$ and $n-4 \in X^{-}(n-5)$, that is, $\{n-2,n-4\}  \subseteq  X^{-}(n-5)$ and $\{n-1,n-3\}  \subseteq  X^{+}$. Let us distinguish
 the following cases.
\begin{description}
\item [$\hookrightarrow $] If $n .\ldots Y$. It is easy to verify that $G-\{n-3,n-4\}$ is indecomposable; contradiction.
If $n \rule[0.4mm]{0.6cm}{0.1mm}\,\,Y$ and $n
\rule[0.4mm]{0.6cm}{0.1mm}\,\,n-1$, then $n .\ldots n-2$. So,
$\{n,n-1\}$ is a non trivial interval of $G$ which is impossible. If
$n \rule[0.4mm]{0.6cm}{0.1mm}\,\,Y$ and $n .\ldots n-1$, it suffices
to verify that $\{n,n-1\}$ is a non trivial interval of $G$ (if $n
.\ldots n-2$) and that $G-\{n-3,n-4\}$ is indecomposable (if $n
\rule[0.4mm]{0.6cm}{0.1mm}\,\,n-2$); contradiction.

\item [$\hookrightarrow $]$n \in Y(u), u \in Y$. We distinguish the following cases.

$\hspace*{0.2cm} \checkmark \, u \in X$. We demonstrate that
$G-\{n-3,n-4\}$ is indecomposable; impossible.

$\hspace*{0.2cm} \checkmark \, u=n-3$. If $n-1 \not\sim \{n-3,n\}$,
then $n .\ldots n-1$. As $G-\{n-3,n-4\}$ is decomposable, then $n
.\ldots n-2$. Besides, $G-\{n-4,n-1\}$ is decomposable implies that
$n \rule[0.4mm]{0.6cm}{0.1mm}\,\, n-3$. But in this case,
$G-\{n-2,n-3\}$ is indecomposable; impossible. If $n-2 \not\sim
\{n-3,n\}$, then , $n-2 .\ldots n$ and we may assume that $n
\rule[0.4mm]{0.6cm}{0.1mm}\,\, n-1$. Since $G-\{n-1,n-4\}$ is
decomposable, then $n \rule[0.4mm]{0.6cm}{0.1mm}\,\, n-3$ and so $G$ is isomorphic to one of the graphs defined in Proposition \ref{propqmoins5}.

$\hspace*{0.2cm} \checkmark \, u=n-4$. If $n-1 \not\sim \{n-4,n\}$,
then as $G-\{n-1,n-4\}$ and $G-\{n-2,n-3\}$ are decomposable, we
have $n .\ldots \{n-2,n-4\}$. Therefore, $G$ is isomorphic to one of the graphs defined in Proposition \ref{propqmoins5}. If $n-2 \not\sim \{n-4,n\}$, then $n
\rule[0.4mm]{0.6cm}{0.1mm}\,\,n-2$ and we obtain that
$G-\{n-1,n-4\}$ is indecomposable which is impossible.
\end{description}

 \item If $n-3 \in X(2)$ and $n-4 \in X^{-}(1)$, then it suffices
to permute $n-3$ and $2$ and we may return to one of the previous
cases.
\end{itemize}

\item [Case $2.3$ :] $n-1 \in X^{-}$ and $n-2 \in X(\mu)$ where $\mu \in
\{n-6,n-5\}$. Similarly to Case $2.2$ and using Theorem
\ref{ThBDInonconnexe}, we have to distinguish the following cases.
\begin{itemize}
  \item If $n-3 \in X(\nu)$ where $\nu \in \{1,2\}$ and $n-4 \in X^{-}(n-5)$,
we have either $\mu=n-6$ and $\nu \in \{1,2\}$ or $\mu=n-5$ and $\nu
\in \{1,2\}$. It is clear that, in both cases, we may return to one
of the cases treated above.

 \item If $n-3 \in X^{-}$ and $n-4 \in X(\nu)$ where $\nu \in
\{n-6,n-5\}$, then we examine the following cases.
\begin{itemize}
  \item If $\mu=\nu=n-6$, we have $\{n-2,n-4\}  \subseteq  X(n-6)$ and $\{n-1,n-3\}  \subseteq  X^{-}$ and we verify that $n-6$ is a non critical
vertex of $G$; contradiction.\\
  \item If $\mu=\nu=n-5$, we get $\{n-2,n-4\}  \subseteq  X^{-}(n-5)$ and $\{n-1,n-3\}  \subseteq
X^{-}$. We distinguish the following cases.
\begin{description}
\item [$\hookrightarrow $] If $n .\ldots Y$. Suppose that $n .\ldots n-1$ (resp. $n \rule[0.4mm]{0.6cm}{0.1mm}\,\,n-1$). So, $n
\rule[0.4mm]{0.6cm}{0.1mm}\,\,n-2$ because otherwise $G$ is
decomposable (resp. because otherwise $G-\{n-3,n-4\}$ is
indecomposable); impossible. But, $\{n,n-1\}$ is a non trivial
interval of $G$ which is impossible. If $n
\rule[0.4mm]{0.6cm}{0.1mm}\,\,Y$. In this case, we verify that
$G-\{n-3,n-4\}$ is indecomposable; contradiction.

\item [$\hookrightarrow $]$n \in Y(u), \, u \in Y$. Let us distinguish the following cases.

$\hspace*{0.2cm} \checkmark \, u \in X$. We proceed in the same
manner as previously.

$\hspace*{0.2cm} \checkmark \, u=n-3$. If $n-1 \not\sim \{n-3,n\}$,
then $n \rule[0.4mm]{0.6cm}{0.1mm}\,\,n-1$. We prove that
$G-\{n-3,n-2\}$ is indecomposable which is impossible. If $n-2
\not\sim \{n-3,n\}$, then , as previously, we have $\{n-1,n-3\}
.\ldots n \rule[0.4mm]{0.6cm}{0.1mm}\,\,n-2$ and $G$ is isomorphic to one of the graphs defined in Proposition \ref{propqmoins5}.

$\hspace*{0.2cm} \checkmark \, u=n-4$. If $n-1 \not\sim \{n-4,n\}$,
then $n \rule[0.4mm]{0.6cm}{0.1mm}\,\, n-1$. Similarly to the
previous case, $n .\ldots \{n-4,n-2\}$ and $G$ is isomorphic to one of the graphs defined in Proposition \ref{propqmoins5}. If $n-2 \not\sim \{n-4,n\}$, then
$n \rule[0.4mm]{0.6cm}{0.1mm}\,\,n-2$ and we obtain as previously
that $G-\{n-1,n-4\}$ is indecomposable which is impossible.
\end{description}

 \item If $\mu=n-6$ and $\nu=n-5$, we have $n-4 \in X^{-}(n-5)$, $n-2 \in X(n-6)$ and $\{n-1,n-3\}  \subseteq  X^{-}$.
Let us distinguish the two cases.%ici N2
\begin{description}
\item [$\hookrightarrow $]$n \in \langle Y \rangle$. We proceed as in the last case and we verify that $G-\{n-3,n-4\}$
is indecomposable which is impossible.

\item [$\hookrightarrow $]$n \in Y(u), u \in Y$. We examine the following cases.

$\hspace*{0.2cm} \checkmark \, u=n-6$. If $n-1 \not\sim \{n-6,n\}$,
then $n \rule[0.4mm]{0.6cm}{0.1mm}\,\, n-1$. As $G-\{n-1,n-4\}$ is
decomposable then $n-6 \sim \{n,n-2\}$. But in this case,
$\{n,n-2\}$ would be a non trivial interval of $G$; impossible. If
$n-2 \not\sim \{n-6,n\}$, we verify that $G-\{n-3,n-4\}$ is
indecomposable; contradiction.

$\hspace*{0.2cm} \checkmark \, u \in \{n-5,n-4\}$. We have either $n
\rule[0.4mm]{0.6cm}{0.1mm}\,\, n-1$ or $n .\ldots n-2$. We verify
that $G-\{n-3,n-4\}$ is indecomposable; contradiction.

$\hspace*{0.2cm} \checkmark \, u \in X \setminus \{n-6,n-5\}$. We
verify that $G-\{n-3,n-4\}$ is indecomposable which is impossible.

$\hspace*{0.2cm} \checkmark \, u=n-3$. If $n-1 \not\sim \{n-3,n\}$,
we prove that $G-\{n-2,n-3\}$ is indecomposable; impossible. If $n-2
\not\sim \{n-3,n\}$, then $n \rule[0.4mm]{0.6cm}{0.1mm}\,\, n-2$. As
$G-\{n-1,n-4\}$ is decomposable, we get $n-3 .\ldots n$. It follows
that $\{n,n-1\}$ would be a non trivial interval of $G$; impossible.
\end{description}

\item If $n-3 \in X(2)$ and $n-4 \in X^{-}(1)$, then it suffices
to permute $n-3$ and $2$ and we may return to one of the previous
cases.
\end{itemize}

\end{itemize}

\end{description}
\end{description}
\end{itemize}
{\hspace*{\fill}$\Box$\medskip}\\

\end{document}